\documentclass[10pt]{amsart}

\usepackage{amsmath,amsthm,amsfonts,amssymb,graphicx,epstopdf,kbordermatrix,scalerel,tikz}
\usepackage[all]{xy}
\usepackage[hidelinks]{hyperref}
\usepackage[margin=1in]{geometry}

\newtheorem{theorem}{Theorem}[subsection]
\newtheorem{lemma}[theorem]{Lemma}
\newtheorem{proposition}[theorem]{Proposition}
\newtheorem{corollary}[theorem]{Corollary}
\newtheorem*{claim*}{Claim}

\theoremstyle{definition} \newtheorem{definition}[theorem]{Definition}
\newtheorem{example}[theorem]{Example}

\theoremstyle{remark} 
\newtheorem{remark}[theorem]{Remark}
\newtheorem{warning}[theorem]{Warning}

\newcommand{\C}{\mathbb{C}}
\newcommand{\F}{\mathbb{F}}
\newcommand{\R}{\mathbb{R}}
\newcommand{\Z}{\mathbb{Z}}

\newcommand{\Ib}{\mathbf{I}}
\newcommand{\Ibrl}{\mathbf{I}_{r/l}}
\newcommand{\x}{\mathbf{x}}
\newcommand{\y}{\mathbf{y}}

\newcommand{\A}{\mathcal{A}}
\newcommand{\B}{\mathcal{B}} 
\newcommand{\Cc}{\mathcal{C}}
\newcommand{\Ccr}{\mathcal{C}_r}
\newcommand{\Ccl}{\mathcal{C}_l}
\newcommand{\Ccrl}{\mathcal{C}_{r/l}}
\newcommand{\D}{\mathcal{D}}
\newcommand{\G}{\mathcal{G}}
\newcommand{\M}{\mathcal{M}}
\newcommand{\Nc}{\mathcal{N}}
\newcommand{\Pc}{\mathcal{P}}
\newcommand{\Sc}{\mathcal{S}}
\newcommand{\T}{\mathcal{T}}
\newcommand{\W}{\mathcal{W}}

\newcommand{\CFK}{\widehat{CFK}}
\newcommand{\HFK}{\widehat{HFK}}
\newcommand{\Uh}{\mathcal{U}_{\hslash}(\mathfrak{gl}(1|1))}
\newcommand{\Uq}{\mathcal{U}_q(\mathfrak{gl}(1|1))}
\newcommand{\UT}{\mathcal{U}_T(\mathfrak{sl}(1|1))}

\newcommand{\OrUUUU}{
\begin{tikzpicture}[scale=0.1]
\draw [->] (2,0) -- (5,3);
\draw (5,0) -- (4,1);
\draw [->] (3,2) -- (2,3);
\draw [->] (0,0) -- (0,3);
\draw [->] (7,0) -- (7,3);
\end{tikzpicture}
}

\newcommand{\OrUUDD}{
\begin{tikzpicture}[scale=0.1]
\draw [->] (2,0) -- (5,3);
\draw [->] (4,1) -- (5,0);
\draw (3,2) -- (2,3);
\draw [->] (0,0) -- (0,3);
\draw [->] (7,3) -- (7,0);
\end{tikzpicture}
}

\newcommand{\OrUU}{
\begin{tikzpicture}[scale=0.1]
\draw [->] (0,0) -- (3,3);
\draw (3,0) -- (2,1);
\draw [->] (1,2) -- (0,3);
\end{tikzpicture}
}

\newcommand{\OrDU}{
\begin{tikzpicture}[scale=0.1]
\draw [->] (3,3) -- (0,0);
\draw (3,0) -- (2,1);
\draw [->] (1,2) -- (0,3);
\end{tikzpicture}
}

\newcommand{\OrUD}{
\begin{tikzpicture}[scale=0.1]
\draw [->] (0,0) -- (3,3);
\draw [->] (2,1) -- (3,0);
\draw (1,2) -- (0,3);
\end{tikzpicture}
}

\newcommand{\OrDD}{
\begin{tikzpicture}[scale=0.1]
\draw [->] (3,3) -- (0,0);
\draw [->] (2,1) -- (3,0);
\draw (1,2) -- (0,3);
\end{tikzpicture}
}

\newcommand{\MinLR}{
\begin{tikzpicture}[scale=0.1]
\draw [->] (3,3) arc [radius=1.5, start angle=180, end angle=360];
\end{tikzpicture}
}

\newcommand{\MinRL}{
\begin{tikzpicture}[scale=0.1]
\draw [->] (3,3) to [out=-90, in=0] (1.5,1.5) to [out=180, in = -90] (0,3);
\end{tikzpicture}
}

\newcommand{\MaxLR}{
\begin{tikzpicture}[scale=0.1]
\draw [->] (0,0) to [out=90, in=180] (1.5,1.5) to [out=0, in=90] (3,0);
\end{tikzpicture}
}

\newcommand{\MaxRL}{
\begin{tikzpicture}[scale=0.1]
\draw [->] (3,0) to [out=90, in=0] (1.5,1.5) to [out=180, in=90] (0,0);
\end{tikzpicture}
}

\DeclareMathOperator{\Alex}{Alex}
\DeclareMathOperator*{\Bigcdot}{\scalerel*{\cdot}{\bigodot}}
\DeclareMathOperator{\coev}{coev}
\DeclareMathOperator{\dgmod}{dgmod}
\DeclareMathOperator{\ev}{ev}
\DeclareMathOperator{\Hom}{Hom}
\DeclareMathOperator{\id}{id}
\DeclareMathOperator{\incl}{incl}
\DeclareMathOperator{\Mas}{Mas}
\DeclareMathOperator{\multi}{multi}

\DeclareMathOperator{\proj}{proj}
\DeclareMathOperator{\quotient}{quotient}
\DeclareMathOperator{\Rest}{Rest}
\DeclareMathOperator{\rk}{rk}
\DeclareMathOperator{\sing}{sing}
\DeclareMathOperator{\SSp}{SS}
\DeclareMathOperator{\tail}{tail}
\DeclareMathOperator{\TypeD}{TypeD}

\title[Decategorification of Ozsv{\'a}th and Szab{\'o}'s bordered theory]{On the decategorification of Ozsv{\'a}th and Szab{\'o}'s bordered theory for knot Floer homology}
\author{Andrew Manion}
\thanks{This material is based upon work supported by the National Science Foundation under Grant No. DGE-1148900 and DMS-1502686.}

\begin{document}

\begin{abstract} We relate decategorifications of Ozsv{\'a}th--Szab{\'o}'s new bordered theory for knot Floer homology to representations of $\Uq$. Specifically, we consider two subalgebras $\Ccr(n,\Sc)$ and $\Ccl(n,\Sc)$ of Ozsv{\'a}th--Szab{\'o}'s algebra $\B(n,\Sc)$, and identify their Grothendieck groups with tensor products of representations $V$ and $V^*$ of $\Uq$, where $V$ is the vector representation. We identify the decategorifications of Ozsv{\'a}th--Szab{\'o}'s DA bimodules for elementary tangles with corresponding maps between representations. Finally, when the algebras are given multi-Alexander gradings, we demonstrate a relationship between the decategorification of Ozsv{\'a}th--Szab{\'o}'s theory and Viro's quantum relative $\A^1$ of the Reshetikhin--Turaev functor based on $\Uq$.
\end{abstract}

\maketitle

\tableofcontents

\section{Introduction}

\subsection{Knot Floer homology and Ozsv{\'a}th--Szab{\'o}'s theory}

In \cite{OSzNew}, Ozsv{\'a}th and Szab{\'o} define knot invariants $H^-(K)$ and $\widehat{H}(K)$ in the spirit of bordered Heegaard Floer homology \cite{LOTBorderedOrig}. In a forthcoming paper \cite{ComputesHFK}, they identify these invariants with the minus version $HFK^-(K)$ and hat version $\HFK$ of knot Floer homology.

\begin{figure} \centering
\includegraphics[scale=0.625]{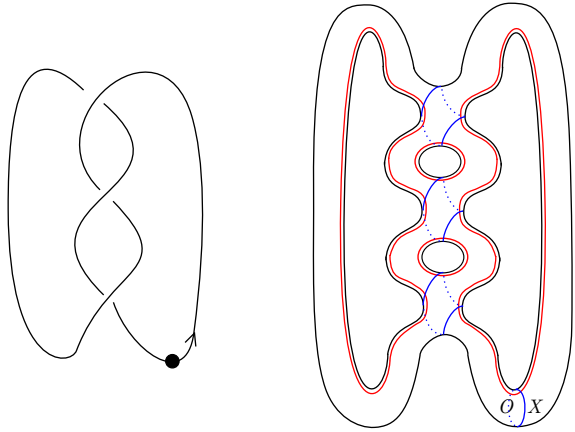}
\caption{Left: projection of the trefoil with a marked point. Right: corresponding Heegaard diagram from \cite{HFAltKnots}.}
\label{FirstHDFig}
\end{figure}

\begin{figure}
\includegraphics[scale=0.5]{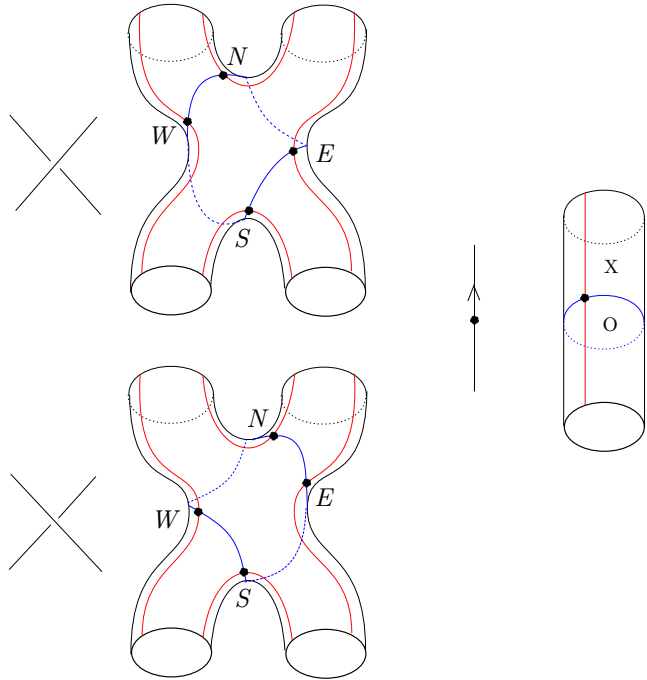}
\caption{Local components of the Heegaard diagram from \cite{HFAltKnots}.}
\label{BasicDiagFig}
\end{figure}

Ozsv{\'a}th and Szab{\'o}'s theory \cite{OSzNew} may be viewed as an algebraic refinement of the model for knot Floer homology provided by the Heegaard diagram defined in \cite{HFAltKnots}; see Figure \ref{FirstHDFig} and Figure \ref{BasicDiagFig}. This Heegaard diagram is naturally associated to an oriented knot projection with a marked point, and generators of the knot Floer complex $CFK^-$ or $\CFK$ of the diagram are in bijection with Kauffman states of the projection as defined in \cite{FKT}.

\subsubsection{Ozsv{\'a}th--Szab{\'o}'s algebras and bimodules}\label{IntroOSzAlgSect}

In \cite{OSzNew}, Ozsv{\'a}th and Szab{\'o} compute the knot Floer complexes by slicing up the Heegaard diagram, or equivalently the knot projection, into simple pieces as in bordered Heegaard Floer homology. We briefly review some properties of their theory; more details are in Section~\ref{OSzReviewSect}.

A generic horizontal slice of an oriented knot projection yields $n$ oriented points arranged along a line. To this configuration, Ozsv{\'a}th and Szab{\'o} assign a differential graded (hereafter dg) algebra $\B(n,\Sc)$ over $\F = \F_2$ with $2^{n+1}$ elementary idempotents, where $\Sc \subset [1,\ldots,n]$ encodes the orientations of the points ($i \in \Sc$ if point $i$ is oriented positively). In contrast with the strands algebras of bordered Floer homology, $\B(n,\Sc)$ is infinite-dimensional over $\F$ when $n \geq 1$. 

As in \cite[Remark 11.13]{OSzNew}, one can define $H^-(K)$ and $\widehat{H}(K)$ using a subalgebra $\A(n,\Sc)$ of $\B(n,\Sc)$ which has $2^{n-1}$ elementary idempotents. We will focus on two intermediate subalgebras $\Ccr(n,\Sc)$ and $\Ccl(n,\Sc)$ which have $2^n$ elementary idempotents. The knot invariants $H^-(K)$ and $\widehat{H}(K)$ can be defined using $\Ccr(n,\Sc)$ or $\Ccl(n,\Sc)$ without modification, just as with $\A(n,\Sc)$.

To local pieces of knot projections, including crossings, maximum points, and minimum points, Ozsv{\'a}th and Szab{\'o} assign finitely generated DA bimodules over the above algebras. These DA bimodules have Maslov gradings (i.e. homological gradings) by $\Z$ and single or multiple intrinsic gradings as described below in Section \ref{OSzReviewSect}. Taking box tensor products of these bimodules over all pieces of a closed knot projection produces a complex whose homology is $H^-(K)$ or $\widehat{H}(K)$ depending on which Type A structure is used for the terminal minimum. In the forthcoming paper \cite{ComputesHFK}, Ozsv{\'a}th--Szab{\'o} show that these invariants agree with knot Floer homology.

\subsection{Decategorification}

The decategorification of Ozsv{\'a}th--Szab{\'o}'s theory should assign an abelian group (more precisely, as we will see, a free module over a group ring) to $n$ oriented points on a line. This group should be the Grothendieck group of a suitable triangulated category associated to $\Ccr(n,\Sc)$ or $\Ccl(n,\Sc)$. The decategorified invariant should also assign linear maps between Grothendieck groups to $(m,n)$--tangles with no closed components.

\subsubsection{Triangulated categories}\label{IntroTriCatsSect}
Given a dg algebra $\A$, one standard triangulated category associated to $\A$ is the compact derived category $\D_c(\A)$, which is the full subcategory of compact objects of the unbounded derived category $\D(\A)$. 

Unfortunately, we are unable to show that the Grothendieck groups of the compact derived categories $\D_c(\Ccr(n,\Sc))$ and $\D_c(\Ccl(n,\Sc))$ have rank $2^n$, which would be ideal. It is possible that they have larger rank, although we believe this is unlikely. 

A more serious issue arises when considering the Type D structure $\Omega^1_r$ over $\Ccr(2,\{1\})$ associated to a single maximum point oriented \MaxLR (the same issue presumably arises for the DA bimodule associated to any maximum or non-terminal minimum). This Type D structure is not operationally bounded. Furthermore, like any dg module, it represents an object in the unbounded derived category, but in Section \ref{CompactSect} we show:
\begin{theorem}[cf. Theorem \ref{NotCompact}]\label{IntroNoncompactThm}
$\Omega^1_r$ is not a compact object of $\D(\Ccr(2,\{1\}))$.
\end{theorem}
Since we want $\Omega^1_r$ to define a class in the Grothendieck group we use, we cannot restrict consideration to the categories $\D_c(\A)$.

We present two ways to work around these issues. The first approach deals only with the DA bimodules representing crossings, not minima or maxima. We can look at homotopy categories $\D^b(\Ccr(n,\Sc))$ and $\D^b(\Ccl(n,\Sc))$ of finitely generated bounded Type D structures. The DA bimodules $X$ for crossings define functors $X \boxtimes -$ between these triangulated categories and linear maps $[X \boxtimes -]$ between their Grothendieck groups. We will compute the Grothendieck groups of these categories and the linear maps for crossing bimodules. 

The second approach applies to all the bordered modules and bimodules of Ozsv{\'a}th--Szab{\'o}'s theory, including the ones for maxima and minima as well as crossings. In Section \ref{SemisimpleSect}, we define ``elementary simple modules'' and ``semisimple modules'' over $\Ccr(n,\Sc)$ and $\Ccr(n,\Sc)$. Isomorphism classes of elementary simple modules are in bijection with the $2^n$ elementary idempotents of $\Ccr(n,\Sc)$ or $\Ccl(n,\Sc)$. 

To any Type D structure (finitely generated but not necessarily bounded) with appropriate gradings, we may associate a class in the Grothendieck group of the homotopy category of semisimple modules. More generally, to a DA bimodule with appropriate gradings, we associate a map between Grothendieck groups.

\subsection{Representations of \texorpdfstring{$\Uq$}{U{\_}q(gl(1|1))} and Viro's quantum relative \texorpdfstring{$\A^1$}{A\textasciicircum{1}}}\label{IntroRepsAndViroSect}

\subsubsection{Representations of \texorpdfstring{$\Uq$}{U{\_}q(gl(1|1))}}\label{IntroUqRepsSect}

If Ozsv{\'a}th--Szab{\'o}'s algebras are given single intrinsic gradings as discussed in Section \ref{OSzSubsect}, we may identify the above Grothendieck groups (after complexification) with tensor products of the vector representation $V$ of the Hopf superalgebra $\Uq$ and its dual $V^*$. The category of finite-dimensional representations of $\Uq$ is a ribbon category; see for example \cite{SartoriAlexander}. Thus, if $\Gamma$ is a crossing, minimum, or maximum, with strands colored by representations of $\Uq$, then $\Gamma$ induces a $\Uq$--linear map between these representations. Using modified bases in tensor products of $V$ and $V^*$, we identify this map with the map on complexified Grothendieck groups associated to one of Ozsv{\'a}th--Szab{\'o}'s DA bimodules.

The $\Uq$--linear map associated to $\Gamma$ goes in the opposite direction compared with the decategorification maps discussed above.  To mediate between the different conventions, we relate the decategorification map for a tangle $\Gamma$ with the representation-theoretic map for a rotated and reversed tangle $\widetilde{\Gamma}$, and vice versa; see Figure \ref{RotateReverseFig} below. 

\subsubsection{Viro's quantum relative \texorpdfstring{$\A^1$}{A\textasciicircum{1}}}

In the multi-graded case, we will relate the decategorification with a variant of this $\Uq$--representation setup introduced by Viro in \cite{Viro}. See also \cite{RozSal} for a similar construction from a physical perspective.

Viro's invariant $\A^1_P$ depends on a choice of ``1-palette'' $P$, described in Section \ref{ViroPaletteSect} below. Given a suitable palette (Definition \ref{UnivPaletteObjs}) and coloring (Definition \ref{UnivColorObjects}), we will show that Viro's invariant for $n$ oriented points along a line may be identified with the Grothendieck groups associated to Ozsv{\'a}th--Szab{\'o}'s algebra $\Ccrl(n,\Sc)$, given multiple intrinsic gradings as in Section \ref{OSzSubsect}. 

Furthermore, the maps Viro associates to tangles (with the palette of Definition \ref{UnivPaletteGraphs} and the coloring of Definition \ref{UnivColorGraphs}) may be identified with the linear maps on Grothendieck groups coming from decategorifying Ozsv{\'a}th--Szab{\'o}'s tangle bimodules. For maxima and minima, the maps from Viro's construction and from decategorification are only identified up to a scalar multiple. As in the singly-graded case, this identification uses modified bases, described below in Definition \ref{MultivarModifiedBasisDef}, as well as the rotation and orientation reversal shown in Figure \ref{RotateReverseFig}.

In \cite[Section 2.7]{Viro}, Viro relates the functor $\A^1_{P_W}$ to maps induced on tensor products of more general irreducible representations of $\Uq$, rather than just the vector representation and its dual. He also shows in \cite[Sections 2.10, 11.7--11.8]{Viro} that the image under $\A^1_P$ of a collection of points may be upgraded to a representation of $U^1 \otimes_{\Z} B$, where $U^1$ is a certain Hopf $\Z$--subalgebra of $\Uq$. A discussion of $U^1$ and its representation theory may be found in \cite[Sections 11.7--11.8]{Viro}.

\subsection{Main results; relationship with other work; outline of paper.}\label{MainResultsSect}

To conclude the introduction, we state our main results somewhat imprecisely, with references to the precise statements below. 

Let $\Ccrl(n,\Sc)$ denote either $\Ccr(n,\Sc)$ or $\Ccl(n,\Sc)$. In Definition \ref{BddTypeDCatsDef}, we associate a triangulated category $\D^b(\Ccrl(n,\Sc))$ to $\Ccrl(n,\Sc)$, defined using finitely generated bounded Type D structures. In Definition \ref{SingleSSCatsDef}, we associate a triangulated category $\SSp^{\sing}(\Ccrl(n,\Sc))$ to $\Ccrl(n,\Sc)$, defined using semisimple modules.

\begin{theorem}[cf. Theorems \ref{SingleAlgDecatBodyThm}, \ref{SSAlgDecatProp}]\label{AlgDecatIntroThm}
The Grothendieck group $K_0(\D^b(\Ccrl(n,\Sc)))$ is a free module of rank $2^n$ over $\Z[t^{\pm 1/2}]$.  Similarly, $K_0(\SSp^{\sing}(\Ccrl(n,\Sc)))$ is a free module of rank $2^n$ over $\Z[t^{\pm 1/2}]$. 
\end{theorem}

\begin{theorem}[cf. Theorem~\ref{MainThmFinalForm}]\label{IntroMainThm} 
Make the identifications
\begin{gather*}
K_0(\D^b(\Ccrl(n,\Sc))) \otimes_{\frac{\Z[t^{\pm 1/2},q^{\pm 1}]}{q=t^{1/2}}} \C(q) \cong V^{\otimes \Sc} \\
\tag*{\text{and}} K_0(\SSp^{\sing}(\Ccrl(n,\Sc))) \otimes_{\frac{\Z[t^{\pm 1/2},q^{\pm 1}]}{q=t^{1/2}}} \C(q) \cong V^{\otimes \Sc}
\end{gather*}
using modified bases.

For a tangle $\Gamma$ that is not the terminal minimum, let $\widetilde{\Gamma}$ be the tangle obtained by rotating $\Gamma$ and reversing the orientations on strands as in Figure \ref{RotateReverseFig}. The decategorification of Ozsv{\'a}th--Szab{\'o}'s DA bimodule for $\Gamma$ agrees with the $\Uq$--linear map associated to $\widetilde{\Gamma}$ under the above identifications.
\end{theorem}
 
In Definition \ref{BddTypeDCatsDef} we associate a triangulated category $\D^b_{M_Y}(\Ccrl(n,\Sc))$ to $\Ccrl(n,\Sc)$. In Definition \ref{MultiSSCatsDef} we associate a triangulated category $\SSp^{\multi}_Y(\Ccrl(n,\Sc))$ to $\Ccrl(n,\Sc)$.

\begin{theorem}[cf. Theorem \ref{MultiAlgDecatBodyThm}, Theorem \ref{SSAlgDecatProp}]\label{MultiAlgDecatIntroThm}
The Grothendieck group 
\[
K_0(\D^b_{M_Y}(\Ccrl(n,\Sc)))
\]
is a free module of rank $2^n$ over $\Z[M_Y]$, where $M_Y := H^0(Y; \frac{1}{4}\Z)$. Similarly, $K_0(\SSp^{\multi}_Y(\Ccrl(n,\Sc)))$ is a free module of rank $2^n$ over $\Z[M_Y]$. 
\end{theorem}

\begin{theorem}[cf. Theorem~\ref{MultiMainThmFinalForm}]\label{MultiIntroMainThm} 
Make the identifications
\begin{gather*}
K_0(\D^b_{M_Y}(\Ccrl(n,\Sc))) \cong \A^1_{P_Y}(Y) \\
\tag*{\text{and}} K_0(\SSp^{\multi}_Y(\Ccrl(n,\Sc))) \cong \A^1_{P_Y}(Y)
\end{gather*}
using modified bases.

For a tangle $\Gamma$ that is not the terminal minimum, let $\widetilde{\Gamma}$ be the tangle obtained by rotating $\Gamma$ and reversing the orientations on strands as in Figure \ref{RotateReverseFig}. The decategorification of Ozsv{\'a}th--Szab{\'o}'s DA bimodule for $\Gamma$ agrees with Viro's map associated to $\widetilde{\Gamma}$ under the above identifications.

\end{theorem}

\subsubsection{Relationship with other work}

As a consequence of these theorems, Ozsv{\'a}th and Szab{\'o}'s construction in \cite{OSzNew} may be regarded as a categorification of tensor powers of the fundamental representation $V$ of $\Uq$ and its dual, together with the maps between these representations induced by crossings, maxima, and minima. 

In \cite{Sartori}, Sartori produces a categorification of tensor powers involving only $V$ and not $V^*$. He constructs functors between his categories which categorify the maps for crossings between positively oriented strands. He does not write down bimodules representing his functors; on the other hand, he has a categorical action of $\Uq$, an element that is not explicitly present in Ozsv{\'a}th--Szab{\'o}'s theory.

In one next-to-extremal weight space of $V^n$, Sartori has explicit bimodules which represent his categorification functors. His algebra categorifying this weight space is the Khovanov--Seidel quiver algebra from \cite{KSQuiver}, and his functors are represented by tensor product with Khovanov--Seidel's chain complexes of bimodules. In \cite{MyKSQuiverPaper} we show that Khovanov--Seidel's algebra is a quotient of Ozsv{\'a}th--Szab{\'o}'s algebra $\Ccl(n,1,\varnothing)$ and obtain a close relationship between Khovanov--Seidel's dg bimodules and Ozsv{\'a}th--Szab{\'o}'s DA bimodules in this case. Thus, \cite{MyKSQuiverPaper} can be seen as a first step toward extending the current paper beyond the level of decategorification. 

In \cite{TianUT}, Tian also categorifies tensor powers of the fundamental representation $V$ of an algebra he calls $\UT$, a close relative of $\Uq$. He has categorical actions of quantum group generators as well, but no bimodules or functors for a crossing. Indecomposable projective modules in his categories correspond to elements of the standard or tensor-product basis of $V^n$, which differs from the modified bases we use when $n > 1$. 

In \cite{PetkovaVertesi}, Petkova and V{\'e}rtesi define a theory with similar properties to Ozsv{\'a}th and Szab{\'o}'s, which they call tangle Floer homology. Ellis, Petkova, and V{\'e}rtesi \cite{EPV} show that the dg algebras of \cite{PetkovaVertesi} categorify a tensor product of irreducible representations of $\Uq$, each of which is $V$ or $V^*$ except for one factor $L(\lambda_{n+1})$ which is neither $V$ nor $V^*$. The Heegaard diagrams motivating these two theories are different; however, it would be interesting to see whether relationships between the theories exist, especially since both theories are known to compute knot Floer homology in some form.

Finally, in \cite{Zibrowius, Zibrowius2}, Zibrowius defines an Alexander invariant for tangles based on counting ``generalized Kauffman states'', which are very similar to Ozsv{\'a}th--Szab{\'o}'s partial Kauffman states (``sites'' in Zibrowius' language correspond more-or-less to idempotents of Ozsv{\'a}th--Szab{\'o}'s algebras). He uses sutured Floer homology to categorify his invariant; as with Petkova--V{\'e}rtesi's theory, it would be interesting to study relationships between his categorification and Ozsv{\'a}th--Szab{\'o}'s algebras and bimodules. 

\subsubsection{Outline of paper.} 
Section \ref{UqRepSect} discusses what we need about representations of $\Uq$, including basic definitions in Section \ref{UqBackgroundSect}, linear maps for tangles in Section \ref{MapsForCrossingsSect} and \ref{MinimaMaximaSect}, modified bases in Section \ref{BasesSect}, and matrix representations in modified bases in Section \ref{MapsInModifiedBasisSect}. Section \ref{ViroSect} does the same for Viro's quantum relative $\A^1$, with a similar outline. In Section \ref{OSzReviewSect}, we review various details of Ozsv{\'a}th--Szab{\'o}'s construction from \cite{OSzNew}. Section \ref{OSzSubsect} discusses the algebras and Sections \ref{OSzCrossingSect}--\ref{OSzTerminalSect} discuss the DA bimodules.

Section \ref{AlgDecatSect} deals with decategorification. After some definitions in Section \ref{AbstractTriCatSect}, we prove Theorem \ref{IntroNoncompactThm} in Section \ref{CompactSect}. We prove Theorems \ref{AlgDecatIntroThm} and \ref{MultiAlgDecatIntroThm} partially in Section \ref{BddTypeDSect} and partially in Section \ref{SemisimpleSect}. Finally, we compute explicit matrices by decategorifying Ozsv{\'a}th--Szab{\'o}'s DA bimodules in Sections \ref{ExplicitDecatCrossingSect}--\ref{TTASect}. In Section \ref{RepThyCompareSect}, we compare the decategorification matrices with modified-basis matrices from Sections \ref{UqRepSect} and \ref{ViroSect}, proving Theorems \ref{IntroMainThm} and \ref{MultiIntroMainThm}.

\section*{Acknowledgments} 
I would especially like to thank Rapha{\"e}l Rouquier for extensive guidance, especially with Section \ref{CompactSect}, and my graduate advisor Zolt{\'a}n Szab{\'o} for teaching me about the new theory. I would also like to thank Robert Lipshitz, Ciprian Manolescu, and Antonio Sartori for useful discussions and suggestions.

\section{\texorpdfstring{$\Uq$}{U{\_}q(gl(1|1))} and its representations}\label{UqRepSect}

\subsection{Background}\label{UqBackgroundSect}
The definition of Hopf superalgebras may be found in \cite[Section 2.3]{SartoriAlexander}. Our conventions, for the most part, will follow Sartori's.

\begin{definition}\label{UqDef} $\Uq$ is the Hopf superalgebra over $\C(q)$ with even generators $K_1^{\pm 1}$ and $K_2^{\pm 1}$, odd generators $E$ and $F$, relations specified in \cite[Section 2.2]{SartoriAlexander} (with $K_i = \mathbf{q}^{h_i}$ in Sartori's notation), and coproduct, counit, and antipode specified in \cite[Section 2.4]{SartoriAlexander}.
\end{definition}

\begin{definition} The vector representation $V$ of $\Uq$ is the $1|1$--dimensional super vector space over $\C(q)$ with one even generator $v_0$, one odd generator $v_1$, and action of $\Uq$ defined in \cite[Section 3.4]{SartoriAlexander}.
\end{definition}

Since $\Uq$ is a Hopf superalgebra, for each finite-dimensional representation $W$ of $\Uq$, there is a dual representation
\[
W^* := \Hom_{\C(q)}(W,\C(q))
\]
whose decomposition into even and odd parts comes from the degrees of linear maps from $W$ to $\C(q)$, with $\C(q)$ placed in even degree. Also, we can construct tensor products of representations using the coproduct. For a representation $W$ of $\Uq$, we will let $W^{\otimes n}$ denote $W$ tensored with itself $n$ times.

In formula (3.4) of \cite{SartoriAlexander}, Sartori describes the action of $\Uq$ on a general family of representations $L(\lambda)$, with two generators $v_0^{\lambda}$ and $v_1^{\lambda}$, depending on a parameter 
\[
\lambda = \lambda_1 \epsilon_1 + \lambda_2 \epsilon_2 \in \Z \epsilon_1 \oplus \Z \epsilon_2
\]
with $\lambda_1 + \lambda_2 \neq 0$. Sartori defines the parity $|\lambda| \in \Z/2\Z$ of the parameter $\lambda \in \Z \epsilon_1 \oplus \Z \epsilon_2$ by 
\[
|\lambda_1 \epsilon_1 + \lambda_2 \epsilon_2| := \lambda_2 \textrm{ modulo } 2.
\]
If $|\lambda|=0$, then $v_0^{\lambda}$ is considered an even generator and $v_1^{\lambda}$ is considered an odd generator. If $|\lambda|=1$, then the reverse is true. In particular, $v_0^{\epsilon_1}$, $(v_0^{\epsilon_1})^*$, and $v_1^{-\epsilon_2}$ are even, and $v_1^{\epsilon_1}$, $(v_1^{\epsilon_1})^*$, and $v_0^{-\epsilon_2}$ are odd. The vector representation $V$ is $L(\epsilon_1)$, and its dual $V^*$ can be identified with $L(-\epsilon_2)$. The explicit identification sends
\begin{equation}\label{Lne2VstarIdent}
\begin{aligned}
&(v_0^{\epsilon_1})^* \leftrightarrow v_1^{-\epsilon_2}, \\
&(v_1^{\epsilon_1})^* \leftrightarrow -q^{-1} v_0^{-\epsilon_2}.
\end{aligned}
\end{equation}

\begin{definition}
Let $\Sc \subset [1, \ldots, n]$, and let
\[
V_i := \begin{cases} V & i \in \Sc \\ V^* & i \notin \Sc. \end{cases}
\] 
Then $V^{\otimes \Sc}$ is the representation 
\[
V^{\otimes \Sc} := V_1 \otimes \cdots \otimes V_n
\]
of $\Uq$.
\end{definition}

\subsection{Maps for crossings}\label{MapsForCrossingsSect}
Let $\Gamma$ be an $(m,n)$--tangle; by convention, $\Gamma$ has $m$ strands on the bottom and $n$ strands on top. 
\begin{definition}
Given $\Gamma$, let $\Sc_{\Gamma} \subset [1, \ldots, n]$ be the indices of the upward-pointing strands at the top of $\Gamma$, and $\Sc'_{\Gamma} \subset [1, \ldots, m]$ be the indices of the upward-pointing strands at the bottom of $\Gamma$. 
\end{definition}
Finite-dimensional representations of $\Uq$ form a ribbon category, as proved e.g. in \cite{SartoriAlexander}, so to $\Gamma$ we may associate a $\Uq$--linear map from $V^{\otimes \Sc'_{\Gamma}}$ to $V^{\otimes \Sc_{\Gamma}}$. We review the definition of this map below. 

\begin{remark}
As mentioned in Section \ref{IntroUqRepsSect}, the decategorification of Ozsv{\'a}th and Szab{\'o}'s DA bimodule for a tangle $\Gamma$ will correspond not to the $\Uq$--linear map induced by $\Gamma$, but rather to the map induced by a related tangle $\widetilde{\Gamma}$, obtained by rotating $\Gamma$ and reversing the orientations on its strands as in Figure~\ref{RotateReverseFig} below. 
\end{remark}

\begin{figure} \centering
\includegraphics[scale=0.625]{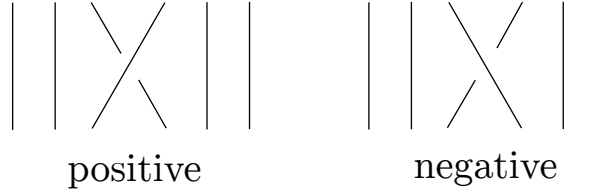}
\caption{Positive and negative braid generators.}
\label{PosNegBraidsFig}
\end{figure}

We begin by defining the $\Uq$--linear maps when $\Gamma$ is a $(2,2)$--tangle consisting of a single crossing that is positive as a braid generator (with any orientations); see Figure \ref{PosNegBraidsFig}. The maps are denoted $\check{R}$ in this case since they come from the $R$--matrix giving the quasitriangular structure of $\Uq$. Crossings that are negative as braid generators are assigned the inverses maps $\check{R}^{-1}$. We will be ambiguous with notation and use $\check{R}$ to denote any of the maps for positive braid generators, regardless of the representations on which they act.

There is a formula in \cite{SartoriAlexander} for $\check{R}^{-1}$ acting on $L(\lambda) \otimes L(\mu)$; here we write out the formula for $\check{R}$ explicitly when $\lambda, \mu \in \{\epsilon_1, -\epsilon_2\}$.
\begin{definition}[cf. \cite{SartoriAlexander}, formula (4.7)]
When the orientations are \OrUU, the map 
\[
\check{R}: L(\epsilon_1) \otimes L(\epsilon_1) \to L(\epsilon_1) \otimes L(\epsilon_1)
\]
is defined by
\[
\check{R} = 
\kbordermatrix{ 
& v_0^{\epsilon_1} \otimes v_0^{\epsilon_1} & v_1^{\epsilon_1} \otimes v_0^{\epsilon_1} & v_0^{\epsilon_1} \otimes v_1^{\epsilon_1} & v_1^{\epsilon_1} \otimes v_1^{\epsilon_1} \\ 
v_0^{\epsilon_1} \otimes v_0^{\epsilon_1} & q & 0 & 0 & 0 \\ 
v_1^{\epsilon_1} \otimes v_0^{\epsilon_1} & 0 & 0 & 1 & 0 \\ 
v_0^{\epsilon_1} \otimes v_1^{\epsilon_1} & 0 & 1 & q-q^{-1} & 0 \\
v_1^{\epsilon_1} \otimes v_1^{\epsilon_1} & 0 & 0 & 0 & -q^{-1}
}.
\]
When the orientations are \OrDU, the map 
\[
\check{R}: L(-\epsilon_2) \otimes L(\epsilon_1) \to L(\epsilon_1) \otimes L(-\epsilon_2)
\]
is defined by
\[
\check{R} = 
\kbordermatrix{ 
& v_0^{-\epsilon_2} \otimes v_0^{\epsilon_1} & v_1^{-\epsilon_2} \otimes v_0^{\epsilon_1} & v_0^{-\epsilon_2} \otimes v_1^{\epsilon_1} & v_1^{-\epsilon_2} \otimes v_1^{\epsilon_1} \\ 
v_0^{\epsilon_1} \otimes v_0^{-\epsilon_2} & 1 & 0 & 0 & 0 \\ 
v_1^{\epsilon_1} \otimes v_0^{-\epsilon_2} & 0 & 0 & -q & 0 \\ 
v_0^{\epsilon_1} \otimes v_1^{-\epsilon_2} & 0 & q^{-1} & q^2-1 & 0 \\
v_1^{\epsilon_1} \otimes v_1^{-\epsilon_2} & 0 & 0 & 0 & 1
}.
\]
When the orientations are \OrUD, the map 
\[
\check{R}: L(\epsilon_1) \otimes L(-\epsilon_2) \to L(-\epsilon_2) \otimes L(\epsilon_1)
\]
is defined by
\[
\check{R} = 
\kbordermatrix{ 
& v_0^{\epsilon_1} \otimes v_0^{-\epsilon_2} & v_1^{\epsilon_1} \otimes v_0^{-\epsilon_2} & v_0^{\epsilon_1} \otimes v_1^{-\epsilon_2} & v_1^{\epsilon_1} \otimes v_1^{-\epsilon_2} \\ 
v_0^{-\epsilon_2} \otimes v_0^{\epsilon_1} & 1 & 0 & 0 & 0 \\ 
v_1^{-\epsilon_2} \otimes v_0^{\epsilon_1} & 0 & 0 & q^{-1} & 0 \\ 
v_0^{-\epsilon_2} \otimes v_1^{\epsilon_1} & 0 & -q & q^{-2} - 1 & 0 \\
v_1^{-\epsilon_2} \otimes v_1^{\epsilon_1} & 0 & 0 & 0 & 1
}.
\]
When the orientations are \OrDD, the map 
\[
\check{R}: L(-\epsilon_2) \otimes L(-\epsilon_2) \to L(-\epsilon_2) \otimes L(-\epsilon_2)
\]
is defined by
\[
\check{R} = 
\kbordermatrix{ 
& v_0^{-\epsilon_2} \otimes v_0^{-\epsilon_2} & v_1^{-\epsilon_2} \otimes v_0^{-\epsilon_2} & v_0^{-\epsilon_2} \otimes v_1^{-\epsilon_2} & v_1^{-\epsilon_2} \otimes v_1^{-\epsilon_2} \\ 
v_0^{-\epsilon_2} \otimes v_0^{-\epsilon_2} & -q^{-1} & 0 & 0 & 0 \\ 
v_1^{-\epsilon_2} \otimes v_0^{-\epsilon_2} & 0 & 0 & 1 & 0 \\ 
v_0^{-\epsilon_2} \otimes v_1^{-\epsilon_2} & 0 & 1 & q-q^{-1} & 0 \\
v_1^{-\epsilon_2} \otimes v_1^{-\epsilon_2} & 0 & 0 & 0 & q
}.
\]
Note that $\check{R}$ sends even vectors to even vectors and odd vectors to odd vectors, so it is an even operator.
\end{definition}

We rewrite the above matrices in terms of $V$ and $V^*$ using (\ref{Lne2VstarIdent}):
\begin{corollary}\label{VVstarRmatrices}
When the orientations are \OrUU, the map 
\[
\check{R}: V \otimes V \to V \otimes V
\]
has matrix
\[
\check{R} = 
\kbordermatrix{ 
& v_0 \otimes v_0 & v_1 \otimes v_0 & v_0 \otimes v_1 & v_1 \otimes v_1 \\ 
v_0 \otimes v_0 & q & 0 & 0 & 0 \\ 
v_1 \otimes v_0 & 0 & 0 & 1 & 0 \\ 
v_0 \otimes v_1 & 0 & 1 & q-q^{-1} & 0 \\
v_1 \otimes v_1 & 0 & 0 & 0 & -q^{-1}
}.
\]
When the orientations are \OrDU, the map 
\[
\check{R}: V^* \otimes V \to V \otimes V^*
\]
has matrix
\[
\check{R} = 
\kbordermatrix{ 
& v_1^* \otimes v_0 & v_0^* \otimes v_0 & v_1^* \otimes v_1 & v_0^* \otimes v_1 \\ 
v_0 \otimes v_1^* & 1 & 0 & 0 & 0 \\ 
v_1 \otimes v_1^* & 0 & 0 & -q & 0 \\ 
v_0 \otimes v_0^* & 0 & q^{-1} & q^{-1} - q & 0 \\
v_1 \otimes v_0^* & 0 & 0 & 0 & 1
}.
\]
When the orientations are \OrUD, the map 
\[
\check{R}: V \otimes V^* \to V^* \otimes V
\]
has matrix
\[
\check{R} = 
\kbordermatrix{ 
& v_0 \otimes v_1^* & v_1 \otimes v_1^* & v_0 \otimes v_0^* & v_1 \otimes v_0^* \\ 
v_1^* \otimes v_0 & 1 & 0 & 0 & 0 \\ 
v_0^* \otimes v_0 & 0 & 0 & q^{-1} & 0 \\ 
v_1^* \otimes v_1 & 0 & -q & q - q^{-1} & 0 \\
v_0^* \otimes v_1 & 0 & 0 & 0 & 1
}.
\]
When the orientations are \OrDD, the map 
\[
\check{R}: V^* \otimes V^* \to V^* \otimes V^*
\]
has matrix
\[
\check{R} = 
\kbordermatrix{ 
& v_1^* \otimes v_1^* & v_0^* \otimes v_1^* & v_1^* \otimes v_0^* & v_0^* \otimes v_0^* \\ 
v_1^* \otimes v_1^* & -q^{-1} & 0 & 0 & 0 \\ 
v_0^* \otimes v_1^* & 0 & 0 & 1 & 0 \\ 
v_1^* \otimes v_0^* & 0 & 1 & q-q^{-1} & 0 \\
v_0^* \otimes v_0^* & 0 & 0 & 0 & q
}. 
\]
\end{corollary}

\begin{definition} Let $1 \leq i \leq n-1$, and let $\Gamma$ be an $(n,n)$--tangle consisting of a single crossing between strands $i$ and $i+1$ that is positive as a braid generator.

To $\Gamma$, associate the map 
\begin{gather*}
\check{R}_{i,i+1}: V^{\otimes \Sc'_{\Gamma}} \to V^{\otimes \Sc_{\Gamma}} \\
\tag*{\text{defined as}} \id \otimes \ldots \otimes \id \otimes \check{R} \otimes \id \otimes \ldots \id,
\end{gather*}
where $\check{R}$ acts on tensor factors $i$ and $i+1$. If $\Gamma$ is negative rather than positive as a braid generator, associate to $\Gamma$ the map
\[
\check{R}^{-1}_{i,i+1} := \id \otimes \ldots \otimes \id \otimes \check{R}^{-1} \otimes \id \otimes \ldots \id.
\]
\end{definition}

\subsection{Maps for minima and maxima}\label{MinimaMaximaSect}
To define the maps for maxima and minima, we recall the following definitions from \cite{SartoriAlexander}:
\begin{definition}[cf. Equations 4.1 and 4.2 of \cite{SartoriAlexander}]\label{EvCoevDef} Let $W$ be a finite-dimensional representation of $\Uq$, and let $\{w_i\}$ be a basis for $W$. The coevaluation maps
\begin{align*}
&\coev_W: \C(q) \to W \otimes W^*, \\
&\widehat{\coev}_W: \C(q) \to W^* \otimes W
\end{align*}
are defined by
\begin{align*}
&\coev_W(1) := \sum_i w_i \otimes w_i^*, \\
&\widehat{\coev}_W(1) := \sum_i (-1)^{|w_i|} w_i^* \otimes w_i;
\end{align*}
one can show that these maps are independent of the chosen basis. The evaluation maps
\begin{align*}
&\ev_W: W^* \otimes W \to \C(q), \\
& \widehat{\ev}_W: W \otimes W^* \to \C(q)
\end{align*}
are defined by
\begin{align*}
&\ev_W(\varphi \otimes w) := \varphi(w), \\
&\widehat{\ev_W}(\varphi \otimes w) := (-1)^{|\varphi||w|} \varphi(w).
\end{align*}
\end{definition}

\begin{example}\label{EvCoevMatrices} When $W$ in Definition~\ref{EvCoevDef} is the vector representation $V$, we write out explicit matrices for the coevaluation and evaluation maps:
$\coev_V$ has matrix
\[
\coev_V = 
\kbordermatrix{ 
& 1 \\ 
v_0 \otimes v_1^* & 0  \\ 
v_1 \otimes v_1^* & 1  \\ 
v_0 \otimes v_0^* & 1  \\
v_1 \otimes v_0^* & 0 
},
\]
$\widehat{\coev}_V$ has matrix
\[
\widehat{\coev}_V = 
\kbordermatrix{ 
& 1 \\ 
v_1^* \otimes v_0 & 0  \\ 
v_0^* \otimes v_0 & 1  \\ 
v_1^* \otimes v_1 & -1  \\
v_0^* \otimes v_1 & 0 
},
\]
$\ev_V$ has matrix
\[
\ev_V = \qquad
\kbordermatrix{ 
& v_1^* \otimes v_0 & v_0^* \otimes v_0 & v_1^* \otimes v_1 & v_0^* \otimes v_1 \\ 
1 & 0 & 1 & 1 & 0 \\ 
},
\]
and $\widehat{\ev}_V$ has matrix
\[
\widehat{\ev}_V = \qquad
\kbordermatrix{ 
& v_0 \otimes v_1^* & v_1 \otimes v_1^* & v_0 \otimes v_0^* & v_1 \otimes v_0^* \\ 
1 & 0 & -1 & 1 & 0 \\ 
}.
\]
\end{example}

\subsubsection{Minima}\label{UqMinimaSect}
Let $\Gamma$ be a $(0,2)$--tangle. If $\Gamma$ is oriented \MinRL, then $\Gamma$ is assigned the map
\[
\coev_V: \C(q) \to V \otimes V^*,
\]
with matrix given above in Example \ref{EvCoevMatrices}.

If $\Gamma$ is oriented \MinLR, we need additional data. Consider the element $u$ of $\Uh$ defined in equation (2.21) of \cite{SartoriAlexander}. We will not need to define $u$ or $\Uh$ precisely; all we will use is that $u$ defines a map $W \to W$ for each representation $W$ of $\Uq$, and that this map is multiplication by $q$ on $V$ and multiplication by $q^{-1}$ on $V^*$. These facts can be checked from \cite[equation (2.23)]{SartoriAlexander}.

Another element $v$ of $\Uh$ gives $\Uh$ the structure of a (topological) ribbon Hopf superalgebra; see \cite[Proposition 2.5]{SartoriAlexander}. The element $v$ also acts in each representation of $\Uq$, and it acts as the identity on $V$ and $V^*$ by \cite[Lemma 4.2]{SartoriAlexander}.

Now, the $(0,2)$--tangle $\Gamma$ oriented \MinLR is assigned the map
\[
(\id \otimes vu^{-1}) \circ \widehat{\coev}_V: \C(q) \to V^* \otimes V.
\]
Since $u^{-1}$ acts as multiplication by $q^{-1}$ on $V$, and $v$ acts as the identity, we can write the matrix for the map associated to $\Gamma$ as
\[
\kbordermatrix{ 
& 1 \\ 
v_1^* \otimes v_0 & 0  \\ 
v_0^* \otimes v_0 & q^{-1}  \\ 
v_1^* \otimes v_1 & -q^{-1}  \\
v_0^* \otimes v_1 & 0 
}.
\]

\begin{definition}
Let $0 \leq i \leq n$, and let $\Gamma$ be an $(n,n+2)$--tangle with no crossings and a single minimum point between strands $i+1$ and $i+2$. If $\Gamma$ is oriented \MinRL, then $\Gamma$ is assigned the map
\[
\id \otimes \cdots \otimes \id \otimes \coev_V \otimes \id \otimes \cdots \otimes \id
\]
from $V^{\otimes \Sc'_{\Gamma}}$ to $V^{\otimes \Sc_{\Gamma}}$, where there are $i$ instances of $\id$ to the left of $\coev_V$. If $\Gamma$ is oriented \MinLR, then $\Gamma$ is assigned the map
\[
\id \otimes \cdots \otimes \id \otimes ((\id \otimes vu^{-1}) \circ \widehat{\coev}_V) \otimes \id \otimes \cdots \otimes \id
\]
from $V^{\otimes \Sc'_{\Gamma}}$ to $V^{\otimes \Sc_{\Gamma}}$, where there are $i$ instances of $\id$ to the left of $(\id \otimes vu^{-1}) \circ \widehat{\coev}_V$.
\end{definition}

\subsubsection{Maxima}\label{UqMaximaSect}
Let $\Gamma$ be a $(2,0)$--tangle. If $\Gamma$ is oriented \MaxRL, then $\Gamma$ is assigned the map
\[
\ev_V: V^* \otimes V \to \C(q),
\]
with matrix given above in Example~\ref{EvCoevMatrices}.

If $\Gamma$ is oriented \MaxLR, then $\Gamma$ is assigned the map
\[
\widehat{\ev}_V \circ (uv^{-1} \otimes \id): V \otimes V^* \to \C(q).
\]
Since $u$ acts as multiplication by $q$ on $V$ and $v$ acts as the identity, we can write the matrix for the map associated to $\Gamma$ as
\[ 
\kbordermatrix{ 
& v_0 \otimes v_1^* & v_1 \otimes v_1^* & v_0 \otimes v_0^* & v_1 \otimes v_0^* \\ 
1 & 0 & -q & q & 0 \\ 
}.
\]

\begin{definition} Let $0 \leq i \leq n$, and let $\Gamma$ be an $(n+2,n)$--tangle with no crossings and a single maximum point between strands $i+1$ and $i+2$. If $\Gamma$ is oriented \MaxRL, then $\Gamma$ is assigned the map
\[
\id \otimes \cdots \otimes \id \otimes \ev_V \otimes \id \otimes \cdots \otimes \id
\]
from $V^{\otimes \Sc'_{\Gamma}}$ to $V^{\otimes \Sc_{\Gamma}}$, where there are $i$ instances of $\id$ to the left of $\ev_V$. If $\Gamma$ is oriented \MaxLR, then $\Gamma$ is assigned the map
\[
\id \otimes \cdots \otimes \id \otimes (\widehat{\ev}_V \circ (uv^{-1} \otimes \id)) \otimes \id \otimes \cdots \otimes \id,
\]
from $V^{\otimes \Sc'_{\Gamma}}$ to $V^{\otimes \Sc_{\Gamma}}$, where there are $i$ instances of $\id$ to the left of $(\widehat{\ev}_V \circ (uv^{-1} \otimes \id))$.
\end{definition}

\subsection{Bases for \texorpdfstring{$V^{\otimes \Sc}$}{V\textasciicircum{S}}}\label{BasesSect}

We now define two modified bases for $V^{\otimes \Sc}$. We start by discussing the standard tensor-product basis in a notation resembling that of \cite[Section 5]{KSFF}. 

\subsubsection{Tensor product basis}

\begin{definition}
Let $W$ be a finite-dimensional vector space over a field $k$. The exterior algebra $\wedge^* W$ can be viewed as a super vector space with even part
\begin{gather*}
\bigoplus_{m \geq 0} \wedge^{2m} W \\
\tag*{\text{and odd part}} \bigoplus_{m \geq 0} \wedge^{2m+1} W.
\end{gather*}
\end{definition} 

\begin{definition}\label{WiDef}
For $1 \leq i \leq n$, let $w_i$ denote the tensor-product element of $V^{\otimes \Sc}$ that, in the $j^{th}$ place for $j \neq i$, has $v_0$ if $j \in \Sc$ and $v_1^*$ if $j \notin \Sc$, and in the $i^{th}$ place, has $v_1$ if $i \in \Sc$ and $v_0^*$ if $i \notin \Sc$.
\end{definition}
Let $\W$ denote the $n$--dimensional vector space over $\C(q)$ formally generated by the elements $w_i$ of Definition \ref{WiDef}. Let $\eta$ be the tensor-product element of $V^{\otimes \Sc}$ that, in each place $i$, has $v_0$ if $i \in \Sc$ and $v_1^*$ if $i \notin \Sc$. 

First consider the case when $\eta$ is an even element of $V^{\otimes \Sc}$ (for instance, $\eta$ is even when $V^{\otimes \Sc} = V^{\otimes n}$). Then, as super vector spaces, we have
\[
V^{\otimes \Sc} \cong \wedge^* \W.
\]
The generator $w_{i_1} \wedge \ldots \wedge w_{i_k}$ of $\wedge^* \W$ is identified with the generator of $V^{\otimes \Sc}$ that has $v_1$ or $v_0^*$ in places $i_1,\ldots,i_k$ and $v_0$ or $v_1^*$ in all other spots.

Now suppose $\eta$ is an odd element of $V^{\otimes \Sc}$. Then we have
\begin{gather*}
(V^{\otimes \Sc})_{odd} = \bigoplus_{m \geq 0} \wedge^{2m} \W \\
\tag*{\text{and}} (V^{\otimes \Sc})_{even} = \bigoplus_{m \geq 0} \wedge^{2m+1} \W.
\end{gather*}
The generators $\{w_{i_1} \wedge \ldots \wedge w_{i_k}\}$ still form a basis for $V^{\otimes \Sc}$ as a super vector space. To keep track of odd and even parts, we must simply remember that if $\eta$ is odd, then even-degree wedge products are odd and odd-degree wedge products are even.

The basis $\{w_{i_1} \wedge \ldots \wedge w_{i_k}\}$ is just the usual basis consisting of tensors of $v_0, v_1, v_0^*$, and $v_1^*$ as appropriate; we have only changed notation at this point. This tensor-product basis, or wedge-product basis, is the first basis for $V^{\otimes \Sc}$ we will consider.

We call $1$, $w_i$, $w_{i+1}$, and $w_i \wedge w_{i+1}$ the ``local generators'' for a crossing between strands $i$ and $i+1$. Rather than write the full $2^n \times 2^n$ matrix for 
\[
\check{R}_{i,i+1} = \id \otimes \ldots \otimes \id \otimes \check{R} \otimes \id \otimes \ldots \id
\]
in the basis $\{w_{i_1} \wedge \ldots \wedge w_{i_k}\}$, we will use (when the orientations are \OrUU) the $4 \times 4$ ``local matrix''
\[ 
\kbordermatrix{ 
& 1 & w_i & w_{i+1} & w_i \wedge w_{i+1} \\ 
1 & q & 0 & 0 & 0 \\ 
w_i & 0 & 0 & 1 & 0 \\ 
w_{i+1} & 0 & 1 & q-q^{-1} & 0 \\
w_i \wedge w_{i+1} & 0 & 0 & 0 & -q^{-1}
}
\]
to represent $\check{R}_{i,i+1}$. For general orientations, the local matrix representing $\check{R}_{i,i+1}$ is the appropriate $4 \times 4$ matrix from Corollary \ref{VVstarRmatrices} with the rows and columns labeled $1$, $w_i$, $w_{i+1}$, and $w_i \wedge w_{i+1}$ in the same order as above.

Because $\check{R}_{i,i+1}$ acts on $V^{\otimes \Sc'}$ as a tensor product of maps, it sends
\[
w_j \wedge w_k \mapsto \check{R}_{i,i+1}(w_j) \wedge \check{R}_{i,i+1}(w_k),
\]
as long as $\{j,k\} \neq \{i,i+1\}$. In other words, the action of $\check{R}_{i,i+1}$ on any ``global'' basis element $w_{i_1} \wedge \ldots \wedge w_{i_k}$ is calculated by extracting the local part, one of the four elements $\{1,w_i,w_{i+1},w_i \wedge w_{i+1}\}$. The local part is transformed by the local matrix, while the rest remains unchanged. The maps for minima and maxima are also tensor products of identity maps and a local piece, so they work the same way. Since $\check{R}$ and the maps for minima and maxima are even operators, no extra minus signs are needed. 

\begin{remark}
The definition of $1$, $w_i$, $w_{i+1}$, and $w_{i} \wedge w_{i+1}$ as specific tensors of $v$ generators depends on the orientations of the strands. However, the ordering $(1, w_i, w_{i+1}, w_i \wedge w_{i+1})$ of the columns and rows of each of the matrices in Corollary \ref{VVstarRmatrices} remains the same for all orientations.
\end{remark}

\begin{remark}
As shown by Kauffman and Saleur in \cite[Section 5]{KSFF}, the maps $\check{R}_{i,i+1}$ from $V^{\otimes n}$ to itself can be obtained, up to changes of normalization, as exterior powers of maps appearing in the Burau representation of the braid group.
\end{remark}

\subsubsection{Modified bases}

We now define the right and left modified bases. First, we define special elements of $V^{\otimes \Sc}$.
\begin{definition}\label{UqSpecialEltsDef}
Let $\Sc \subset [1,\ldots,n]$ and let $1 \leq i \leq n-1$. 
\begin{itemize}
\item
If $i, i+1 \in \Sc$, define
\[
l_i := w_i - q^{-1} w_{i+1}.
\]
\item
If $i \in \Sc$ but $i+1 \notin \Sc$, define
\[
l_i = w_i + w_{i+1}.
\]
\item
If $i \notin \Sc$ but $i+1 \in \Sc$, define
\[
l_i = q^{-1} w_i - q^{-1} w_{i+1}.
\]
\item
If $i, i+1 \notin \Sc$, define
\[
l_i = q^{-1} w_i + w_{i+1}.
\]
\end{itemize}

Let
\begin{equation}\label{L0SpecialCase}
l_0 := 
\begin{cases}
-q^{-1} w_1 & 1 \in \Sc \\
w_1 & 1 \notin \Sc
\end{cases}
\end{equation}
and
\begin{equation}\label{LNSpecialCase}
l_n :=
\begin{cases}
w_n & n \in \Sc \\
q^{-1} w_n & n \notin \Sc.
\end{cases}
\end{equation}
\end{definition}
Similar elements, in a slightly different tensor-product representation of $\Uq$, are defined in Ellis--Petkova--V{\'e}rtesi \cite[Section 2.3.3]{EPV}. 

\begin{definition}\label{ModifiedBasisDef} The right modified basis for $V^{\otimes \Sc}$ is
\[
\{ l_{i_1} \wedge \ldots \wedge l_{i_k} \}, \,\,\, 1 \leq i_i, \ldots, i_k \leq n.
\]
The left modified basis for $V^{\otimes \Sc}$ is
\[
\{ l_{i_1} \wedge \ldots \wedge l_{i_k} \}, \,\,\, 0 \leq i_i, \ldots, i_k \leq n-1.
\]
We write
\[
l_{\x} := l_{i_1} \wedge \ldots \wedge l_{i_k}, \quad \x = \{i_1,\ldots,i_k\} \subset [1,\ldots,n] \textrm{ or } [0,\ldots,n-1].
\]
\end{definition}

\begin{remark}\label{CanonicalBasesRem}
In general, we refrain from calling these bases ``canonical'' bases for $V^{\otimes \Sc}$. In the special case where all strands point up, so $V^{\otimes \Sc} = V^{\otimes n}$, Zhang \cite{ZhangCanonical} defines canonical bases which are used by Sartori in \cite{Sartori} to guide his categorification of $V^{\otimes n}$. One can see from \cite[Proposition 5.5]{Sartori}, plus a computation, that our right modified basis for $V^{\otimes n}$ may be transformed into the canonical basis by:
\begin{itemize}
\item replacing $v_0$ with $v_1$ and $v_1$ with $v_0$ everywhere,
\item mirroring the order of the tensor factors from left to right, so that $v_{i_1} \otimes \cdots \otimes v_{i_n}$ is replaced with $v_{i_n} \otimes \cdots \otimes v_{i_1}$, and
\item replacing $q$ with $-q^{-1}$.
\end{itemize}

\end{remark}

\subsection{Maps in the modified basis}\label{MapsInModifiedBasisSect}
We want to express the local matrices for $\check{R}_{i,i+1}$ and the maximum and minimum maps in the right and left modified bases.
\subsubsection{Generic case}
First, we consider the generic case. For a crossing between strands $i$ and $i+1$, with $n$ total strands, the generic case is $1 < i < n-1$. 

We use the same shorthand as above. The ``local part'' of $l_{i_1} \wedge \ldots \wedge l_{i_k}$ is an element of the set 
\[
\{1, \,\,\, l_{i-1}, \,\,\, l_i, \,\,\, l_{i+1}, \,\,\, l_{i-1} \wedge l_i, \,\,\, l_{i-1} \wedge l_{i+1}, \,\,\, l_{i-1} \wedge l_i \wedge l_{i+1}\}.
\]
There are twice as many possibilities as before; since $l_{i-1}$ is a linear combination of $w_{i-1}$ and $w_i$, and $w_i$ is a local generator, we must include $l_{i-1}$ in the list of local basis elements.

\begin{proposition}\label{DualCanonicalRForm}
Let $n \geq 3$ and $1 < i < n-1$. The local matrix for $\check{R}_{i,i+1}$ in the right or left modified basis is given as follows, where $\Gamma$ is the tangle corresponding to $\check{R}_{i,i+1}$:
\begin{itemize}
\item
If $\Gamma$ is oriented \OrUU, then $\check{R}_{i,i+1} =$
\[
\kbordermatrix{
& 1 & l_{i-1} & l_i & l_{i+1} & l_{i-1} \wedge l_i & l_{i-1} \wedge l_{i+1} & l_i \wedge l_{i+1} & l_{i-1} \wedge l_i \wedge l_{i+1} \\
1 & q &  &  &  &  &  &  &  \\
l_{i-1} &  & q & 0 & 0 &  &  &  &  \\
l_i &  & 1 & -q^{-1} & 1 &  &  &  &  \\
l_{i+1} &  & 0 & 0 & q &  &  &  &  \\
l_{i-1} \wedge l_i &  &  &  &  & -q^{-1} & 1 & 0 &  \\
l_{i-1} \wedge l_{i+1} &  &  &  &  & 0 & q & 0 &  \\
l_i \wedge l_{i+1} &  &  &  &  & 0 & 1 & -q^{-1} &  \\
l_{i-1} \wedge l_i \wedge l_{i+1} &  &  &  &  &  &  &  & -q^{-1}
}.
\]
\item
If $\Gamma$ is oriented \OrDU, then $\check{R}_{i,i+1} =$
\[
\kbordermatrix{
& 1 & l_{i-1} & l_i & l_{i+1} & l_{i-1} \wedge l_i & l_{i-1} \wedge l_{i+1} & l_i \wedge l_{i+1} & l_{i-1} \wedge l_i \wedge l_{i+1} \\
1 & 1 &  &  &  &  &  &  &  \\
l_{i-1} &  & 1 & 0 & 0 &  &  &  &  \\
l_i &  & q^{-1} & 1 & -q &  &  &  &  \\
l_{i+1} &  & 0 & 0 & 1 &  &  &  &  \\
l_{i-1} \wedge l_i &  &  &  &  & 1 & -q & 0 &  \\
l_{i-1} \wedge l_{i+1} &  &  &  &  & 0 & 1 & 0 &  \\
l_i \wedge l_{i+1} &  &  &  &  & 0 & q^{-1} & 1 &  \\
l_{i-1} \wedge l_i \wedge l_{i+1} &  &  &  &  &  &  &  & 1
}.
\]
\item
If $\Gamma$ is oriented \OrUD, then $\check{R}_{i,i+1} =$
\[
\kbordermatrix{
& 1 & l_{i-1} & l_i & l_{i+1} & l_{i-1} \wedge l_i & l_{i-1} \wedge l_{i+1} & l_i \wedge l_{i+1} & l_{i-1} \wedge l_i \wedge l_{i+1} \\
1 & 1 &  &  &  &  &  &  &  \\
l_{i-1} &  & 1 & 0 & 0 &  &  &  &  \\
l_i &  & -q & 1 & q^{-1} &  &  &  &  \\
l_{i+1} &  & 0 & 0 & 1 &  &  &  &  \\
l_{i-1} \wedge l_i &  &  &  &  & 1 & q^{-1} & 0 &  \\
l_{i-1} \wedge l_{i+1} &  &  &  &  & 0 & 1 & 0 &  \\
l_i \wedge l_{i+1} &  &  &  &  & 0 & -q & 1 &  \\
l_{i-1} \wedge l_i \wedge l_{i+1} &  &  &  &  &  &  &  & 1
}.
\]
\item
If $\Gamma$ is oriented \OrDD, then $\check{R}_{i,i+1} =$
\[
\kbordermatrix{
& 1 & l_{i-1} & l_i & l_{i+1} & l_{i-1} \wedge l_i & l_{i-1} \wedge l_{i+1} & l_i \wedge l_{i+1} & l_{i-1} \wedge l_i \wedge l_{i+1} \\
1 & -q^{-1} &  &  &  &  &  &  &  \\
l_{i-1} &  & -q^{-1} & 0 & 0 &  &  &  &  \\
l_i &  & 1 & q & 1 &  &  &  &  \\
l_{i+1} &  & 0 & 0 & -q^{-1} &  &  &  &  \\
l_{i-1} \wedge l_i &  &  &  &  & q & 1 & 0 &  \\
l_{i-1} \wedge l_{i+1} &  &  &  &  & 0 & -q^{-1} & 0 &  \\
l_i \wedge l_{i+1} &  &  &  &  & 0 & 1 & q &  \\
l_{i-1} \wedge l_i \wedge l_{i+1} &  &  &  &  &  &  &  & q
}.
\]
\end{itemize}

\end{proposition}

\begin{proof} The proof amounts to writing things out by hand using the dual-basis matrices in Corollary \ref{VVstarRmatrices}. We will show the details for one matrix entry and one choice of orientations on strands $i-1$ and $i+2$.

Assume that $\Gamma$ is oriented \OrUUUU. To compute $\check{R}_{i,i+1}(l_{i-1} \wedge l_i)$, we first write
\begin{align*}
l_{i-1} \wedge l_i &= (w_{i-1} - q^{-1}w_i) \wedge (w_i - q^{-1} w_{i+1}) \\
&= w_{i-1} \wedge w_i - q^{-1} w_{i-1} \wedge w_{i+1} + q^{-2} w_i \wedge w_{i+1}.
\end{align*}

Now, $w_{i-1} \wedge w_i$ has local form $w_i$, so 
\[
\check{R}_{i,i+1}(w_{i-1} \wedge w_i) = w_{i-1} \wedge w_{i+1}.
\]
Similarly, $w_{i-1} \wedge w_{i+1}$ has local form $w_{i+1}$, so
\[
\check{R}_{i,i+1}(w_{i-1} \wedge w_{i+1}) = w_{i-1} \wedge w_i + (q-q^{-1}) w_{i-1} \wedge w_{i+1}.
\]
Finally,
\[
\check{R}_{i,i+1}(w_i \wedge w_{i+1}) = -q^{-1} w_i \wedge w_{i+1}.
\]
Thus,
\begin{align*}
\check{R}_{i,i+1}(l_{i-1} \wedge l_i) &= w_{i-1} \wedge w_{i+1} \\
& \qquad -q^{-1}(w_{i-1} \wedge w_i + (q-q^{-1}) w_{i-1} \wedge w_{i+1}) \\
& \qquad +q^{-2}(-q^{-1} w_i \wedge w_{i+1}) \\
&= -q^{-1}(w_{i-1} \wedge w_i - q^{-1} w_{i-1} \wedge w_{i+1} + q^{-2} w_i \wedge w_{i+1}) \\
&= -q^{-1} l_{i-1} \wedge l_i.
\end{align*}

The rest of the matrix entries can be computed similarly. There are many cases to check; besides the orientations on strands $i$ and $i+1$, the orientations on strands $i-1$ and $i+2$ must also be considered since these orientations affect the definition of $l_{i-1}$ and $l_{i+1}$. However, many entries in the matrices are zero automatically, making the computation somewhat more reasonable.
\end{proof}

\begin{proposition}\label{MinMaxModifiedBasis}
Let $n \geq 2$ and $1 \leq i \leq n-1$. The local matrices for the minimum and maximum maps in the right or left modified basis are given as follows:
\begin{itemize}
\item
A minimum with matched strands $\{i+1,i+2\}$ at top, oriented \MinRL or \MinLR, is assigned the map with matrix
\[
\kbordermatrix{
& 1 & l_i \\
1 & 0 & 0 \\
l_i & 0 & 0 \\
l_{i+1} & 1 & 0 \\
l_{i+2} & 0 & 0  \\
l_i \wedge l_{i+1} & 0 & 1 \\
l_i \wedge l_{i+2} & 0 & 0 \\
l_{i+1} \wedge l_{i+2} & 0 & 1 \\
l_i \wedge l_{i+1} \wedge l_{i+2} & 0 & 0 
}.
\]

\item
A maximum with matched strands $\{i+1,i+2\}$ on bottom, oriented \MaxRL or \MaxLR, is assigned the map with matrix
\[
\kbordermatrix{
& 1 & l_i & l_{i+1} & l_{i+2} & l_i \wedge l_{i+1} & l_i \wedge l_{i+2} & l_{i+1} \wedge l_{i+2} & l_i \wedge l_{i+1} \wedge l_{i+2} \\
1 & 0 & 1 & 0 & 1 & 0 & 0 & 0 & 0 \\
l_{i} & 0 & 0 & 0 & 0 & 0 & 1 & 0 & 0 \\
}.
\]

\end{itemize}
\end{proposition}

\begin{proof} The proof is another entry-by-entry computation, this time using the matrices from Section \ref{UqMinimaSect} and Section \ref{UqMaximaSect}. The orientations on the two strands adjacent to the critical point must be considered. The cases $i = 1$ and $i = n-1$ are not special, but we must also consider $i = 0$ and $i = n$ as special cases; see below.
\end{proof}

\subsubsection{Special cases}\label{FirstSpecialCasesSec}
First we consider a crossing between strands $i$ and $i+1$. When $i = 1$, there is no local generator $l_{i-1}$ in the right modified basis, whose list of local generators is
\[
\{1, \,\,\, l_1, \,\,\, l_2, \,\,\, l_1 \wedge l_2\}
\]
(unless $n = 1$ in which case the list is $\{1, l_1\}$, or $n = 0$ in which case the list is $\{1\}$, but there are no maps for crossings to consider in these cases). The left modified basis has the usual eight local generators if $n \geq 2$, with $l_0$ defined by formula (\ref{L0SpecialCase}). It has local generators $\{1,l_0\}$ if $n = 1$ and $\{1\}$ if $n = 0$.

\begin{proposition}\label{RightSpecialCaseRProp}
If $n \geq 3$, the local matrix for $\check{R}_{1,2}$ in the left modified basis is the appropriate $8 \times 8$ matrix from Proposition \ref{DualCanonicalRForm}, without modification.

If $n \geq 2$, the local matrix for $\check{R}_{1,2}$ in the right modified basis is a submatrix of this $8 \times 8$ matrix. One discards all columns and rows involving $l_0$ to obtain a $4 \times 4$ matrix representing $\check{R}_{1,2}$.
\end{proposition}

\begin{proof} 
For the left modified basis, this proposition is proved by checking matrix entries individually as in Proposition \ref{DualCanonicalRForm}. Some of these checks differ from computations in Proposition \ref{DualCanonicalRForm} because formula (\ref{L0SpecialCase}) is being used. However, the resulting $8 \times 8$ matrices are the same.

For the right modified basis, the computations are the same as in Proposition \ref{DualCanonicalRForm}; the absence of $l_{i-1} = l_0$ does not change the other values. The only exception is when $n = 2$, in which case we need to use formula (\ref{LNSpecialCase}) for $l_n = l_2$.
\end{proof}

Similarly, when $i = n-1$, the local generator $l_{i+1} = l_n$ of the right modified basis is defined by formula (\ref{LNSpecialCase}). When $n \geq 2$, the right modified basis has the usual eight local generators; when $n = 1$, the local generators are $\{1,l_1\}$, and when $n = 0$ the local generators are $\{1\}$.

There is no local generator $l_n$ in the left modified basis, so the list of local generators for $i = n-1$ is
\[
\{1, \,\,\, l_{n-2}, \,\,\, l_{n-1}, \,\,\, l_{n-2} \wedge l_{n-1}\}
\]
as long as $n \geq 2$ (if $n = 1$ the list is $\{1,l_0\}$, and if $n = 0$ the list is $\{1\}$).

\begin{proposition}\label{LeftSpecialCaseRProp}
If $n \geq 3$, the local matrix for $\check{R}_{n-1,n}$ in the right modified basis is the appropriate $8 \times 8$ matrix from Proposition \ref{DualCanonicalRForm}, without modification.

If $n \geq 2$, the local matrix for $\check{R}_{n-1,n}$ in the left modified basis is a submatrix of this $8 \times 8$ matrix. One discards all columns and rows involving $l_n$ to obtain a $4 \times 4$ matrix representing $\check{R}_{n-1,n}$.
\end{proposition}

\begin{proof}
The proof is another series of matrix entry checks as described in Proposition \ref{RightSpecialCaseRProp}, using formula (\ref{LNSpecialCase}) for $l_n$. When $n = 2$, we must also use formula (\ref{L0SpecialCase}) for $l_{n-2} = l_0$.
\end{proof}

Now we consider the minimum and maximum maps with strands $i+1$ and $i+2$ matched. The cases $i = 1$ and $i = n-1$ are not exceptional; however, we do need to consider $i = 0$ and $i = n$ as special cases.

\begin{proposition}
Let $\Gamma$ be an $(n,n+2)$--tangle or $(n+2,n)$--tangle representing a minimum or maximum with strands $1$ and $2$ matched. For $n \geq 1$, the local matrix in the left modified basis for the $\Uq$--linear map $V^{\otimes \Sc'_{\Gamma}} \to V^{\otimes \Sc_{\Gamma}}$ associated to $\Gamma$ is the appropriate $8 \times 2$ or $2 \times 8$ matrix from Proposition \ref{MinMaxModifiedBasis}, without modification.

For $n \geq 0$, the local matrix in the right modified basis for this map is a $4 \times 1$ or $1 \times 4$ submatrix of the $8 \times 2$ or $2 \times 8$ matrix, obtained by by discarding all columns and rows involving $l_0$.
\end{proposition}

\begin{proof}
The computations are the same as in Proposition~\ref{MinMaxModifiedBasis}, except that formula (\ref{L0SpecialCase}) is used in the left modified basis. When $n = 0$, we must also use formula (\ref{LNSpecialCase}) for the local generator $l_{n+2} = l_2$ of the right modified basis.
\end{proof}

\begin{proposition}\label{RSpecialCaseMinMaxProp}
Let $\Gamma$ be an $(n,n+2)$--tangle consisting of a single minimum point between strands $n+1$ and $n+2$. For $n \geq 1$, the local matrix in the right modified basis for the $\Uq$--linear map $V^{\otimes \Sc'_{\Gamma}} \to V^{\otimes \Sc_{\Gamma}}$ associated to $\Gamma$ is the appropriate $8 \times 2$ matrix from Proposition \ref{MinMaxModifiedBasis}, without modification.

For $n \geq 0$, the local matrix for this map in the left modified basis is a $4 \times 1$ submatrix of the $8 \times 2$ matrix, obtained by discarding all columns involving $l_n$ and rows involving $l_{n+2}$.

Similarly, let $\Gamma$ be an $(n+2,n)$--tangle consisting of a single maximum point between strands $n+1$ and $n+2$. For $n \geq 1$, the local matrix in the right modified basis for the $\Uq$-linear map associated to $\Gamma$ is the appropriate $2 \times 8$ matrix from Proposition~\ref{MinMaxModifiedBasis}, without modification.

For $n \geq 0$, the local matrix for this map in the left modified basis is a $1 \times 4$ submatrix of the $2 \times 8$ matrix, obtained by discarding all columns involving $l_{n+2}$ and rows involving $l_n$. 
\end{proposition}

\begin{proof}
Again, the computations are the same as in Proposition~\ref{MinMaxModifiedBasis}, except that formula (\ref{LNSpecialCase}) is used in the right modified basis. When $n = 0$, formula (\ref{L0SpecialCase}) is used for the local generator $l_0 = l_n$ of the left modified basis.
\end{proof}

\section{Viro's quantum relative \texorpdfstring{$\A^1$}{A\textasciicircum{1}}}\label{ViroSect}

\subsection{Palettes and colorings}\label{ViroPaletteSect}

\begin{definition}[cf. Section 2.8 of \cite{Viro}]\label{PaletteDef}
A $1$--palette is a quadruple 
\[
P = (B,M,V, \phi: M \times V \to M)
\]
where $B$ is a commutative ring, $M$ is a subgroup of the group of multiplicative units of $B$, $V$ is a subgroup of the underlying additive group of $B$, and $\phi: M \times V \to M$ is a bilinear map (i.e. both multiplication in $M$ and addition in $V$ get sent to multiplication in $M$). 
\end{definition}

Given $P = (B,M,V, \phi: M \times V \to M)$, let $\M_B$ be the category of finite-dimensional modules over $B$ and $B$--linear maps between them. 

\begin{remark} If we wanted, we could still require the objects of $\M_B$ to be supermodules, which come with a decomposition as a direct sum of an even part and an odd part. The whole commutative ring $B$ is viewed as even, so supermodule structures on $B$--modules are relatively easy to arrange (any decomposition of a module as a direct sum will work). 

As noted in \cite[Section 2.10]{Viro}, one can upgrade the target category of the functor $\A^1_P: \G_P \to \M_B$ (defined below) to a category of representations of a Hopf superalgebra $U^1 \otimes_{\Z} B$, where $U^1$ has both odd and even generators. In this case, as in the singly graded case, it is important to work with supermodules rather than just ordinary modules. However, since we will not discuss $U^1$ here, we will not worry about supermodule structures.
\end{remark}

\begin{definition}\label{UnivPaletteObjs}
Let $Y$ be an oriented zero-manifold consisting of $n$ points. Let $M_Y := H^0(Y; \frac{1}{4}\Z)$, viewed multiplicatively. Let $B_Y := \Z[M_Y]$ and $V := \Z$. The map $\phi$ from $M_Y \times V$ to $M_Y$ is defined to be exponentiation. The $1$--palette $(B_Y,M_Y,V,\phi)$ will be denoted $P_Y$. 
\end{definition}

\begin{definition}\label{UnivPaletteGraphs}
Let $W$ be an compact oriented $1$--manifold. Let 
\[
M_W := H^1\bigg(W,\partial W; \frac{1}{4}\Z\bigg),
\]
viewed multiplicatively. Let $B_W := \Z[M]$ and let $V := \Z$. The map $\phi$ from $M_W \times V$ to $M_W$ is again defined to be exponentiation. The $1$--palette $(B_W,M_W,V,\phi)$ will be denoted $P_W$. 
\end{definition}

Let $\Gamma$ be a tangle with underlying $1$--manifold $W$, such that $\partial W = -Y_2 \sqcup Y_1$. We define colorings of $Y_1$ and $Y_2$ by the $1$--palettes $P_{Y_1}$ and $P_{Y_2}$ of Definition \ref{UnivPaletteObjs}, and we define a coloring of $\Gamma$ by the $1$--palette $P_W$ of Definition \ref{UnivPaletteGraphs}:
\begin{definition}\label{UnivColorObjects}
For each point $p_i$ of $Y_1$, we have an element $\pi_i$ of $H^0(Y_1;\frac{1}{4}\Z)$ such that $\pi_i([p_j]) = \delta_{i,j}$, where $[p_j] \in H_0(Y_1; \frac{1}{4}\Z)$ is the fundamental class of $p_j$. Color $p_i$ as $(\pi_i^{1/4}, 0)$; this assignment defines a $P_{Y_1}$--coloring on $Y_1$. We can define a $P_{Y_2}$--coloring on $Y_2$ in the same way.
\end{definition}

\begin{definition}\label{UnivColorGraphs}
For each strand $s_i$ of $\Gamma$, we have an element $\sigma_i$ of $H^1(W,\partial W; \frac{1}{4}\Z)$ such that $\sigma_i([s_j]) = \delta_{i,j}$, where $[s_j] \in H_1(W,\partial W; \frac{1}{4}\Z)$ is the fundamental class of $s_j$. Color $s_i$ as $(\sigma_i^{1/4}, 0)$. For any choice of framings on $\Gamma$, this assignment defines a $P_W$--coloring on $\Gamma$.
\end{definition}

\begin{remark}\label{UniversalPaletteRemark} $P_W$ is a relative version, allowing non-closed graphs (but no trivalent vertices), of the ``universal'' palette discussed in \cite[Section 7.4]{Viro}. Viro also takes coefficients in $\Z$, rather than $\frac{1}{4}\Z$, and he takes the quotient field of $B_W$ for convenience in discussing skein relations with trivalent vertices. The coloring of $\Gamma$ defined above for a given choice of orientations and framings corresponds to Viro's ``universal coloring'' of $\Gamma$ with $T$--components chosen to be zero. 
\end{remark}

The coboundary map $d^*: H^0(\partial W;\frac{1}{4}\Z) \to H^1(W,\partial W; \frac{1}{4}\Z)$ gives us ring homomorphisms 
\begin{gather*}
B_{Y_1} = \Z \bigg[ H^0 \bigg( Y_1;\frac{1}{4}\Z \bigg) \bigg] \xrightarrow{d^*} \Z \bigg[ H^1 \bigg( W,\partial W; \frac{1}{4}\Z \bigg) \bigg] = B_W \\
\tag*{\text{and}} B_{Y_2} = \Z \bigg[ H^0 \bigg( Y_2;\frac{1}{4}\Z \bigg) \bigg] \xrightarrow{-d^*} \Z \bigg[ H^1 \bigg( W,\partial W; \frac{1}{4}\Z \bigg) \bigg] = B_W.
\end{gather*}
These homomorphisms satisfy the requirements for functorial change of colors as discussed in \cite[Section 7.3]{Viro}. Thus, we get $P_W$--colorings from the $P_{Y_i}$--colorings on $Y_1$ and $Y_2$ by applying the homomorphisms to the components $t$ and $T$ of the colors $(t,T)$ of the points in $Y_i$. These colorings agree with the ones obtained by restricting the $P_W$--coloring of $\Gamma$ to its boundary $-Y_2 \sqcup Y_1$. By the definition of $\A^1_P$ below, we will have
\[
\A^1_{P_W}(Y_i) \cong \A^1_{P_{Y_i}}(Y_i) \otimes_{\Z[M_{Y_i}]} \Z[M_W] 
\]
as $\Z[M_W]$--modules, where $\Z[M_{Y_i}]$ acts on $\Z[M_W]$ via the above homomorphisms. 

\subsection{Generic graphs and the functor \texorpdfstring{$\A^1$}{A\textasciicircum{1}}}\label{ViroA1DefSect}

Let $P$ be a $1$--palette. In \cite[Section 2.8]{Viro}, Viro defines a category $\G^1_P$ of framed generic graphs (i.e. framed oriented graphs with vertices of valence $1$ and $3$, generically embedded in $\R^2 \times [0,1]$, whose intersection with $\R^2 \times \{0,1\}$ is their set of $1$--valent vertices) that are colored by $P$ as defined below. Objects of $\G^1_P$ are finite collections $Y$ of oriented points arranged along a line in $\R^2$, each colored with a triple $(t,T,\epsilon)$ where $t \in M$ and $T \in V$. The $\epsilon$--component of the color simply encodes the orientation of the point, so we will omit it and use pairs $(t,T)$ instead. The elements $t$ of $M$ are required to satisfy $t^4 \neq 1$; since we will only consider the case where $M$ is a free abelian group, this condition will not arise in this paper. 

Morphisms of $\G^1_P$ are colored framed generic oriented graphs $\Gamma$, viewed as morphisms from their bottom endpoints to their top endpoints. They may have trivalent vertices as well as crossings between strands; we will consider only crossings and not trivalent vertices. The orientation of a graph $\Gamma$ must be compatible with the orientations of its endpoints, viewed as objects of $\G^1_P$: bottom endpoints of $\Gamma$ must have the opposite of the boundary-induced orientation and top endpoints must have the boundary-induced orientation. 

The morphisms $\Gamma$ of $\G^1_P$ are colored in the sense that every edge is labeled with some $(t,T)$ just like the vertices in the objects of $\G^1_P$. If an edge is adjacent to a vertex, it must have the same label as the vertex. There are admissibility conditions on the colorings around trivalent vertices which we can ignore here. See \cite{Viro} for more details.

Viro defines a functor
\[
\A^1_P: \G^1_P \to \M_B.
\] 
As we will see below, an oriented zero-manifold $Y$ consisting of $n$ points along a line in $\R^2$ is assigned a free $B$--module $\A^1_P(Y)$ of rank $2^n$. A graph $\Gamma \subset \R^2 \times [0,1]$ whose underlying oriented $1$--manifold is $W$, with $\partial W = -Y_2 \sqcup Y_1$, is assigned a $B$--linear map
\[
\A^1_P(\Gamma): \A^1_P(Y_1) \to \A^1_P(Y_2)
\]
which depends on the orientations, framings, and colorings $(t,T)$ of $\Gamma$. However, changing the framings or the $T$--component of the color on any edge only affects the associated morphism under $\A^1_P$ by scalar multiplication. When $T = 0$, the morphisms are framing-independent. 

\subsection{Assigning modules to objects}\label{ViroObjectInvsSect}
In this subsection, we will discuss the values $\A^1_P(Y)$ of $\A^1_P$ on objects of $\G^1_P$. Objects $Y$ of $\G_P$ are collections of $n$ oriented points along a line, for some $n \geq 0$, each colored with $(t_i,T_i,\epsilon_i)$ for some $t_i \in M$ and $T_i \in V$. The sign $\epsilon_i$ (plus or minus) of the color must agree with the orientation of the point, so we will omit it from the notation. 

Let $\Lambda_B$ be a free supermodule over $B$ of rank $(1|1)$ with one even generator $v_0$ and one odd generator $v_1$. The value of $\A^1_P$ on an object $Y$ of $\G_P$ is defined to be $\Lambda_B^{\otimes n}$, where the tensor products are taken over $B$. Ignoring the supermodule structure, $\Lambda_B^{\otimes n}$ is a free $B$--module of rank $2^n$. 

Suppose $P = P_Y$ is the palette of Definition \ref{UnivPaletteObjs}, where $Y$ consists of $n$ oriented points along a line. Give $Y$ the coloring of Definition \ref{UnivColorObjects}. We have $B_Y = \Z[H^0(Y;\frac{1}{4}\Z)]$, so
\begin{equation}\label{ObjPaletteEqn}
\A^1_{P_Y}(Y) = \Lambda_{B_Y}^{\otimes n} = (B_Y \oplus B_Y)^{\otimes n} = B_Y^{\oplus (2^n)} = \bigg( \Z \bigg[ H^0 \bigg(Y; \frac{1}{4}\Z \bigg) \bigg] \bigg)^{\oplus (2^n)}.
\end{equation}

Now suppose $P = P_W$ is the palette of Definition \ref{UnivPaletteGraphs} for an oriented $1$--manifold $W$, where $Y$ is one of the boundary components of $W$. Give $Y$ the $P_W$--coloring of Definition \ref{UnivColorGraphs} restricted to the boundary of $W$, or equivalently the coloring induced by the functorial change of colors from $P_Y$ to $P_W$, as discussed in Section \ref{ViroPaletteSect}. We have $B_W = \Z[H^1(W,\partial W; \frac{1}{4}\Z)]$, so
\begin{equation}\label{WPaletteEqn}
\A^1_{P_W}(Y) = \Lambda_{B_W}^{\otimes n} = (B_W \oplus B_W)^{\otimes n} = B_W^{\oplus (2^n)} = \bigg( \Z \bigg[ H^1 \bigg( W,\partial W; \frac{1}{4}\Z \bigg) \bigg] \bigg)^{\oplus (2^n)}.
\end{equation}

Comparing (\ref{ObjPaletteEqn}) and (\ref{WPaletteEqn}), we see that
\begin{equation}\label{ChangingPalettesIsTensorObjs}
\A^1_{P_W}(Y) =  \A^1_{P_Y}(Y) \otimes_{\Z[H^0(Y;\frac{1}{4}\Z)]} \Z \bigg[ H^1 \bigg(W, \partial W; \frac{1}{4}\Z \bigg) \bigg].
\end{equation}

\subsection{Maps for crossings}\label{ViroMapsForCrossingsSect}
Let $\Gamma$ be an oriented $(2,2)$--tangle consisting of a single crossing (positive as a braid generator, without regard for orientations). Let $W$ be the oriented $1$--manifold underlying $\Gamma$; then $W$ is a disjoint union of two arcs. Label these arcs as strand $1$ and strand $2$; strand $1$ is the one whose bottom endpoint is leftmost. 

Give $\Gamma$ any framing. As in Definition \ref{UnivColorGraphs} above, color strand $i$ of $\Gamma$ as $(t_i,0)$ where 
\begin{equation}\label{TiOnStrandsEqn}
t_i := \sigma_i^{1/4} \in H^1\bigg(W,\partial W; \frac{1}{4}\Z\bigg).
\end{equation}
With this choice of coloring on $\Gamma$, we compute 
\[
\A^1_{P_W}(\Gamma): \Lambda_{B_W} \otimes_{B_W} \Lambda_{B_W} \to \Lambda_{B_W} \otimes_{B_W} \Lambda_{B_W}.
\]
Let $v_0 \otimes v_0$, $v_1 \otimes v_0$, $v_0 \otimes v_1$, and $v_1 \otimes v_1$ denote the four natural generators of $\Lambda_{B_W} \otimes_{B_W} \Lambda_{B_W}$. In Viro's terminology, $v_0$ is the bosonic generator and $v_1$ is the fermionic generator. We may represent $\A^1_{P_W}(\Gamma)$ as a matrix in this basis; the matrix entries are elements of $B_W$. 

These entries are the Boltzmann weights defined in \cite[Table 3]{Viro}; we may read off the matrix for $\A^1_{P_W}(\Gamma)$ from this table. One should replace Viro's $(t,T)$ with $(t_1,T_1)$ and $(u,U)$ with $(t_2,T_2)$, where $t_1$ and $t_2$ are defined by (\ref{TiOnStrandsEqn}) above. Since we choose $T_1 = T_2 = 0$, the factor $t_1^{-T_2} t_2^{-T_1}$ appearing in Viro's Boltzmann weights is equal to $1$ and we may ignore it.

\begin{proposition}\label{BoltzmannWeightsProp} 
If $\Gamma$ is oriented \OrUU, then $\A^1_{P_W}(\Gamma)$ has matrix
\[
\A^1_{P_W}(\Gamma) = 
\kbordermatrix{ 
& v_0 \otimes v_0 & v_1 \otimes v_0 & v_0 \otimes v_1 & v_1 \otimes v_1 \\ 
v_0 \otimes v_0 & t_1 t_2 & 0 & 0 & 0 \\ 
v_1 \otimes v_0 & 0 & 0 & t_1^{-1} t_2 & 0 \\ 
v_0 \otimes v_1 & 0 & t_1 t_2^{-1}  & (t_1^4 - 1)t_1^{-1}t_2^{-1} & 0 \\
v_1 \otimes v_1 & 0 & 0 & 0 & -t_1^{-1}t_2^{-1}
}.
\]

If $\Gamma$ is oriented \OrDU, then $\A^1_{P_W}(\Gamma)$ has matrix
\[
\A^1_{P_W}(\Gamma) = 
\kbordermatrix{ 
& v_1 \otimes v_0 & v_0 \otimes v_0 & v_1 \otimes v_1 & v_0 \otimes v_1 \\ 
v_0 \otimes v_1 & t_1^{-1} t_2 & 0 & 0 & 0 \\ 
v_1 \otimes v_1 & 0 & 0 & -t_1 t_2 & 0 \\ 
v_0 \otimes v_0 & 0 & t_1^{-1} t_2^{-1} & (1-t_1^{-4})t_1 t_2^{-1} & 0 \\
v_1 \otimes v_0 & 0 & 0 & 0 & t_1 t_2^{-1}
}.
\]

If $\Gamma$ is oriented \OrUD, then $\A^1_{P_W}(\Gamma)$ has matrix
\[
\A^1_{P_W}(\Gamma) = 
\kbordermatrix{ 
& v_0 \otimes v_1 & v_1 \otimes v_1 & v_0 \otimes v_0 & v_1 \otimes v_0 \\ 
v_1 \otimes v_0 & t_1 t_2^{-1} & 0 & 0 & 0 \\ 
v_0 \otimes v_0 & 0 & 0 & t_1^{-1} t_2^{-1} & 0 \\ 
v_1 \otimes v_1 & 0 & -t_1 t_2  & (1 - t_1^4)t_1^{-1}t_2 & 0 \\
v_0 \otimes v_1 & 0 & 0 & 0 & t_1^{-1} t_2
}.
\]

Finally, if $\Gamma$ is oriented \OrDD, then $\A^1_{P_W}(\Gamma)$ has matrix
\[
\A^1_{P_W}(\Gamma) = 
\kbordermatrix{ 
& v_1 \otimes v_1 & v_0 \otimes v_1 & v_1 \otimes v_0 & v_0 \otimes v_0 \\ 
v_1 \otimes v_1 & -t_1^{-1} t_2^{-1} & 0 & 0 & 0 \\ 
v_0 \otimes v_1 & 0 & 0 & t_1 t_2^{-1} & 0 \\ 
v_1 \otimes v_0 & 0 & t_1^{-1} t_2  & (1 - t_1^{-4})t_1 t_2  & 0 \\
v_0 \otimes v_0 & 0 & 0 & 0 & t_1 t_2
}.
\]
\end{proposition}

\begin{proof}
The matrix entries are given in \cite[Table 3]{Viro}; see also Remark \ref{ViroTable1Vs3Remark}. 
\end{proof}

\begin{remark}\label{ViroTable1Vs3Remark} The top-left and lower-right entries of the second and third matrices above are actually the inverses of the Boltzmann weights in \cite[Table 3]{Viro}. Comparing with \cite[Table 1]{Viro}, which is meant to reduce to Table 3 after substituting $t$ for $q^i$, $u$ for $q^j$, $T$ for $I$, and $U$ for $J$, we see that these Boltzmann weights are also inverted in Table 1 relative to Table 3. 
\end{remark}

For $1 \leq i \leq n$, define $t_i \in H^0(Y; \frac{1}{4}\Z)$ by 
\begin{equation}\label{TiOnPointsEqn}
t_i = \pi_i^{1/4},
\end{equation}
where $\pi_i$ is the cohomology class from Definition \ref{UnivColorObjects}. We may change to a dual basis for $\Lambda_{B_W}$ on each downward-pointing strand (say strand $i$) by defining $v_0^* = v_0$ and $v_1^* = -t_i^2 v_1$; see \cite[11.8.A]{Viro}. 

\begin{warning}\label{FirstTiWarning}
In (\ref{TiOnStrandsEqn}), we defined elements $t_i := \sigma_i^{1/4} \in H^1(W,\partial W; \frac{1}{4}\Z)$, where $i$ indexed the connected components of $W$ (i.e. the strands of $\Gamma$). In (\ref{TiOnPointsEqn}), we have $t_i := \pi_i^{1/4}$ of $H^0(Y;\frac{1}{4}\Z)$, where $i$ indexes the points of $Y$. We find this notation to be the most convenient, at the expense of possible confusion about the meaning of $t_i$. In the typical situation of Section \ref{ViroMapsModifiedBasisSect}, one will have $t_i \in H^0(Y;\frac{1}{4}\Z)$ as a coefficient of a modified basis element, and one will want to obtain a corresponding element in $H^1(W,\partial W;\frac{1}{4}\Z)$. The corresponding element, obtained either through functorial change of colors or through restriction, will always be 
\[
t_{\kappa(i)} \in H^1\bigg(W, \partial W; \frac{1}{4}\Z\bigg),
\]
where $\kappa(i)$ is the index of the component of $W$ containing the point with index $i$. When this component is the unique component of $W$ containing a minimum or maximum point, we will write $t$ in place of $t_{\kappa(i)}$.
\end{warning}

We may write dual-basis matrices for $\A^1_{P_W}(\Gamma)$:
\begin{corollary}\label{ViroCrossingsDualBases}
If $\Gamma$ is oriented \OrUU, then $\A^1_{P_W}(\Gamma)$ has dual-basis matrix
\[
\A^1_{P_W}(\Gamma) = 
\kbordermatrix{ 
& v_0 \otimes v_0 & v_1 \otimes v_0 & v_0 \otimes v_1 & v_1 \otimes v_1 \\ 
v_0 \otimes v_0 & t_1 t_2 & 0 & 0 & 0 \\ 
v_1 \otimes v_0 & 0 & 0 & t_1^{-1} t_2 & 0 \\ 
v_0 \otimes v_1 & 0 & t_1 t_2^{-1}  & (t_1^4 - 1)t_1^{-1}t_2^{-1} & 0 \\
v_1 \otimes v_1 & 0 & 0 & 0 & -t_1^{-1}t_2^{-1}
}.
\]

If $\Gamma$ is oriented \OrDU, then $\A^1_{P_W}(\Gamma)$ has dual-basis matrix
\[
\A^1_{P_W}(\Gamma) = 
\kbordermatrix{ 
& v_1^* \otimes v_0 & v_0^* \otimes v_0 & v_1^* \otimes v_1 & v_0^* \otimes v_1 \\ 
v_0 \otimes v_1^* & t_1^{-1} t_2 & 0 & 0 & 0 \\ 
v_1 \otimes v_1^* & 0 & 0 & -t_1 t_2 & 0 \\ 
v_0 \otimes v_0^* & 0 & t_1^{-1} t_2^{-1} & (1 - t_1^4)t_1^{-1} t_2^{-1} & 0 \\
v_1 \otimes v_0^* & 0 & 0 & 0 & t_1 t_2^{-1}
}.
\]

If $\Gamma$ is oriented \OrUD, then $\A^1_{P_W}(\Gamma)$ has dual-basis matrix
\[
\A^1_{P_W}(\Gamma) = 
\kbordermatrix{ 
& v_0 \otimes v_1^* & v_1 \otimes v_1^* & v_0 \otimes v_0^* & v_1 \otimes v_0^* \\ 
v_1^* \otimes v_0 & t_1 t_2^{-1} & 0 & 0 & 0 \\ 
v_0^* \otimes v_0 & 0 & 0 & t_1^{-1} t_2^{-1} & 0 \\ 
v_1^* \otimes v_1 & 0 & -t_1 t_2  & (t_1^4 - 1)t_1^{-1}t_2^{-1} & 0 \\
v_0^* \otimes v_1 & 0 & 0 & 0 & t_1^{-1} t_2
}.
\]

Finally, if $\Gamma$ is oriented \OrDD, then $\A^1_{P_W}(\Gamma)$ has dual-basis matrix
\[
\A^1_{P_W}(\Gamma) = 
\kbordermatrix{ 
& v_1^* \otimes v_1^* & v_0^* \otimes v_1^* & v_1^* \otimes v_0^* & v_0^* \otimes v_0^* \\ 
v_1^* \otimes v_1^* & -t_1^{-1} t_2^{-1} & 0 & 0 & 0 \\ 
v_0^* \otimes v_1^* & 0 & 0 & t_1 t_2^{-1} & 0 \\ 
v_1^* \otimes v_0^* & 0 & t_1^{-1} t_2  & (t_1^4 - 1)t_1^{-1} t_2^{-1}  & 0 \\
v_0^* \otimes v_0^* & 0 & 0 & 0 & t_1 t_2
}.
\]
In each of the matrix entries, $t_i$ is defined as in (\ref{TiOnStrandsEqn}).
\end{corollary}

As with the maps $\check{R}_{i,i+1}$ above, if $\Gamma$ is an $(n,n)$--tangle with a crossing between strands $i$ and $i+1$, then
\[
\A^1_{P_W}(\Gamma): \Lambda_{B_W}^{\otimes n} \to \Lambda_{B_W}^{\otimes n}
\]
is the tensor product of the appropriate map from Proposition \ref{BoltzmannWeightsProp} on tensor factors $i$ and $i+1$ with the identity map on all other tensor factors. The strands are ordered so that strand $i$ has positive slope and strand $i+1$ has negative slope. 

Finally, if the crossing $\Gamma$ is negative as a braid generator, rather than positive, then we can compute $\A^1_{P_W}(\Gamma)$ as the inverse of one of the matrices from Proposition \ref{BoltzmannWeightsProp}. We will not write out these inverse matrices here.

\subsection{Maps for minima and maxima}\label{ViroMinMaxSect}
First we consider minima. Let $\Gamma$ be a $(0,2)$--tangle consisting of a single minimum point. In this case, the underlying $1$--manifold $W$ of $\Gamma$ is just a single oriented arc $s$. Let $\sigma \in H^1(W, \partial W; \frac{1}{4}\Z)$ satisfy $\sigma([s]) = 1$. The $P_W$--color of $s$ under Definition \ref{UnivColorGraphs} is $(t,0)$ where $t := \sigma^{1/4}$. 

We compute $\A^1_{P_W}(\Gamma)$ using the Boltzmann weights defined in \cite[Table 3]{Viro}, just like when $\Gamma$ is a crossing.

\begin{proposition}\label{ViroMinimaMapsProp}
If $\Gamma$ is oriented \MinRL, then $\A^1_{P_W}(\Gamma)$ has matrix
\[
\A^1_{P_W}(\Gamma) = 
\kbordermatrix{ 
& 1 \\ 
v_0 \otimes v_1 & 0  \\ 
v_1 \otimes v_1 & -t^2  \\ 
v_0 \otimes v_0 & 1  \\
v_1 \otimes v_0 & 0 
}
\]
and dual-basis matrix
\[
\A^1_{P_W}(\Gamma) = 
\kbordermatrix{ 
& 1 \\ 
v_0 \otimes v_1^* & 0  \\ 
v_1 \otimes v_1^* & 1  \\ 
v_0 \otimes v_0^* & 1  \\
v_1 \otimes v_0^* & 0 
}.
\]
If $\Gamma$ is oriented \MinLR, then $\A^1_{P_W}(\Gamma)$ has matrix
\[
\A^1_{P_W}(\Gamma) = 
\kbordermatrix{ 
& 1 \\ 
v_1 \otimes v_0 & 0  \\ 
v_0 \otimes v_0 & t^{-2}  \\ 
v_1 \otimes v_1 & 1  \\
v_0 \otimes v_1 & 0 
}
\]
and dual-basis matrix
\[
\A^1_{P_W}(\Gamma) = 
\kbordermatrix{ 
& 1 \\ 
v_1^* \otimes v_0 & 0  \\ 
v_0^* \otimes v_0 & t^{-2}  \\ 
v_1^* \otimes v_1 & -t^{-2}  \\
v_0^* \otimes v_1 & 0 
}.
\]
\end{proposition}

\begin{proof}
The matrix entries are given in \cite[Table 3]{Viro}.
\end{proof}

In general, if $\Gamma$ is an $(n,n+2)$--tangle consisting of a single minimum point between strands $i+1$ and $i+2$, then
\[
\A^1_{P_W}(\Gamma): \Lambda_{B_W}^{\otimes n} \to \Lambda_{B_W}^{\otimes (n+2)}
\]
is the tensor product of identity transformations on $\Lambda_{B_W}$ with the appropriate map from Proposition \ref{ViroMinimaMapsProp}, inserted to the right of $i$ of the identity transformations.

Now we consider maxima. Let $\Gamma$ be a $(2,0)$--tangle consisting of a single maximum point. Again, the $P_W$--color of the one arc $s$ of $\Gamma$ under Definition \ref{UnivColorGraphs} is $(t,0)$, where $t := \sigma^{1/4}$ and $\sigma([s]) = 1$.

\begin{proposition}\label{ViroMaximaMapsProp}
If $\Gamma$ is oriented \MaxRL, then $\A^1_{P_W}(\Gamma)$ has matrix
\[
\A^1_{P_W}(\Gamma) = 
\kbordermatrix{ 
& v_1 \otimes v_0 & v_0 \otimes v_0 & v_1 \otimes v_1 & v_0 \otimes v_1 \\ 
1 & 0 & 1 & -t^{-2} & 0 
}
\]
and dual-basis matrix
\[
\A^1_{P_W}(\Gamma) = 
\kbordermatrix{ 
& v_1^* \otimes v_0 & v_0^* \otimes v_0 & v_1^* \otimes v_1 & v_0^* \otimes v_1 \\ 
1 & 0 & 1 & 1 & 0 
}.
\]

If $\Gamma$ is oriented \MaxLR, then $\A^1_{P_W}(\Gamma)$ has matrix
\[
\A^1_{P_W}(\Gamma) =
\kbordermatrix{ 
& v_0 \otimes v_1 & v_1 \otimes v_1 & v_0 \otimes v_0 & v_1 \otimes v_0 \\ 
1 & 0 & 1 & t^2 & 0 
}
\]
and dual-basis matrix
\[
\A^1_{P_W}(\Gamma) =
\kbordermatrix{ 
& v_0 \otimes v_1^* & v_1 \otimes v_1^* & v_0 \otimes v_0^* & v_1 \otimes v_0^* \\ 
1 & 0 & -t^2 & t^2 & 0 
}.
\]
\end{proposition}

\begin{proof}
The matrix entries are again given in \cite[Table 3]{Viro}.
\end{proof}

In general, if $\Gamma$ is an $(n+2,n)$--tangle consisting of a single maximum point between strands $i+1$ and $i+2$, then
\[
\A^1_{P_W}(\Gamma): \Lambda_{B_W}^{\otimes (n+2)} \to \Lambda_{B_W}^{\otimes n}
\]
is the tensor product of identity transformations on $\Lambda_{B_W}$ with the appropriate map from Proposition \ref{ViroMinimaMapsProp} applied to tensor factors $i+1$ and $i+2$ of $\Lambda_{B_W}^{\otimes (n+2)}$.

\subsection{Modified bases}

\begin{definition}\label{ViroWjDef}
Let $Y$ be a set of $n$ points along a line, oriented by $\Sc$. For $1 \leq i \leq n$, let $w_i$ denote the tensor-product element of $\Lambda_{B_Y}^{\otimes n}$ that, in the $j^{th}$ place for $j \neq i$, has $v_0$ if $j \in \Sc$ and $v_1^*$ if $j \notin \Sc$, and in the $i^{th}$ place, has $v_1$ if $i \in \Sc$ and $v_0^*$ if $i \notin \Sc$. 
\end{definition}

\begin{remark}\label{ViroWjRemark}
As in Section \ref{BasesSect}, the definition of $1$, $w_i$, $w_{i+1}$, and $w_{i} \wedge w_{i+1}$ as specific tensors of $v$ and $v^*$ generators depends on the orientation, but the ordering $(1, w_i, w_{i+1}, w_i \wedge w_{i+1})$ of the columns and rows of each of the matrices in Corollary \ref{ViroCrossingsDualBases} remains the same for all orientations.
\end{remark}

Recall that $M_Y = H^0(Y;\frac{1}{4}\Z)$ and $B_Y = \Z[M_Y]$. Let $\W$ denote the free module over $B_Y$ formally generated by the $w_i$. As $B_Y$--modules, we have $\Lambda_{B_Y}^{\otimes n} \cong \wedge^* \W$. Since we will not keep track of the supermodule structure in the multi-graded case, we will not say more here (one could deal with supermodules as in Section \ref{BasesSect}).

\begin{definition}\label{MultivarModifiedBasisDef}
We define elements of $\Lambda_{B_Y}^{\otimes n}$ as in Definition \ref{UqSpecialEltsDef}. The variables $t_i$ below are the ones defined in (\ref{TiOnPointsEqn}). First suppose $1 \leq i \leq n-1$.
\begin{itemize}
\item If $i, i+1 \in \Sc$, define
\[
l_i := t_i^2 w_i - w_{i+1}.
\]
\item If $i \in \Sc$ but $i+1 \notin \Sc$, define
\[
l_i := t_i^2 w_i + t_{i+1}^2 w_{i+1}.
\]
\item If $i \notin \Sc$ but $i+1 \in \Sc$, define
\[
l_i := w_i - w_{i+1}.
\]
\item If $i, i+1 \notin \Sc$, define
\[
l_i := w_i + t_{i+1}^2 w_{i+1}.
\]
\end{itemize}
Define 
\begin{equation}\label{ViroL0SpecialCase}
l_0 := 
\begin{cases}
-w_1 & 1 \in \Sc \\
t_1^2 w_1 &1 \notin \Sc. 
\end{cases}
\end{equation}
Similarly, define 
\begin{equation}\label{ViroLNSpecialCase}
l_n := 
\begin{cases}
t_n^2 w_n & n \in \Sc \\
w_n & n \notin \Sc.
\end{cases}
\end{equation}
\end{definition}

\begin{definition}\label{ViroModifiedBasisDef} As in Definition \ref{ModifiedBasisDef}, the right modified basis for $\Lambda_{B_Y}^{\otimes n}$ is 
\[
\{ l_{i_1} \wedge \ldots \wedge l_{i_k} \}, 1 \leq i_i, \ldots, i_k \leq n.
\]
The left modified basis for $\Lambda_{B_Y}^{\otimes n}$ is
\[
\{ l_{i_1} \wedge \ldots \wedge l_{i_k} \}, 0 \leq i_i, \ldots, i_k \leq n-1.
\]
For $\x = \{i_1,\ldots,i_k\} \subset [1,\ldots,n]$ or $[0,\ldots,n-1]$, write
\[
l_{\x} := l_{i_1} \wedge \ldots \wedge l_{i_k}.
\]
\end{definition}

\subsection{Maps in the modified bases.}\label{ViroMapsModifiedBasisSect}
Let $\Gamma$ be a crossing, minimum, or maximum point, colored as in Definition \ref{UnivColorGraphs}. Let $W$ be the oriented $1$--manifold underlying $\Gamma$ and let $\partial W = -Y_2 \sqcup Y_1$. We want to compute the matrix for $\A^1_{P_W}(\Gamma)$ in the right and left modified bases. 

Strictly speaking, we have only defined modified bases for $\A^1_{P_{Y_i}}(Y_i)$ so far. By (\ref{ChangingPalettesIsTensorObjs}), we can take our modified basis elements in $\A^1_{P_{Y_i}}(Y_i)$ and tensor them with $1_{\Z[H^1(W,\partial W;\frac{1}{4}\Z)]}$ to obtain modified bases for $\A^1_{P_W}(Y_i)$. First, we consider the case when $\Gamma$ is a crossing. Recall that by convention, we label the crossing strands $i$ and $i+1$, where strand $i$ has positive slope and strand $i+1$ has negative slope. The generators $t_j$ of 
\begin{equation}\label{TiEquationForStrands}
\Z\bigg[H^1\bigg(W, \partial W; \frac{1}{4}\Z\bigg)\bigg] \cong \Z[t_1^{\pm 1}, \ldots, t_n^{\pm 1}]
\end{equation}
have labels corresponding to the labels of the strands. We also write
\[
\Z\bigg[H^0\bigg(Y_2, \partial Y_2; \frac{1}{4}\Z\bigg)\bigg] \cong \Z[t_1^{\pm 1}, \ldots, t_n^{\pm 1}];
\]
these $t_j$ are sent to the $t_j$ of (\ref{TiEquationForStrands}) under the above tensor product with $1_{\Z[H^1(W,\partial W;\frac{1}{4}\Z)]}$. Finally, we write
\[
\Z\bigg[H^0\bigg(Y_1, \partial Y_1; \frac{1}{4}\Z\bigg)\bigg] \cong \Z[t_1^{\pm 1}, \ldots, t_n^{\pm 1}]
\]
as well. Under the above tensor product with $1_{\Z[H^1(W,\partial W;\frac{1}{4}\Z)]}$, all of these $t_j$ except $\{t_i, t_{i+1}\}$ are sent to the $t_j$ of (\ref{TiEquationForStrands}). The variable $t_i$ is sent to $t_{i+1}$ and vice-versa because, in going left to right along the top endpoints of $W$, we encounter strand $i+1$ before strand $i$. This permutation of variables is similar to the homomorphism $\tau^{\mathbf{gr}}_i$ described in \cite[Section 4.3]{OSzNew}; see also \cite[Equation 5.1]{OSzNew}.

Now let $\Gamma$ be a minimum or maximum between strands $i+1$ and $i+2$. There is a special element $t$ of $H^1(W,\partial W; \frac{1}{4}\Z)$ corresponding to the strand containing the minimum or maximum point. The rest of the generators of $H^1(W, \partial W; \frac{1}{4}\Z)$ are labeled $t_j$, and they correspond to the strands that touch both the top and bottom boundaries of $\Gamma$.

If $t_j \in H^0(Y; \frac{1}{4}\Z)$ corresponds to a point $p_j$ in the incoming or outgoing boundary of $\Gamma$, and if $p_j$ is on the side (top or bottom) of $\Gamma$ that does not contain the critical point, then $t_j$ is sent to $t_j \in H^1(W, \partial W; \frac{1}{4}\Z)$. If $p_j$ is on the side of $\Gamma$ with the critical point, and $p_j$ lies to the left of the critical point (i.e. $j \leq i$), then $t_j \in H^0(Y; \frac{1}{4}\Z)$ is sent to $t_j \in H^1(W, \partial W; \frac{1}{4}\Z)$. If $p_j$ lies to the right of the critical point (i.e. $j \geq i+2$), then $t_j \in H^0(Y; \frac{1}{4}\Z)$ is sent to $t_{j-2} \in H^1(W, \partial W; \frac{1}{4}\Z)$. Finally, if $j = i$ or $j = i+1$, then $t_j \in H^0(Y; \frac{1}{4}\Z)$ is sent to
\[
t \in H^1\bigg(W, \partial W; \frac{1}{4}\Z\bigg).
\]

\subsubsection{Generic case}
Let $n \geq 3$ and $1 < i < n$, and suppose $\Gamma$ is a crossing between strands $i$ and $i+1$ that is positive as a braid generator.

\begin{proposition}\label{ViroCrossingMatrices} The local matrix for $\A^1_{P_W}(\Gamma)$ in the right or left modified basis is given as follows:
\begin{itemize}
\item 
If $\Gamma$ is oriented \OrUU, then $\A^1_{P_W}(\Gamma) =$
\[
\kbordermatrix{
& 1 & l_{i-1} & l_i & l_{i+1} & l_{i-1} \wedge l_i & l_{i-1} \wedge l_{i+1} & l_i \wedge l_{i+1} & l_{i-1} \wedge l_i \wedge l_{i+1} \\
1 & t_i t_{i+1} &  &  &  &  &  &  &  \\
l_{i-1} &  & t_i t_{i+1} & 0 & 0 &  &  &  &  \\
l_i &  & t_i t_{i+1}^{-1} & -t_i^{-1} t_{i+1}^{-1} & t_i^{-1} t_{i+1} &  &  &  &  \\
l_{i+1} &  & 0 & 0 & t_i t_{i+1} &  &  &  &  \\
l_{i-1} \wedge l_i &  &  &  &  & -t_i^{-1} t_{i+1}^{-1} & t_i^{-1} t_{i+1} & 0 &  \\
l_{i-1} \wedge l_{i+1} &  &  &  &  & 0 & t_i t_{i+1} & 0 &  \\
l_i \wedge l_{i+1} &  &  &  &  & 0 & t_i t_{i+1}^{-1} & -t_i^{-1} t_{i+1}^{-1} &  \\
l_{i-1} \wedge l_i \wedge l_{i+1} &  &  &  &  &  &  &  & -t_i^{-1} t_{i+1}^{-1}
}.
\]
\item 
If $\Gamma$ is oriented \OrDU, then $\A^1_{P_W}(\Gamma) =$
\[
\kbordermatrix{
& 1 & l_{i-1} & l_i & l_{i+1} & l_{i-1} \wedge l_i & l_{i-1} \wedge l_{i+1} & l_i \wedge l_{i+1} & l_{i-1} \wedge l_i \wedge l_{i+1} \\
1 & t_i^{-1} t_{i+1} &  &  &  &  &  &  &  \\
l_{i-1} &  & t_i^{-1} t_{i+1} & 0 & 0 &  &  &  &  \\
l_i &  & t_i^{-1} t_{i+1}^{-1} & t_i t_{i+1}^{-1} & -t_i t_{i+1} &  &  &  &  \\
l_{i+1} &  & 0 & 0 & t_i^{-1} t_{i+1} &  &  &  &  \\
l_{i-1} \wedge l_i &  &  &  &  & t_i t_{i+1}^{-1} & -t_i t_{i+1} & 0 &  \\
l_{i-1} \wedge l_{i+1} &  &  &  &  & 0 & t_i^{-1} t_{i+1} & 0 &  \\
l_i \wedge l_{i+1} &  &  &  &  & 0 & t_i^{-1} t_{i+1}^{-1} & t_i t_{i+1}^{-1} &  \\
l_{i-1} \wedge l_i \wedge l_{i+1} &  &  &  &  &  &  &  & t_i t_{i+1}^{-1}
}.
\]
\item
If $\Gamma$ is oriented \OrUD, then $\A^1_{P_W}(\Gamma) =$
\[
\kbordermatrix{
& 1 & l_{i-1} & l_i & l_{i+1} & l_{i-1} \wedge l_i & l_{i-1} \wedge l_{i+1} & l_i \wedge l_{i+1} & l_{i-1} \wedge l_i \wedge l_{i+1} \\
1 & t_i t_{i+1}^{-1} &  &  &  &  &  &  &  \\
l_{i-1} &  & t_i t_{i+1}^{-1} & 0 & 0 &  &  &  &  \\
l_i &  & -t_i t_{i+1} & t_i^{-1} t_{i+1} & t_i^{-1} t_{i+1}^{-1} &  &  &  &  \\
l_{i+1} &  & 0 & 0 & t_i t_{i+1}^{-1} &  &  &  &  \\
l_{i-1} \wedge l_i &  &  &  &  & t_i^{-1} t_{i+1} & t_i^{-1} t_{i+1}^{-1} & 0 &  \\
l_{i-1} \wedge l_{i+1} &  &  &  &  & 0 & t_i t_{i+1}^{-1} & 0 &  \\
l_i \wedge l_{i+1} &  &  &  &  & 0 & -t_i t_{i+1} & t_i^{-1} t_{i+1} &  \\
l_{i-1} \wedge l_i \wedge l_{i+1} &  &  &  &  &  &  &  & t_i^{-1} t_{i+1}
}.
\]
\item 
If $\Gamma$ is oriented \OrDD, then $\A^1_{P_W}(\Gamma) =$
\[
\kbordermatrix{
& 1 & l_{i-1} & l_i & l_{i+1} & l_{i-1} \wedge l_i & l_{i-1} \wedge l_{i+1} & l_i \wedge l_{i+1} & l_{i-1} \wedge l_i \wedge l_{i+1} \\
1 & -t_i^{-1} t_{i+1}^{-1} &  &  &  &  &  &  &  \\
l_{i-1} &  & -t_i^{-1} t_{i+1}^{-1} & 0 & 0 &  &  &  &  \\
l_i &  & t_i^{-1} t_{i+1} & t_i t_{i+1} & t_i t_{i+1}^{-1} &  &  &  &  \\
l_{i+1} &  & 0 & 0 & -t_i^{-1} t_{i+1}^{-1} &  &  &  &  \\
l_{i-1} \wedge l_i &  &  &  &  & t_i t_{i+1} & t_i t_{i+1}^{-1} & 0 &  \\
l_{i-1} \wedge l_{i+1} &  &  &  &  & 0 & -t_i^{-1} t_{i+1}^{-1} & 0 &  \\
l_i \wedge l_{i+1} &  &  &  &  & 0 & t_i^{-1} t_{i+1} & t_i t_{i+1} &  \\
l_{i-1} \wedge l_i \wedge l_{i+1} &  &  &  &  &  &  &  & t_i t_{i+1}
}.
\]
\end{itemize}
\end{proposition}

\begin{proof}
As with Proposition \ref{DualCanonicalRForm}, we will only show the computation for one matrix entry and one choice of orientations. Suppose $\Gamma$ is oriented \OrUUDD. We can write
\begin{align*}
l^{in}_{i-1} \wedge l^{in}_{i+1} &= (t_{i-1}^2 w_{i-1} - w_i) \wedge (w_{i+1} + t_{i+2}^2 w_{i+2}) \\
&= t_{i-1}^2 w_{i-1} \wedge w_{i+1} - w_i \wedge w_{i+1} + t_{i-1}^2 t_{i+2}^2 w_{i-1} \wedge w_{i+2} - t_{i+2}^2 w_i \wedge w_{i+2} ,
\end{align*}
where $l^{in}_j$ denotes the element $l_j$ of the incoming modified basis.

Now, $w_{i-1} \wedge w_{i+1}$ has local form $w_{i+1}$, so by the dual-basis matrix in Corollary \ref{ViroCrossingsDualBases}, we have
\[
\A^1_{P_W}(\Gamma)(w_{i-1} \wedge w_{i+1}) = t_i^{-1} t_{i+1}^{-1} w_{i-1} \wedge w_i + (t_i^4 - 1)t_i^{-1} t_{i+1}^{-1} w_{i-1} \wedge w_{i+1}
\]
(recall the definition of $w_j$ in Definition \ref{ViroWjDef}, as well as Remark \ref{ViroWjRemark}).

Similarly, 
\[
\A^1_{P_W}(\Gamma)(w_i \wedge w_{i+1}) = t_i^{-1} t_{i+1} w_i \wedge w_{i+1}.
\]

The element $w_{i-1} \wedge w_{i+2}$ has local form $1$, so
\[
\A^1_{P_W}(\Gamma)(w_{i-1} \wedge w_{i+2}) = t_i t_{i+1}^{-1} w_{i-1} \wedge w_{i+2}. 
\]

Finally, the element $w_i \wedge w_{i+2}$ has local form $w_i$, so
\[
\A^1_{P_W}(\Gamma)(w_i \wedge w_{i+2}) = -t_i t_{i+1} w_{i+1} \wedge w_{i+2}.
\]

Putting everything together, we get
\begin{align*}
\A^1_{P_W}(\Gamma)&(l^{in}_{i-1} \wedge l^{in}_{i+1}) = t_{i-1}^2 t_i^{-1} t_{i+1}^{-1} w_{i-1} \wedge w_i \\
&\qquad+ t_{i-1}^2 (t_i^4 - 1)t_i^{-1} t_{i+1}^{-1} w_{i-1} \wedge w_{i+1} \\
&\qquad- t_i^{-1} t_{i+1} w_i \wedge w_{i+1} \\
&\qquad+ t_{i-1}^2 t_i t_{i+1}^{-1} t_{i+2}^2 w_{i-1} \wedge w_{i+2} \\
&\qquad + t_i t_{i+1} t_{i+2}^2 w_{i+1} \wedge w_{i+2} \\
&= t_{i-1}^2 t_i^{-1} t_{i+1}^{-1} w_{i-1} \wedge w_i - t_{i-1}^2 t_i^{-1} t_{i+1}^{-1} w_{i-1} \wedge w_{i+1} - t_i^{-1} t_{i+1} w_i \wedge w_{i+1} \\
&\qquad+ t_{i-1}^2 t_i^3 t_{i+1}^{-1} w_{i-1} \wedge w_{i+1} + t_i^3 t_{i+1} w_i \wedge w_{i+1} \\
&\qquad \qquad + t_{i-1}^2 t_i t_{i+1}^{-1} t_{i+2}^2 w_{i-1} \wedge w_{i+2} + t_i t_{i+1} t_{i+2}^2 w_i \wedge w_{i+2} \\
&\qquad - t_i^3 t_{i+1} w_i \wedge w_{i+1} - t_i t_{i+1} t_{i+2}^2 w_i \wedge w_{i+2} + t_i t_{i+1} t_{i+2}^2 w_{i+1} \wedge w_{i+2} \\
&= t_i^{-1} t_{i+1}^{-1} (t_{i-1}^2 w_{i-1} + t_{i+1}^2 w_i) \wedge (w_i - w_{i+1}) \\
&\qquad+ t_i t_{i+1}^{-1} (t_{i-1}^2 w_{i-1} + t_{i+1}^2 w_i) \wedge (t_i^2 w_{i+1} + t_{i+2}^2 w_{i+2}) \\
&\qquad- t_i t_{i+1} (w_i - w_{i+1}) \wedge (t_i^2 w_{i+1} + t_{i+2}^2 w_{i+2}) \\
&= t_i^{-1} t_{i+1}^{-1} l_{i-1}^{out} \wedge l_i^{out} + t_i t_{i+1}^{-1} l_{i-1}^{out} \wedge l_{i+1}^{out} - t_i t_{i+1} l_i^{out} \wedge l_{i+1}^{out},
\end{align*}
where $l^{out}_j$ denotes the element $l_j$ of the outgoing modified basis. Note that $t_i$ and $t_{i+1}$ are swapped in the elements $l^{out}_j$ as a consequence of the discussion at the beginning of Section \ref{ViroMapsModifiedBasisSect}. The other computations are similar.
\end{proof}

If $\Gamma$ represents a negative braid generator, rather than a positive one, then the map associated to $\Gamma$ is the inverse of one of the above maps.

\begin{proposition}\label{ViroMinMaxProp}
If $\Gamma$ is a minimum or maximum between strands $i+1$ and $i+2$ for $0 < i < n$, then the local matrix for $\A^1_{P_W}(\Gamma)$ in the right or left modified basis is given as follows:
\begin{itemize}
\item
A minimum with matched strands $\{i+1,i+2\}$ at top, oriented \MinRL or \MinLR, is assigned the map with matrix
\[
\kbordermatrix{
& 1 & l_i \\
1 & 0 & 0 \\
l_i & 0 & 0 \\
l_{i+1} & t^{-2} & 0 \\
l_{i+2} & 0 & 0  \\
l_i \wedge l_{i+1} & 0 & t^{-2} \\
l_i \wedge l_{i+2} & 0 & 0 \\
l_{i+1} \wedge l_{i+2} & 0 & t^{-2} \\
l_i \wedge l_{i+1} \wedge l_{i+2} & 0 & 0 
}.
\]

\item
A maximum with matched strands $\{i+1,i+2\}$ on bottom, oriented \MaxRL or \MaxLR, is assigned the map with matrix
\[
\kbordermatrix{
& 1 & l_i & l_{i+1} & l_{i+2} & l_i \wedge l_{i+1} & l_i \wedge l_{i+2} & l_{i+1} \wedge l_{i+2} & l_i \wedge l_{i+1} \wedge l_{i+2} \\
1 & 0 & t^2 & 0 & t^2 & 0 & 0 & 0 & 0 \\
l_{i} & 0 & 0 & 0 & 0 & 0 & t^2 & 0 & 0 \\
}.
\]

\end{itemize}
\end{proposition}

\begin{proof} As with Proposition \ref{MinMaxModifiedBasis}, we will omit these computations to save space. They follow from Proposition \ref{ViroMinimaMapsProp}, Proposition \ref{ViroMaximaMapsProp}, and expanding out modified basis elements.
\end{proof}

\subsubsection{Special cases}

There is not much to say here; Propositions \ref{RightSpecialCaseRProp} through \ref{RSpecialCaseMinMaxProp} of Section \ref{FirstSpecialCasesSec} hold \emph{mutatis mutandis}. 

\section{Ozsv{\'a}th-Szab{\'o}'s theory}\label{OSzReviewSect}

\subsection{Algebras}\label{OSzSubsect}

For $n \geq 1$ and a subset $\Sc$ of $[1,\ldots,n]$, Ozsv{\'a}th and Szab{\'o} define a dg algebra $\B(n,\Sc)$ in Section 3 of \cite{OSzNew}. This algebra is associated to a horizontal slice across $n$ strands of an oriented knot projection. The subset $\Sc$ encodes the orientations of the strands; $i \in \Sc$ if and only if strand $i$ is oriented upwards.

The algebra $\B(n,\Sc)$ has an intrinsic grading by $(\frac{1}{4}\Z)^n$, called the multi-Alexander grading, and a homological grading by $\Z$, called the Maslov grading (in Section \ref{CompactSect} we will multiply the Maslov grading by $-1$). Abstractly, we may view the multi-Alexander grading group $(\frac{1}{4}\Z)^n$ as $H^0(Y;\frac{1}{4}\Z)$ where $Y$ is a zero-manifold consisting of $n$ points, oriented according to $\Sc$. Regardless of $\Sc$,  however, the standard elements $e_i$ of
\[
\bigg(\frac{1}{4}\Z\bigg)^n = \frac{1}{4}\Z \bigg\langle e_1, \ldots e_n\bigg\rangle
\]
are chosen to correspond to the cohomology classes of the points of $Y$ oriented negatively.

From the definitions in \cite{OSzNew}, one sees that the multi-Alexander grading takes values in $H^0(Y;\frac{1}{2}\Z) \subset H^0(Y;\frac{1}{4}\Z)$. We can obtain a single grading from the multi-Alexander gradings by pairing with the fundamental homology class $[Y] \in H_0(Y; \frac{1}{4}\Z)$, where $Y$ is oriented by $\Sc$. The single grading takes values in $\frac{1}{2}\Z \subset \frac{1}{4}\Z$.

The dg algebra $\B(n,\Sc)$ has $2^{n+1}$ elementary idempotents $\Ib_{\x}$ which correspond to subsets $\x$ of $[0,\ldots,n]$. The idempotents all have Maslov degree $0$ and multi-Alexander degree $0$. Graphically, the idempotent $\Ib_{\x}$ is depicted by drawing $n$ parallel vertical strands, labeling the regions between them as $0, \ldots, n$, and putting a dot in region $i$ if $i \in \x$.

The idempotent ring $\Ib(n)$ of $\B(n,\Sc)$ has the idempotents $\Ib_{\x}$ as generators, with multiplication defined by
\[
\Ib_{\x} \cdot \Ib_{\y} := 
\begin{cases}
\Ib_{\x}, & \x = \y \\
0, & \x \neq \y.
\end{cases}
\]
We may view $\B(n,\Sc)$ as a dg algebra over the subring $\Ib_{\x} \subset \B(n,\Sc)$.

\begin{example}
The dg algebra $\B(2,\Sc)$ for any $\Sc$ has $8$ idempotents, depicted as follows:
\begin{alignat*}{8}
&\Ib_{\varnothing} & = &\quad \hskip 0.221in | \hskip 0.221in | \\
&\Ib_{\{0\}} & = &\quad \Bigcdot | \hskip 0.221in | \\
&\Ib_{\{1\}} & = &\quad \hskip 0.221in | \Bigcdot | \\
&\Ib_{\{2\}} & = &\quad \hskip 0.221in | \hskip 0.221in | \Bigcdot \\
&\Ib_{\{0,1\}} & = &\quad \Bigcdot | \Bigcdot | \\
&\Ib_{\{0,2\}} & = &\quad \Bigcdot | \hskip 0.221in | \Bigcdot \\
&\Ib_{\{1,2\}} & = &\quad \hskip 0.221in | \Bigcdot | \Bigcdot \\
&\Ib_{\{0,1,2\}} & = &\quad \Bigcdot | \Bigcdot | \Bigcdot \\
\end{alignat*}
\end{example}

\noindent The algebras $\Cc_r(n,\Sc)$ and $\Cc_l(n,\Sc)$ are defined by truncating the idempotent ring $\Ib(n)$ of $\B(n,\Sc)$:
\begin{definition}
$\Ib_r(n)$ is defined to be the subring of $\Ib(n)$ generated by $\Ib_{\x}$ with $0 \notin \x$ (we do not require subrings to have the same multiplicative identity as the larger ring). Similarly, $\Ib_l(n)$ is the subring of $\Ib(n)$ generated by $\Ib_{\x}$ with $n \notin \x$.
\end{definition}

\begin{definition}\label{BodyCAlgDef} The dg algebra $\Cc_r(n,\Sc)$ is the algebra over $\Ib_r(n)$ obtained by restricting $\B(n,\Sc)$ to idempotents $\Ib_{\x} \in \Ib_r(n)$, i.e. 
\[
\Ccr(n,\Sc) := \bigg(\sum_{\x: 0 \notin \x} \Ib_{\x} \bigg) \cdot \B(n,\Sc) \cdot \bigg( \sum_{\x: 0 \notin \x} \Ib_{\x} \bigg).
\]
Similarly, $\Cc_l(n,\Sc)$ is obtained by restricting $\B(n,\Sc)$ to idempotents in $\Ib_l(n)$:
\[
\Ccl(n,\Sc) := \bigg(\sum_{\x: n \notin \x} \Ib_{\x} \bigg) \cdot \B(n,\Sc) \cdot \bigg( \sum_{\x: n \notin \x} \Ib_{\x} \bigg)
\]
\end{definition}

\begin{example}
The idempotents of $\Ccr(2,\Sc)$ for any $\Sc$ are
\begin{alignat*}{4}
&\Ib_{\varnothing} & = &\quad \hskip 0.221in | \hskip 0.221in | \\
&\Ib_{\{1\}} & = &\quad \hskip 0.221in | \Bigcdot | \\
&\Ib_{\{2\}} & = &\quad \hskip 0.221in | \hskip 0.221in | \Bigcdot \\
&\Ib_{\{1,2\}} & = &\quad \hskip 0.221in | \Bigcdot | \Bigcdot \\
\end{alignat*}
\end{example}

We may group the elementary idempotents according to their size as subsets of $[0,\ldots,n]$:
\begin{definition}
For $0 \leq k \leq n+1$, $\Ib(n,k)$ is the subring of $\Ib(n)$ generated by $\Ib_{\x}$ with $|\x| = k$. For $0 \leq k \leq n$, $\Ib_r(n,k)$ is the subring of $\Ib_r(n)$ generated by $\Ib_{\x}$ with $|\x| = k$, and $\Ib_l(n,k)$ is defined similarly.
\end{definition}

The following proposition follows immediately from Ozsv{\'a}th and Szab{\'o}'s construction:
\begin{proposition}\label{OSzAlgDecompProp}
We have
\begin{equation}\label{BDecompEqn}
\B(n,\Sc) = \bigoplus_{k=0}^{n+1} \B(n,k,\Sc),
\end{equation}
where $\B(n,k,\Sc)$ is the algebra over $\Ib(n,k)$ defined by restricting $\B(n,\Sc)$ to idempotents $\Ib_{\x}$ for $k$--element subsets $\x$ of $[0,\ldots,n]$, i.e.
\[
\B(n,k,\Sc) := \bigg(\sum_{\x:\,\,|\x| = k} \Ib_{\x} \bigg) \cdot \B(n,\Sc) \cdot \bigg( \sum_{\x:\,\, |\x| = k} \Ib_{\x} \bigg).
\]
Similarly,
\begin{equation}\label{CrDecompEqn}
\Ccr(n,\Sc) = \bigoplus_{k=0}^n \Ccr(n,k,\Sc)
\end{equation}
and
\begin{equation}\label{ClDecompEqn}
\Ccl(n,\Sc) = \bigoplus_{k=0}^n \Ccl(n,k,\Sc).
\end{equation}
\end{proposition}

\subsection{DA bimodules for crossings}\label{OSzCrossingSect}

Let $\Gamma$ be an oriented $(n,n)$--tangle with no maxima or minima and only one crossing between strands $i$ and $i+1$. If $\Gamma$ is positive as a braid generator, then to $\Gamma$, Ozsv{\'a}th and Szab{\'o} associate a DA bimodule over $(\B(n,\Sc'),\B(n,\Sc))$ called
\[
{^{\B(n,\Sc')}}\Pc^i_{\B(n,\Sc)},
\]
or just $\Pc^i$ when the algebras are clear from context. If $\Gamma$ is negative as a braid generator, then $\Gamma$ is assigned a DA bimodule
\[
{^{\B(n,\Sc')}}\Nc^i_{\B(n,\Sc)},
\]
or $\Nc^i$. These bimodules have a Maslov grading by $\Z$ and a multi-Alexander grading by $H^1(W, \partial W; \frac{1}{4}\Z)$, where $W$ is the abstract $1$--manifold underlying $\Gamma$. The single Alexander grading is obtained by pairing classes in $H^1(W, \partial W; \frac{1}{4}\Z)$ with the fundamental homology class of $W$ in $H_1(W, \partial W; \frac{1}{4}\Z)$.

We will not define $\Pc^i$ or $\Nc^i$ here; the definition occupies nearly all of Section 5 of \cite{OSzNew}. Fortunately, most of the intricate structure of these DA bimodules will be irrelevant for decategorification. The following is all that we need:
\begin{definition}\label{PiBasicDefB}
As $(\Ib(n),\Ib(n))$--bimodules, ${^{\B(n,\Sc')}}\Pc^i_{\B(n,\Sc)}$ and ${^{\B(n,\Sc')}}\Nc^i_{\B(n,\Sc)}$ are both defined to be the submodule of $\Ib(n) \otimes_{\F} \Ib(n)$ generated by the following elements:
\begin{itemize}
\item $\Ib_{\x} \otimes \Ib_{\x}$ for $\x \subset [0,\ldots,n]$,
\item $\Ib_{\x} \otimes \Ib_{\y}$ if $i \in \x \setminus \y$ and either $i-1 \in \y \setminus \x$ or $i+1 \in \y \setminus \x$, and $\x$ and $\y$ agree otherwise.
\end{itemize}
The bimodules $\Pc^i$ and $\Nc^i$ differ in their gradings (discussed below) and their DA operation maps $\delta^1_j$, which are defined in \cite[Section 5]{OSzNew}.
\end{definition}

\begin{figure}
\includegraphics[scale=0.625]{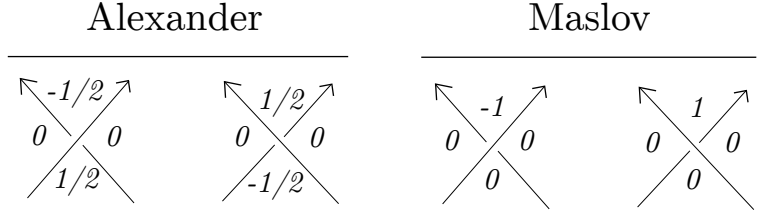}
\caption{The single Alexander and Maslov gradings of Kauffman corners.}
\label{AlexMaslovDefFig}
\end{figure}

For a general tangle diagram with no closed components, Ozsv{\'a}th and Szab{\'o} give a geometric description of the generators of the associated DA bimodule as ``partial Kauffman states.'' This description can be found in \cite[Definition 5.1]{OSzNew}. We will be concerned here with the case where the tangle diagram is a single crossing.

\begin{figure} \centering
\includegraphics[scale=0.625]{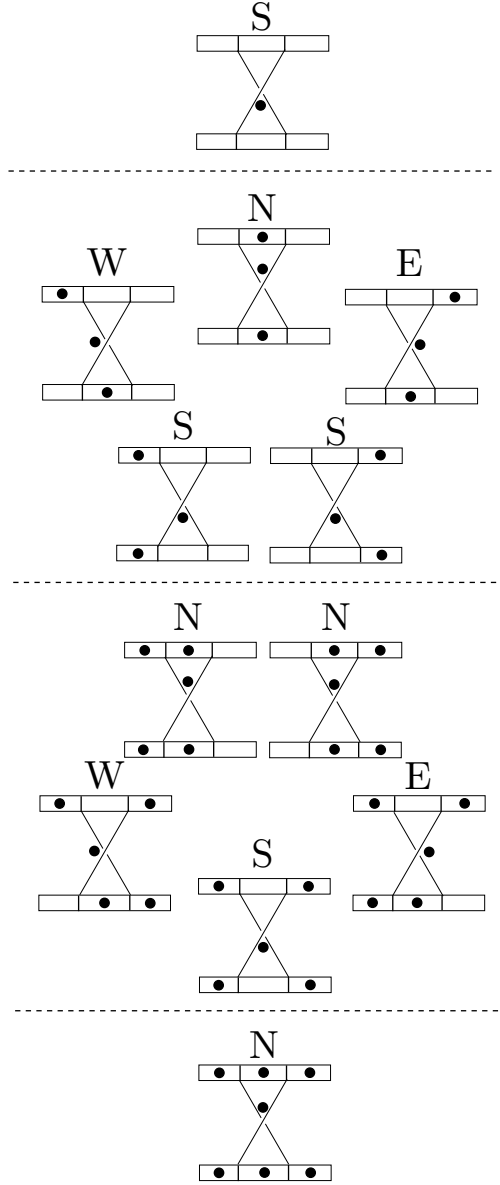}
\caption{Partial Kauffman states for a single crossing.}
\label{PartialKSFig}
\end{figure}

A partial Kauffman state in the case of a single crossing is determined by its incoming idempotent $\y$ and outgoing idempotent $\x$, as required by Definition \ref{PiBasicDefB}. If $\xi$ is such a partial Kauffman state, or equivalently if $\xi$ is a generator of $\Pc^i$ or $\Nc^i$, then we may label $\xi$ as $N$ (``north''), $S$ (``south''), $E$ (``east''), or $W$ (``west''). The label of $\xi$ is:
\[
\begin{cases}
N & \textrm{ if } \xi = \Ib_{\x} \otimes \Ib_{\x} \textrm{ with } i \in \x; \\
S & \textrm{ if } \xi = \Ib_{\x} \otimes \Ib_{\x} \textrm{ with } i \notin \x; \\
W & \textrm{ if } \xi = \Ib_{\x} \otimes \Ib_{\y} \textrm{ with } i \in \x \setminus \y, i-1 \in \y \setminus \x; \\
E & \textrm{ if } \xi = \Ib_{\x} \otimes \Ib_{\y} \textrm{ with } i \in \x \setminus \y, i+1 \in \y \setminus \x.
\end{cases}
\]
A partial Kauffman state for a single crossing is depicted by putting a dot (i.e. a Kauffman corner) in the north, south, east, or west planar region adjacent to the one crossing. See Figure \ref{PartialKSFig}; in each picture, the incoming idempotent is drawn on top and the outgoing idempotent is drawn on bottom.  

\begin{figure} \centering
\includegraphics[scale=0.625]{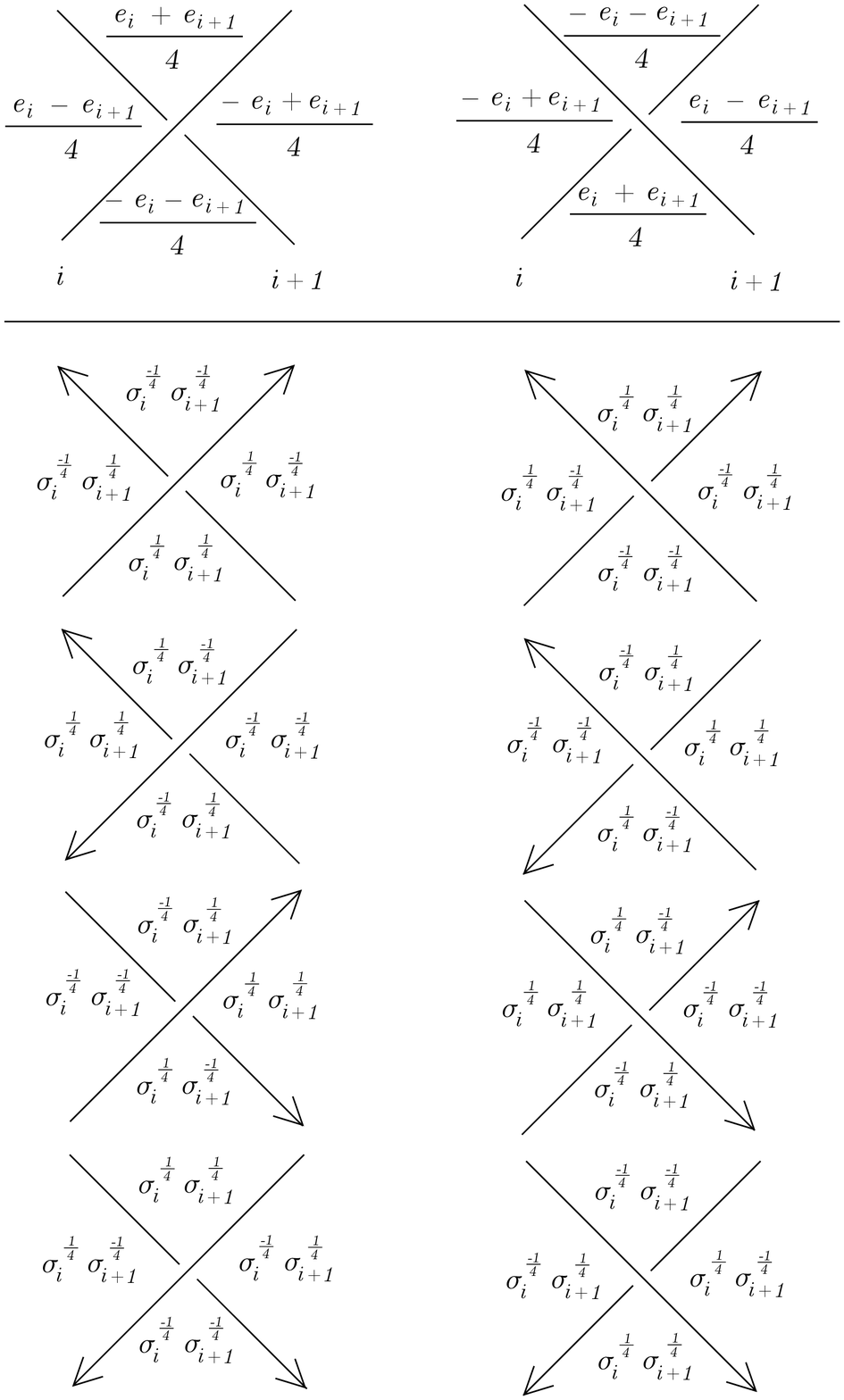}
\caption{Alexander multigradings of DA bimodule generators; see \cite[Section 4.3]{OSzNew}.}
\label{OSzMultiFig}
\end{figure}

\begin{definition}\label{LocalBigradingDef} Each generator $\xi$ of $\Pc^i$ or $\Nc^i$ has a Maslov grading $\Mas(\xi)$ and single Alexander grading $\Alex^{\sing}(\xi)$ defined by the local formulas for Kauffman corners shown in Figure~\ref{AlexMaslovDefFig}. The multi-Alexander grading $\Alex^{\multi}(\xi) \in H^1(W,\partial W; \frac{1}{4}\Z)$ is shown in Figure \ref{OSzMultiFig}. In this figure, $e_j$ is the cohomology class of the strand $s_j$, oriented downward, in $H^1(W, \partial W; \frac{1}{4}\Z)$ viewed as an additive group (as usual). In contrast, $\sigma_j$ is the cohomology class of the strand $s_j$, with orientation as depicted in the diagram (not necessarily downward), in $H^1(W, \partial W; \frac{1}{4}\Z)$ viewed as a multiplicative group. Note that these gradings depend only on the label $N$, $S$, $E$, or $W$ and whether the crossing is positive or negative as a braid generator; the formulas in terms of $\sigma_{i}^{\pm 1/4}$ and $\sigma_{i+1}^{\pm 1/4}$ also depend on the orientations of the strands (i.e. $\Sc$), because $\sigma_i$ and $\sigma_{i+1}$ depend on these orientations.
\end{definition}

\begin{warning}\label{DifferentNamesForEjWarning}
As in Warning \ref{FirstTiWarning}, $e_j$ is now a name for two slightly different objects: the cohomology class of the $j^{th}$ point in $H^0(Y, \frac{1}{4}\Z)$, oriented negatively, and the cohomology class of the $j^{th}$ strand in $H^1(W, \partial W; \frac{1}{4}\Z)$, oriented downward. The downward orientation makes sense since we are considering crossings here, not minima or maxima. As described in \cite[Equation 5.1]{OSzNew}, the generators $e_j$ of the algebra gradings on the bottom of $W$ correspond to the corresponding generators $e_j$ for $H^1(W, \partial W; \frac{1}{4}\Z)$, and the same is true for the generators $e_j$ of the algebra gradings on the top of $W$ after transforming them by the permutation $\tau_i^{\mathbf{gr}}$ that switches $e_i$ and $e_{i+1}$. 
\end{warning}

\begin{remark}
The labeling of strands as $i$ and $i+1$ in Figure \ref{OSzMultiFig} is the opposite of the labeling in \cite[Figure 20]{OSzNew}. We chose the labeling in Figure \ref{OSzMultiFig} based on \cite[Equation 5.1]{OSzNew}, as described above in Warning \ref{DifferentNamesForEjWarning}. 
\end{remark}

\subsubsection{Truncated DA bimodules}

\begin{definition}\label{PirDef}
As an $(\Ib_r(n),\Ib_r(n))$--bimodule, ${^{\Ccr(n,\Sc')}}\Pc^i_{\Ccr(n,\Sc)}$ and ${^{\Ccr(n,\Sc')}}\Nc^i_{\Ccr(n,\Sc)}$ (or $\Pc^i_r$ and $\Nc^i_r$ for shorthand) are spanned by all generators $\Ib_{\x} \otimes \Ib_{\y}$ of ${^{\B(n,\Sc')}}\Pc^i_{\B(n,\Sc)}$ and ${^{\B(n,\Sc')}}\Nc^i_{\B(n,\Sc)}$ with $0 \notin \y$, including the case $\Ib_{\x} \otimes \Ib_{\x}$ with $0 \notin \x$.
\end{definition}

Note that, by the conditions on the generators of $\Pc^i$ and $\Nc^i$ in Definition~\ref{PiBasicDefB}, we also have $0 \notin \x$ for all generators $\Ib_{\x} \otimes \Ib_{\y}$ of $\Pc^i_r$ and $\Nc^i_r$ (there is no strand $0$, so there is no crossing between strands $0$ and $1$).

From the fact that $\Pc^i$ and $\Nc^i$ are DA bimodules over $(\B(n,\Sc'),\B(n,\Sc))$, we can give $\Pc^i_r$ and $\Nc^i_r$ the structure of DA bimodules over $(\Ccr(n,\Sc'),\Ccr(n,\Sc))$. Indeed, since the generators of $\Pc^i_r$ and $\Nc^i_r$ are subsets of the generators of $\Pc^i$ and $\Nc^i$, and there are DA actions
\begin{gather*}
\delta^1_j: \Pc^i \otimes T^{j-1}(\B(n,\Sc)) \to \B(n,\Sc') \otimes \Pc^i \\
\tag*{\text{and}} \delta^1_j: \Nc^i \otimes T^{j-1}(\B(n,\Sc)) \to \B(n,\Sc') \otimes \Nc^i,
\end{gather*}
we may restrict these actions to $\Pc^i_r \otimes T^{j-1}(\Ccr(n,\Sc))$ and $\Nc^i_r \otimes T^{j-1}(\Ccr(n,\Sc))$. Since the DA actions are linear with respect to the right action of $\Ib(n)$, we see that the values of the restricted DA actions lie in the subspace of $\B(n,\Sc') \otimes \Pc^i$ or $\B(n,\Sc') \otimes \Nc^i$ spanned by elements $b \otimes (\Ib_{\x} \otimes \Ib_{\y})$ with $0 \notin \y$. By definition, $\Ib_{\x} \otimes \Ib_{\y}$ is thus an element of $\Pc^i_r$ or $\Nc^i_r$.

It follows that $0 \notin \x$ as well, so the right idempotent $\x$ of $b$ is in $\Ib_r(n)$. Since the DA actions are linear with respect to the left action of $\Ib(n)$, the left idempotent of $b$ is the left idempotent of a generator of $\Pc^i_r$ or $\Nc^i_r$. Thus, the left idempotent of $b$ is also an element of $\Ib_r(n)$, so $b \in \Ccr(n,\Sc')$.

Hence the restricted DA actions give values in $\Ccr(n,\Sc') \otimes \Pc^i_r$ or $\Ccr(n,\Sc') \otimes \Nc^i_r$, and they satisfy the DA relations because the unrestricted actions satisfy these relations.

Similarly, we get DA bimodules $\Pc^i_l$ and $\Nc^i_l$: 
\begin{definition}\label{PilDef}
As $(\Ib_l(n),\Ib_l(n))$--bimodules, ${^{\Ccl(n,\Sc')}}\Pc^i_{\Ccl(n,\Sc)}$ and ${^{\Ccl(n,\Sc')}}\Nc^i_{\Ccl(n,\Sc)}$ (or $\Pc^i_l$ and $\Nc^i_l$ for shorthand) are spanned by all generators $\Ib_{\x} \otimes \Ib_{\y}$ of ${^{\B(n,\Sc')}}\Pc^i_{\B(n,\Sc)}$ and ${^{\B(n,\Sc')}}\Nc^i_{\B(n,\Sc)}$ with $n \notin \y$, including the case $\Ib_{\x} \otimes \Ib_{\x}$ with $n \notin \x$. We give $\Pc^i_l$ and $\Nc^i_l$ the structure of DA bimodules over $(\Ccl(n,\Sc'),\Ccl(n,\Sc))$, as we did above with $\Pc^i_r$, $\Nc^i_r$, and the algebras $(\Ccr(n,\Sc'),\Ccr(n,\Sc))$.
\end{definition}

The DA bimodules $\Pc^i$, $\Pc^i_r$, and $\Pc^i_l$ respect the decompositions (\ref{BDecompEqn}), (\ref{CrDecompEqn}), and (\ref{ClDecompEqn}): for each $k$, there are DA bimodules
\begin{gather*}
{^{\B(n,k,\Sc')}}\Pc^i_{\B(n,k,\Sc)} \\
\tag*{\text{as well as}} {^{\Ccr(n,k,\Sc')}}\Pc^i_{\Ccr(n,k,\Sc)} \\
\tag*{\text{and}} {^{\Ccl(n,k,\Sc')}}\Pc^i_{\Ccl(n,k,\Sc)},
\end{gather*}
and the full DA bimodules are direct sums of the bimodules for each $k$. The same is true for $\Nc^i$, $\Nc^i_r$, and $\Nc^i_l$.

The discussion of partial Kauffman states above applies equally well to $\Pc^i_r$, $\Pc^i_l$, $\Nc^i_r$, and $\Nc^i_l$:
\begin{proposition}\label{MSigmaConcreteProp} 
The generators of $\Pc^i_r$ and $\Nc^i_r$ are partial Kauffman states, as described above and (in more generality) in \cite[Definition 5.1]{OSzNew}, in which no Kauffman corner or ``idempotent dot'' lies in the leftmost (unbounded) region of the partial knot diagram. The generators of $\Pc^i_l$ and $\Nc^i_l$ are partial Kauffman states in which no Kauffman corner or ``idempotent dot'' lies in the rightmost (also unbounded) region of the partial knot diagram.  The Maslov grading, as well as the single and multiple Alexander gradings, are as specified in Definition \ref{LocalBigradingDef}.
\end{proposition}

\begin{proof}
This proposition follows from the definitions of $\Pc^i_{r/l}$ and $\Nc^i_{r/l}$ as restrictions of $\Pc^i$ and $\Nc^i$.
\end{proof}

\subsection{DA bimodules for maxima}\label{OSzMaximaSect}
Let $\Gamma$ be an oriented $(n+2,n)$--tangle with no crossings and a maximum between strands $i+1$ and $i+2$, where $0 \leq i \leq n$. To $\Gamma$, Ozsv{\'a}th and Szab{\'o} associate a DA bimodule over $(\B(n,\Sc'),\B(n,\Sc))$ called
\[
{^{\B(n+2,\Sc')}}\Omega^{i+1}_{\B(n,\Sc)},
\]
or just $\Omega^{i+1}$, in Section 8 of \cite{OSzNew}. We recall the relevant parts of the definition:
\begin{definition}\label{OmegaIBasicDefB}
As an $(\Ib(n+2),\Ib(n))$--bimodule, ${^{\B(n+2,\Sc')}}\Omega^{i+1}_{\B(n,\Sc)}$ is the submodule of $\Ib(n+2) \otimes_{\F} \Ib(n)$ generated by the elements
\[
\Ib_{\x} \otimes \Ib_{\psi'(\x)}
\]
where $\x$ is an ``allowed idempotent state,'' i.e. $\x \cap \{i,i+1,i+2\} = \{i+1\}$, $\{i,i+1\}$, or $\{i+1,i+2\}$, and $\psi'$ is defined by
\[
\psi'(\x) := (\x \cap \{0,\ldots,i-1\}) \cup \{j - 2: j \in \x \cap \{i+3,\ldots,n\}\}
\]
if $\x \cap \{i,i+1,i+2\} = \{i+1\}$ and
\[
\psi'(\x) := (\x \cap \{0,\ldots,i-1\}) \cup \{i\} \cup \{ j - 2: j \in \x \cap \{i+3,\ldots,n\} \}
\]
if $\x \cap \{i,i+1,i+2\} = \{i,i+1\}$ or $\{i+1,i+2\}$. Note that we have changed the notation from \cite{OSzNew} slightly; our maximum point is between strands $i+1$ and $i+2$ while Ozsv{\'a}th--Szab{\'o} have it between $c$ and $c+1$. The bigrading of each generator $\Ib_{\x} \otimes \Ib_{\psi'(\x)}$ is $(0,0)$, as discussed in Section 7.1 of \cite{OSzNew} (which applies in the DA setting as well as the DD setting). The DA operations on $\Omega^{i+1}$ are defined in \cite[Section 8]{OSzNew}.
\end{definition}

As before, we define truncated versions $\Omega^{i+1}_r$ and $\Omega^{i+1}_l$:
\begin{definition}
As an $(\Ib_r(n+2),\Ib_r(n))$--bimodule, ${^{\Ccr(n+2,\Sc')}}\Omega^{i+1}_{\Ccr(n,\Sc)}$ (or $\Omega^{i+1}_r$ for shorthand) is spanned by all generators $\Ib_{\x} \otimes \Ib_{\psi'(\x)}$ of ${^{\B(n+2,\Sc')}}\Omega^{i+1}_{\B(n,\Sc)}$ with $0 \notin \psi'(\x)$. Similarly, as an $(\Ib_l(n+2),\Ib_l(n))$--bimodule, ${^{\Ccl(n+2,\Sc')}}\Omega^{i+1}_{\Ccl(n,\Sc)}$ (or $\Omega^{i+1}_l$ for shorthand) is spanned by all generators $\Ib_{\x} \otimes \Ib_{\psi'(\x)}$ of ${^{\B(n+2,\Sc')}}\Omega^{i+1}_{\B(n,\Sc)}$ with $n \notin \psi'(\x)$.
\end{definition}
If $\x \otimes \psi'(\x)$ is a generator of $\Omega^{i+1}_r$, then we also have $0 \notin \x$, and if $\x \otimes \psi'(\x)$ is a generator of $\Omega^{i+1}_l$, then we also have $n \notin \x$. Thus, as with $\Pc^i_{r/l}$ and $\Nc^i_{r/l}$, we can give $\Omega^{i+1}_{r/l}$ the structure of a DA bimodule over $(\Cc_{r/l}(n+2,\Sc'), \Cc_{r/l}(n,\Sc))$ by restricting the DA actions on $\Omega^{i+1}$.

In Section \ref{CompactSect} below, we will define the DA operations on $\Omega^1_r$ in the special case where $n = 0$ and $\Sc' = \{1\}$.

\subsection{DA bimodules for minima}\label{OSzMinimaSect}
First, let $\Gamma$ be an oriented $(n,n+2)$--tangle with no crossings and a minimum between strands $1$ and $2$. To $\Gamma$, Ozsv{\'a}th and Szab{\'o} associate a DA bimodule over $(\B(n,\Sc'),\B(n+2,\Sc))$ called
\[
{^{\B(n,\Sc')}}\mho^{1}_{\B(n+2,\Sc)},
\]
or just $\mho^1$, in Sections 9.1 and 9.2 of \cite{OSzNew}.

\begin{definition}
As an $(\Ib(n),\Ib(n+2))$--bimodule, ${^{\B(n,\Sc')}}\mho^1_{\B(n+2,\Sc)}$ is the submodule of $\Ib(n) \otimes_{\F} \Ib(n+2)$ generated by the elements
\[
\Ib_{\psi(\y)} \otimes \Ib_{\y}
\]
where $\y$ is a ``preferred idempotent state,'' i.e. 
\[
\y \cap \{0,1,2\} = \{0\}, \{2\}, \textrm{ or } \{0,2\},
\]
and $\psi$ is defined by
\[
\psi(\y) = \{j - 2: j \in \y \cap \{3,\ldots,n\}\}
\]
if $\y \cap \{0,1,2\} = \{0\}$ or $\{2\}$ and
\[
\psi(\y) = \{0\} \cup \{j - 2: j \in \y \cap \{3,\ldots,n\}\}
\]
if $\y \cap \{0,1,2\} = \{0,2\}$. As before, the bigrading of each generator $\Ib_{\psi(\y)} \otimes \Ib_{\y}$ is $(0,0)$. The DA bimodule operations on $\mho^1$ are defined in \cite[Section 9.1]{OSzNew}.
\end{definition}

For $n > 0$, the truncated bimodule $\mho^1_l$ is defined in the same way as $\mho^i$, except that it is a submodule of $\Ib_l(n) \otimes_{F} \Ib_l(n+2)$. For the other truncated bimodule $\mho^1_r$, $\y \in \Ib_r(n+2)$ is defined to be a preferred idempotent state if $\y \cap \{1,2\} = \{2\}$. The map $\psi$ and the bimodule $\mho^1_r$ are then defined like in the non-truncated case. This definition of $\mho^1_r$ works for all $n \geq 0$. Finally, for $n = 0$, $\mho^1_l$ is defined by letting $y \in \Ib_l(2)$ be a preferred idempotent state if $\y \cap \{0,1\} = \{0\}$. The map $\psi$ and the bimodule $\mho^1_l$ are again defined as in the non-truncated case.

\begin{figure} \centering
\includegraphics[scale=0.625]{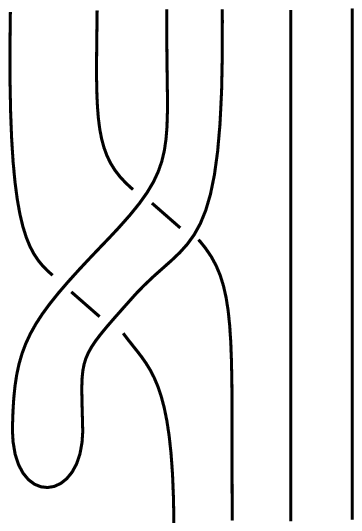}
\caption{$\mho^3 := \mho^1 \boxtimes \Pc^2 \boxtimes \Pc^1 \boxtimes \Pc^3 \boxtimes \Pc^2$.}
\label{GeneralMinFig}
\end{figure}

Now for general $0 \leq i \leq n$, let $\Gamma$ be an oriented $(n,n+2)$--tangle with no crossings and a minimum between strands $i+1$ and $i+2$. In \cite[Section 9.3]{OSzNew}, Ozsv{\'a}th and Szab{\'o} associate a DA bimodule $\mho^{i+1}$ to $\Gamma$, defined as a box tensor product of DA bimodules $\Pc^j$ for positive crossings with the bimodule $\mho^1$. This tensor product is well-defined by \cite[Proposition 3.19]{OSzNew}. The resulting DA bimodule has generators in bijection with partial Kauffman states for a tangle $\Gamma'$ isotopic to $\Gamma$, in which the minimum point has been isotoped all the way to the left above the other strands in $\Gamma$. See \cite[Figure 29]{OSzNew}, or Figure \ref{GeneralMinFig} for a minimum between strands $3$ and $4$ (the case $i = 2$). 

The truncated versions are defined using $\mho^1_{r/l}$ and $\Pc^j_{r/l}$ instead: 
\begin{definition}
The $(\Ccrl(n,\Sc'), \Ccrl(n+2,\Sc))$--bimodule $\mho^{i+1}_{r/l}$ is defined as a box tensor product of DA bimodules $\Pc^j_{r/l}$ and $\mho^1_{r/l}$, determined by isotoping the minimum point all the way to the left above the other strands, as in the non-truncated case. This tensor product is again well-defined by \cite[Proposition 3.19]{OSzNew}, which expresses the finiteness of certain sums; to show this tensor product of truncated bimodules is well-defined, we only need to check that a subset of the sums in the non-truncated case are finite, so no additional work is required.
\end{definition}

\subsection{The terminal Type A structure}\label{OSzTerminalSect}

For completeness, we briefly review the hat version $\widehat{t\mho}$ of the terminal Type A structure. This Type A structure is associated to a $(0,2)$--tangle $\Gamma$ with a marked point, viewed as the global minimum of a marked knot diagram. 

\begin{definition}[cf. Section 9.4 of \cite{OSzNew}]
For the orientation pattern \MinLR, $\widehat{t\mho}_{\B(2,\{2\})}$ has two generators $\mathbf{X}$ and $\mathbf{Y}$ as a right $\Ib(2)$--module, with 
\[
\mathbf{X} = \mathbf{X} \cdot \Ib_{\{1\}}, \,\,\,\, \mathbf{Y} = \mathbf{Y} \cdot \Ib_{\{0\}}.
\]
Similarly, for the orientation pattern \MinRL, $\widehat{t\mho}_{\B(2,\{1\})}$ has two generators $\mathbf{X}$ and $\mathbf{Y}$ with
\[
\mathbf{X} = \mathbf{X} \cdot \Ib_{\{1\}}, \,\,\,\, \mathbf{Y} = \mathbf{Y} \cdot \Ib_{\{2\}}.
\]
All generators are given Maslov grading zero and Alexander multi-grading zero.
\end{definition}
It is easy to define the Type A operations on both versions of $\widehat{t\mho}$: all non-idempotent algebra generators act as zero. The filtered version $t\mho$ and the minus version $t\mho^-$ have more structure.

We can truncate $\widehat{t\mho}$ as usual:
\begin{definition}
The truncated Type A structure $(\widehat{t\mho}_{\B(2,\{2\})})_r$ has only one generator $\mathbf{X}$, with $\mathbf{X} = \mathbf{X} \cdot \Ib_{\{1\}}$. The other truncation $(\widehat{t\mho}_{\B(2,\{2\})})_l$ has the same generators and idempotent action as $\widehat{t\mho}_{\B(2,\{2\})}$.

The truncated Type A structure $(\widehat{t\mho}_{\B(2,\{1\})})_r$ has the same generators and idempotent action as $\widehat{t\mho}_{\B(2,\{1\})}$. The other truncation $(\widehat{t\mho}_{\B(2,\{1\})})_l$ has only one generator $\mathbf{X}$, with $\mathbf{X} = \mathbf{X} \cdot \Ib_{\{1\}}$. 

The Type A operations on $\widehat{t\mho}$ are zero, so they descend to the truncated versions $\widehat{t\mho}_r$ and $\widehat{t\mho}_l$ without issue.
\end{definition}

\section{Decategorification}\label{AlgDecatSect}

\subsection{Triangulated categories and Grothendieck groups}\label{AbstractTriCatSect}

Let $\A$ be a dg algebra with a homological grading by $\Z$ and an intrinsic grading by either $\frac{1}{2}\Z$ or a free abelian group $M$, preserved by multiplication and the differential. To decategorify $\A$, one associates to $\A$ a triangulated category $\T$ and then computes the Grothendieck group $K_0(\T)$. Since $\A$ has an intrinsic grading, $K_0(\T)$ will have additional module structure.

\begin{remark} If $\A$ is a bordered algebra for a pointed matched circle representing a parametrized surface $F$, then $\A$ has an intrinsic grading by $M = H^1(F;\Z)$ (viewed multiplicatively). However, $\A$ does not have a homological grading by $\Z$ in general; instead, $\A$ is graded by a nonabelian group $\Z \ltimes H^1(F;\Z)$, so things are more difficult. See \cite{PetkovaBordered}, where homological gradings by $\Z/2\Z$ are used. 

In this paper, we will have either $M = H^0(Y;\frac{1}{4}\Z)$ for a zero-manifold $Y$ or $M = H^1(W,\partial W;\frac{1}{4}\Z)$ for a one-manifold $W$. For a zero-manifold $Y$, the isomorphism from $H^0(Y)$ to $H^1(F)$ where $F$ is the surface $D^2 \setminus \{n \textrm{ points}\}$ is suggestive in the context of the bordered surface algebras.
\end{remark}

A detailed treatment of triangulated categories may be found in \cite[Section 10.1]{KashiwaraSchapira}. Here we briefly review the Grothendieck group of a triangulated category. Recall that a category is called essentially small if it has a set of isomorphism classes of objects.

\begin{definition}\label{GGDef}
Let $\T$ be a triangulated category with shift functor $\Sigma_{\T}$, and assume $\T$ is essentially small. The Grothendieck group $K_0(\T)$ is the quotient of the free $\Z$--module generated by isomorphism classes $[X]$ of objects $X$ of $\T$ by the relations 
\[
[X]+[Z]=[Y]
\]
whenever $X \to Y \to Z \to \Sigma_{\T}(X)$ is a distinguished triangle in $\T$. It follows from these relations that $[\Sigma_{\T}(X)] = -[X]$ for all objects $X$ of $\T$.
\end{definition}

\begin{definition}\label{ZLinInducedK0Map}
Let $F: \T \to \T'$ be a triangulated functor between essentially small triangulated categories. The $\Z$--linear map
\begin{gather*}
[F]: K_0(\T) \to K_0(\T') \\
\tag*{\text{is defined by}} [F]([X]) := [F(X)].
\end{gather*}
\end{definition}

\subsection{Compact derived categories}\label{CompactSect}

Given a dg algebra $\A$ to be decategorified, there is a standard choice of triangulated category $\T$, namely the compact derived category $\D_c(\A)$ of compact objects in the unbounded derived category $\D(\A)$. It is a thick subcategory of $\D(\A)$ which is essentially small. We have
\[
D_c(\A) = \A\textrm{-perf}
\]
by e.g. \cite[Corollary 3.7]{KellerOnDG}, where $\A$-perf is the category of perfect objects of $\D(\A)$, defined as the smallest full subcategory of $\D(\A)$ containing $\A$ as a left $\A$--module and closed under grading shifts, mapping cones, and direct summands. We also have
\[
D_c(\A) \sim \D^{\pi}(\A),
\]
where $\D^{\pi}(\A)$ is the split-closed derived category defined in \cite[Section I.4c]{SeidelPL} (in a more general $\A_{\infty}$ setting). In homological mirror symmetry, the category $\D^{\pi}(\A)$ appears on the A-side when $\A$ is a Fukaya category; see \cite[Theorem 1.1]{SeidelGenus2} and \cite[Theorem 1.3]{SeidelQuartic}. By contrast, the B-side involves the bounded derived category of an abelian category rather than a dg or $\A_{\infty}$ category.

The goal of this subsection is to show why $\D_c(\A)$ is insufficient for decategorifying Ozsv{\'a}th-Szab{\'o}'s theory. This subsection will not be needed elsewhere in this paper. It is also the only place where we require knowledge of Ozsv{\'a}th-Szab{\'o}'s actual algebras and Type DA structure operations, rather than just the algebra idempotents and DA bimodule generators. See \cite{OSzNew} for these definitions, or see below for the ones we will need.

Recall that the unbounded derived category of $\A$, denoted $\D(\A)$, is the homotopy category of dg modules over $\A$, localized at quasi-isomorphisms (see \cite{Keller} for a precise definition). If $\A$ has a single intrinsic grading by either $\frac{1}{2}\Z$ or a free abelian group $M$, then the objects of $\D(\A)$ are required to have intrinsic gradings by the same group, as well as homological gradings by $\Z$. The homological shift functor $\Sigma_{\D(\A)}$ is the downward shift $[1]$ in the homological grading; see Warning \ref{PlusOrMinusOneWarn} below. Morphisms in $\D(\A)$ are required to preserve all gradings.

\begin{definition} Let $\T$ be a triangulated category with arbitrary set-indexed direct sums. An object $C$ of $\T$ is compact if the natural map
\[
\bigoplus_{\alpha \in S} \Hom_{\T}(C,X_{\alpha}) \to \Hom_{\T}\bigg(C, \bigoplus_{\alpha \in S} X_{\alpha}\bigg)
\]
is an isomorphism for all families $X_{\alpha}$ of objects of $\T$.
\end{definition}

Consider the Type DA bimodule ${^{\Ccr(2,\Sc)}}\Omega^1_{\Ccr(0,\varnothing)} =: \Omega^1_r$, where $\Sc = \{1\} \subset \{1,2\}$ (in other words, the orientation pattern is \MaxLR). While $\Omega^1_r$ is technically a DA bimodule, it is effectively a Type D structure over $\Ccr(2,\{1\})$ because $\Ccr(0,\varnothing)$ is the ground field $\F$. As a Type D structure, $\Omega^1_r$ represents a dg module over $\Ccr(2,\{1\})$ and thus an object of $\D(\Ccr(2,\{1\}))$. We show below that this object is not compact, implying that we cannot work only with compact objects when discussing Ozsv{\'a}th and Szab{\'o}'s theory.

For simplicity, we will focus on the direct summand $\Ccr(2,1,\{1\})$ of $\Ccr(2,\{1\})$; recall from Proposition \ref{OSzAlgDecompProp} that
\[
\displaystyle \Ccrl(n,\Sc) = \bigoplus_{k = 0}^n \Ccrl(n,k,\Sc).
\]

\begin{theorem}\label{NotCompact} The Type D structure $\Omega^1_{1,r}$ over $\Ccr(2,1,\{1\})$, defined to be the $k = 1$ summand of the Type D structure $\Omega^1_r$ over $\Ccr(2,\{1\})$, is not a compact object of the derived category $\D(\Ccr(2,1,\Sc))$, where the intrinsic gradings on the algebra and dg modules are by $H^1(\MaxLR, \partial(\MaxLR); \frac{1}{4}\Z)$.
\end{theorem}
\noindent It follows from Theorem \ref{NotCompact} that $\Omega^1_r$ is not a compact object of $\D(\Ccr(2,\Sc))$; in other words, Theorem \ref{IntroNoncompactThm} follows from Theorem \ref{NotCompact}.

\begin{figure} \centering
\includegraphics[scale=0.625]{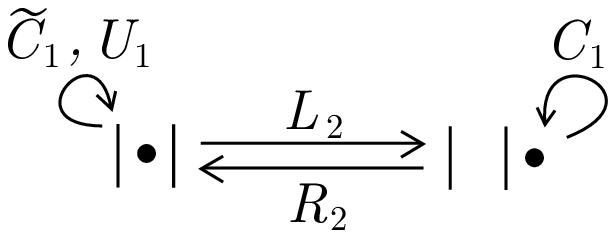}
\caption{Quiver generators of $\Ccr(2,1,\{1\})$.}
\label{DgQuiverFig}
\end{figure}

To prove Theorem \ref{NotCompact}, we must first say what $\Ccr(2,1,\{1\})$ and the Type D structure $\Omega^1_{1,r}$ are:
\begin{definition}\label{ActualAlgDef}
$\Ccr(2,1,\{1\})$ is the algebra of the quiver shown in Figure \ref{DgQuiverFig}, over $\F = \F_2$, with relations
\begin{itemize}
\item $C_1^2 = 0$
\item $\widetilde{C}_1^2 = 0$
\item $\widetilde{C}_1 R_2 = R_2 C_1$
\item $L_2 \widetilde{C}_1 = C_1 L_2$
\item $[U_1,\widetilde{C}_1] = 0$
\item $U_1 R_2 = 0$
\item $L_2 U_1 = 0$
\end{itemize}
and differential $d(\widetilde{C}_1) = U_1$; the differential is zero on the other quiver generators.

This algebra has a homological grading by $\Z$ and an intrinsic grading by 
\[
H^0\bigg(Y;\frac{1}{4}\Z\bigg) \cong \frac{1}{4}\Z \bigg\langle e_1 \bigg\rangle \oplus \frac{1}{4}\Z \bigg\langle e_2 \bigg\rangle,
\] 
where $Y$ consists of two points oriented $+-$ (recall that by convention, $e_i$ is the cohomology of the $i^{th}$ point oriented negatively, regardless of the orientations on $Y$). The generators $R_2$ and $L_2$ have homological grading zero and intrinsic grading $\frac{e_2}{2}$. The generators $C_1$ and $\widetilde{C}_1$ have homological grading $1$ and intrinsic grading $e_1$ (the homological gradings here are the negatives of those in \cite{OSzNew}, as in Warning \ref{PlusOrMinusOneWarn} below). The generator $U_1$ has homological grading $2$ and intrinsic grading $e_1$. This quiver description is not exactly how Ozsv{\'a}th and Szab{\'o} define their algebras; we leave it as an exercise to check that the quiver algebra defined here is the appropriate idempotent truncation of the algebra $\B(2,1,\{1\})$ defined in \cite[Section 3.3]{OSzNew}.
\end{definition}

\begin{warning}\label{PlusOrMinusOneWarn}
Our convention is that the differential on a dg algebra $\A$ increases the homological degree by $1$. Correspondingly, the homological shift functor $\Sigma_{\D(\A)}$ on $\D(\A)$ should be thought of as a grading shift downwards by $1$. Readers who prefer $-1$ differentials and (therefore) upwards homological shifts, as in \cite{OSzNew}, may simply reverse all the homological gradings in this section; this convention does not influence decategorification, which depends only on the parity of the homological gradings.
\end{warning}

\begin{definition}\label{OmegaOneOneDef} The Type D structure $\Omega^1_{1,r}$ has one generator $x$ with all gradings zero and (left) idempotent equal to $|\Bigcdot|$. The Type D structure operation on $\Omega^1_{1,r}$ is defined by
\[
\delta^1(x) = R_2 C_1 L_2 \otimes x.
\]
Again, this definition is equivalent to the slightly different-looking one in \cite[Section 8]{OSzNew}; note that $R_2 C_1 L_2$ is the same as $\widetilde{C}_1 R_2 L_2$, which is called $C_1 U_2$ by Ozsv{\'a}th--Szab{\'o}, but we prefer to write it in the first way.

There is an intrinsic grading by $H^1(W,\partial W; \frac{1}{4}\Z) \cong \frac{1}{4}\Z$ on $\Omega^1_{1,r}$, where $W$ is the oriented arc \MaxLR. If we let $e$ denote the cohomology class of this arc, then the $H^1(W,\partial W; \frac{1}{4}\Z)$--gradings of the algebra generators are:
\begin{itemize}
\item $\deg R_2 = \frac{1}{2}e$;
\item $\deg L_2 = \frac{1}{2}e$;
\item $\deg C_1 = -e$;
\item $\deg \widetilde{C}_1 = -e$;
\item $\deg U_1 = -e$.
\end{itemize} 
We see that $R_2 C_1 L_2$ has intrinsic degree zero and homological degree $1$, so the Type D operation on $\Omega^1_{1,r}$ is compatible with the grading structure.
\end{definition}

\begin{warning}
In this section, we implicitly identify Type D structures with their associated (left) dg modules.
\end{warning}

\begin{figure} \centering
\includegraphics[scale=0.625]{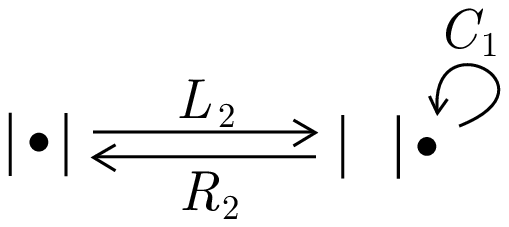}
\caption{Quiver generators of $H(\Ccr(2,1,\{1\}))$.}
\label{HomologyQuiverFig}
\end{figure}

\begin{lemma}
The homology of $\Ccr(2,1,\{1\})$ is isomorphic to the algebra of the quiver shown in Figure \ref{HomologyQuiverFig} modulo the relations:
\begin{itemize}
\item $C_1^2 = 0$
\item $C_1 L_2 R_2 = L_2 R_2 C_1$
\end{itemize}
The gradings are the same as in Definition \ref{ActualAlgDef}; in particular, $R_2$ and $L_2$ have homological degree $0$ while $C_1$ has homological degree $1$.
\end{lemma}

\begin{proof}\
A basis for the indecomposable projective left module $Q_{|\Bigcdot|}$ over $\Ccr(2,1,\{1\})$ is
\begin{align*}
\{&|\Bigcdot|, \,\, U_1^{k+1}, \,\, \widetilde{C}_1 U_1^k , \\
&(R_2 L_2)^{k+1} , \,\, L_2 (R_2 L_2)^k , \,\, C_1L_2(R_2L_2)^k , \,\, R_2C_1L_2(R_2L_2)^k  \,\,\vert \,\, k \geq 0\},
\end{align*}
and a basis for the other indecomposable projective module $Q_{| \hskip 0.07in |\Bigcdot}$ is 
\[
\{ | \hskip 0.14in |\Bigcdot, \,\, (L_2R_2)^{k+1}, \,\,  R_2(L_2R_2)^k, \,\, C_1(L_2R_2)^k, \,\, R_2C_1(L_2R_2)^k \,\, \vert \,\, k \geq 0\}.
\]
Together, these elements form a basis for $\Ccr(2,1,\{1\})$. The differential of each basis element except $\widetilde{C}_1 U_1^k$ is zero; the differential of $\widetilde{C}_1 U_1^k$ is $U_1^{k+1}$. Thus, a basis for the homology is obtained by removing $\widetilde{C}_1 U_1^k$ and $U_1^{k+1}$ from the basis elements above. The algebra of Figure \ref{HomologyQuiverFig} has the same basis as $H(\Ccr(2,1,\{1\}))$; the relations $C_1^2 = 0$ and $C_1 L_2 R_2 = L_2 R_2 C_1$ suffice to reduce any monomial in the quiver generators uniquely to one of the basis elements of $H(\Ccr(2,1,\{1\}))$. The multiplication in the quiver algebra is the same as in $H(\Ccr(2,1,\{1\}))$, so the quiver algebra is isomorphic to $H(\Ccr(2,1,\{1\}))$.
\end{proof}

\begin{lemma}\label{FormalityLemma}
The dg algebra $\Ccr(2,1,\{1\})$ is formal. 
\end{lemma}

\begin{proof}
The inclusion map $\iota$ from $H(\Ccr(2,1,\{1\}))$, viewed as the above quiver algebra, into $\Ccr(2,1,\{1\})$ is a homomorphism of dg algebras which induces an isomorphism on homology, so it is a quasi-isomorphism.
\end{proof}

From Lemma \ref{FormalityLemma}, we get an equivalence
\[
\D(\Ccr(2,1,\{1\})) \sim \D(H(\Ccr(2,1,\{1\})))
\]
of triangulated categories. In the notation of \cite[Proposition 2.4.10]{LOTBimod}, the restriction functor 
\[
\Rest_{\iota}: \D(\Ccr(2,1,\{1\})) \to \D(H(\Ccr(2,1,\{1\})))
\]
is an equivalence.

Let $\Omega$ be the Type D structure over $H(\Ccr(2,1,\{1\}))$ with one generator $x$, with grading zero and left idempotent $| \Bigcdot |$, and
\[
\delta^1(x) = R_2 C_1 L_2 \otimes x.
\]

\begin{lemma}\label{RestOmegaLemma}
The inclusion of $\Omega$ into $\Rest_{\iota}(\Omega^1_{1,r})$ is a quasi-isomorphism.
\end{lemma}

\begin{proof} 
Again, it helps to use bases. An $\F$--basis for $\Rest_{\iota}(\Omega^1_{1,r})$ (viewed as a dg module over $H(\Ccr(2,1,\{1\}))$) is
\begin{align*}
\{&x, \,\, U_1^{k+1} \cdot x, \,\, \widetilde{C}_1 U_1^k \cdot x, \\
&(R_2 L_2)^{k+1} \cdot x, \,\, L_2 (R_2 L_2)^k \cdot x, \,\, C_1L_2(R_2L_2)^k \cdot x, \,\, R_2C_1L_2(R_2L_2)^k \cdot x \,\,\vert \,\, k \geq 0\}
\end{align*}
and an $\F$--basis for $\Omega$ (viewed as a dg module over $H(\Ccr(2,1,\{1\}))$) is 
\[
\{x, \,\, (R_2 L_2)^{k+1} \cdot x, \,\,  L_2 (R_2 L_2)^k \cdot x, \,\, C_1L_2(R_2L_2)^k \cdot x, \,\, R_2C_1L_2(R_2L_2)^k \cdot x \,\,\vert\,\, k \geq 0\}.
\]
The inclusion of bases respects the algebra actions, since these actions are restricted to the smaller quiver algebra $H(\Ccr(2,1,\{1\}))$, and it yields an isomorphism on homology (the homology of each side is one-dimensional, generated by $C_1 L_2$ over $\F$).
\end{proof}

\begin{lemma}\label{ProjRes}
As a dg module over $H(\Ccr(2,1,\{1\}))$, $\Omega$ is quasi-isomorphic to the dg module $\Omega^{\proj}$ defined by
\begin{equation}\label{ProjResEqn}
\Omega^{\proj} := 
\scalebox{0.8}{\parbox{\linewidth}{
\xymatrix{
\cdots \ar[r]^-{C_1} \ar[dr]^-{1} & Q'_{| \hskip 0.07in |\Bigcdot}\{\frac{5e}{2}\} \ar[r]^-{C_1} \ar[dr]^-{1} \ar[d]^-{L_2 R_2} & Q'_{| \hskip 0.07in |\Bigcdot}\{\frac{3e}{2}\} \ar[r]^-{C_1} \ar[dr]^-{1} \ar[d]^-{L_2 R_2} & Q'_{| \hskip 0.07in |\Bigcdot}\{\frac{e}{2}\} \ar[r]^-{C_1 L_2} \ar[dr]^-{1} \ar[d]^-{L_2 R_2} & Q'_{|\Bigcdot|} \ar[d]^-{ R_2} & \\
\cdots \ar[r]_-{C_1} & Q'_{| \hskip 0.07in |\Bigcdot}\{ \frac{7e}{2} \}[-1] \ar[r]_-{C_1} & Q'_{| \hskip 0.07in |\Bigcdot}\{ \frac{5e}{2} \}[-1] \ar[r]_-{C_1} & Q'_{| \hskip 0.07in |\Bigcdot}\{ \frac{3e}{2} \}[-1] \ar[r]_-{C_1} & Q'_{| \hskip 0.07in |\Bigcdot}\{ \frac{e}{2} \}[-1]
}
}}
\end{equation}
where $Q'_{|\Bigcdot|}$ and $Q'_{| \hskip 0.07in |\Bigcdot}$ are the indecomposable projective left $H(\Ccr(2,1,\{1\}))$--modules corresponding to the two vertices of Figure \ref{HomologyQuiverFig}, and the algebra elements on the arrow labels define the differential on this dg module via right multiplication. The dg module $\Rest_{\iota}(\Omega^1_{1,r})$ is also quasi-isomorphic to $\Omega^{\proj}$.
\end{lemma}

\begin{proof}
First we compute the homology of $\Omega^{\proj}$. We note that the differential admits a filtration in which the $n^{th}$ generator in the bottom row and the $(n-1)^{st}$ generator in the top row (numbering from right to left) are in filtration level $n$. All arrows with labels other than $1$ decrease the filtration level by $1$; the arrows labeled by $1$ decrease the filtration level by $2$.

This filtration induces a spectral sequence that we can use to compute the homology. The differential on the $E_0$ page is trivial; the differential on the $E_1$ page is obtained by removing the arrows labeled $1$ from the diagram. We claim that the homology of $(E_1,d_1)$ is already one-dimensional over $\F$, generated by the generator of the bottom-right term. Given this claim, the sequence collapses at $E_2$ and the total homology is one-dimensional with the same generator (or equivalently $C_1 L_2$ times the generator of the top-right term plus $L_2 R_2$ times the generator of the bottom second-to-rightmost term, since in the total homology this sum equals the generator of the bottom-right term).

To verify the claim, note that $\F$--bases for $Q'_{|\Bigcdot|}$ and $Q'_{| \hskip 0.07in |\Bigcdot}$ respectively are
\begin{gather*}
\{| \Bigcdot |, \,\, (R_2 L_2)^{k+1}, \,\, L_2 (R_2 L_2)^k, \,\, C_1L_2(R_2L_2)^k, \,\, R_2C_1L_2(R_2L_2)^k \,\, \vert \,\, k \geq 0\} \\
\tag*{\text{and}} \{ | \hskip 0.14in |\Bigcdot, \,\, (L_2R_2)^{k+1}, \,\,  R_2(L_2R_2)^k, \,\, C_1(L_2R_2)^k, \,\, R_2C_1(L_2R_2)^k \,\, \vert \,\, k \geq 0\}.
\end{gather*}
Using these bases, the $E_1$ differential from filtration level $2$ to filtration level $1$ has matrix

\scalebox{0.9}{\parbox{\linewidth}{
\[
\kbordermatrix{ 
& (R_2 L_2)^k & L_2(R_2 L_2)^k & C_1 L_2 (R_2 L_2)^k & R_2 C_1 L_2 (R_2 L_2)^k & (L_2 R_2)^k & R_2(L_2 R_2)^k & C_1(L_2 R_2)^k & R_2 C_1 (L_2 R_2)^k \\ 
(L_2 R_2)^k & 0 & L_2 R_2 & 0 & 0 & 0 & 0 & 0 & 0 \\ 
 R_2(L_2 R_2)^k & 1 & 0 & 0 & 0 & 0 & 0 & 0 & 0 \\ 
C_1(L_2 R_2)^k & 0 & 0 & L_2 R_2 & 0 & 1 & 0 & 0 & 0 \\
R_2 C_1 (L_2 R_2)^k & 0 & 0 & 0 & L_2 R_2 & 0 & 1 & 0 & 0 
},
\]
}}

\noindent the $E_1$ differential from filtration level $3$ to filtration level $2$ has matrix

\scalebox{0.9}{\parbox{\linewidth}{
\[
\kbordermatrix{ 
& (L_2 R_2)^k & R_2(L_2 R_2)^k & C_1(L_2 R_2)^k & R_2 C_1 (L_2 R_2)^k & (L_2 R_2)^k & R_2(L_2 R_2)^k & C_1(L_2 R_2)^k & R_2 C_1 (L_2 R_2)^k \\ 
(R_2 L_2)^k & 0 & 0 & 0 & 0 & 0 & 0 & 0 & 0 \\
L_2(R_2 L_2)^k & 0 & 0 & 0 & 0 & 0 & 0 & 0 & 0 \\
C_1 L_2 (R_2 L_2)^k & 1 & 0 & 0 & 0 & 0 & 0 & 0 & 0 \\
R_2 C_1 L_2 (R_2 L_2)^k & 0 & 1 & 0 & 0 & 0 & 0 & 0 & 0 \\
(L_2 R_2)^k & L_2 R_2 & 0 & 0 & 0 & 0 & 0 & 0 & 0 \\ 
 R_2(L_2 R_2)^k & 0 & L_2 R_2 & 0 & 0 & 0 & 0 & 0 & 0 \\ 
C_1(L_2 R_2)^k & 0 & 0 & L_2 R_2 & 0 & 1 & 0 & 0 & 0 \\
R_2 C_1 (L_2 R_2)^k & 0 & 0 & 0 & L_2 R_2 & 0 & 1 & 0 & 0 
},
\]
}}

\noindent and, for $n > 2$, the $E_1$ differential from filtration level $n$ to filtration level $n-1$ has matrix

\scalebox{0.9}{\parbox{\linewidth}{
\[
\kbordermatrix{ 
& (L_2 R_2)^k & R_2(L_2 R_2)^k & C_1(L_2 R_2)^k & R_2 C_1 (L_2 R_2)^k & (L_2 R_2)^k & R_2(L_2 R_2)^k & C_1(L_2 R_2)^k & R_2 C_1 (L_2 R_2)^k \\ 
(L_2 R_2)^k & 0 & 0 & 0 & 0 & 0 & 0 & 0 & 0 \\
R_2(L_2 R_2)^k & 0 & 0 & 0 & 0 & 0 & 0 & 0 & 0 \\
C_1(L_2 R_2)^k & 1 & 0 & 0 & 0 & 0 & 0 & 0 & 0 \\
R_2 C_1 (L_2 R_2)^k & 0 & 1 & 0 & 0 & 0 & 0 & 0 & 0 \\
(L_2 R_2)^k & L_2 R_2 & 0 & 0 & 0 & 0 & 0 & 0 & 0 \\ 
 R_2(L_2 R_2)^k & 0 & L_2 R_2 & 0 & 0 & 0 & 0 & 0 & 0 \\ 
C_1(L_2 R_2)^k & 0 & 0 & L_2 R_2 & 0 & 1 & 0 & 0 & 0 \\
R_2 C_1 (L_2 R_2)^k & 0 & 0 & 0 & L_2 R_2 & 0 & 1 & 0 & 0 
}.
\]
}}

The kernel of the matrix from level $n$ to level $n-1$ is the image of the matrix from level $n+1$ to level $n$, except when $n = 1$ in which case the kernel is the whole bottom-right term and the image is everything except $1$ times the generator of the bottom-right term, proving the claim. 

The homology of $\Omega$ is also one-dimensional, with generator $C_1 L_2 \otimes x$. We have a homomorphism of dg modules $f: \Omega^{\proj} \to \Omega$ defined by the dotted arrows below:
\[
\scalebox{0.85}{\parbox{\linewidth}{
\xymatrix{
 & & & & & \Omega = Q'_{|\Bigcdot|} \save !R(.9) \ar@(ur,dr)^-{R_2 C_1 L_2} \restore \\
\cdots \ar[r]^-{C_1} \ar[dr]^-{1} & Q'_{| \hskip 0.07in |\Bigcdot}\{\frac{5e}{2}\} \ar[r]^-{C_1} \ar[dr]^-{1} \ar[d]^-{L_2 R_2} & Q'_{| \hskip 0.07in |\Bigcdot}\{\frac{3e}{2}\} \ar[r]^-{C_1} \ar[dr]^-{1} \ar[d]^-{L_2 R_2} & Q'_{| \hskip 0.07in |\Bigcdot}\{\frac{e}{2}\} \ar[r]^-{C_1 L_2} \ar[dr]^-{1} \ar[d]^-{L_2 R_2} & Q'_{|\Bigcdot|} \ar[d]_-{ R_2} \ar@{.>}[ur]^-{x} & \\
\cdots \ar[r]_-{C_1} & Q'_{| \hskip 0.07in |\Bigcdot}\{ \frac{7e}{2} \}[-1] \ar[r]_-{C_1} & Q'_{| \hskip 0.07in |\Bigcdot}\{ \frac{5e}{2} \}[-1] \ar[r]_-{C_1} & Q'_{| \hskip 0.07in |\Bigcdot}\{ \frac{3e}{2} \}[-1] \ar[r]_-{C_1} & Q'_{| \hskip 0.07in |\Bigcdot}\{ \frac{e}{2} \}[-1] \ar@{.>}[uur]_-{C_1 L_2 \otimes x} & 
}
}}
\]
The homomorphism $f$ induces an isomorphism on homology, so $f$ is a quasi-isomorphism. Finally, composing $f$ with the quasi-isomorphism from $\Omega$ to $\Rest_{\iota}(\Omega^1_{1,r})$ of Lemma \ref{RestOmegaLemma}, we get a quasi-isomorphism from $\Omega^{\proj}$ to $\Rest_{\iota}(\Omega^1_{1,r})$.  
\end{proof}

\begin{definition}
For $d \geq 1$:
\begin{itemize}
\item Let $X_d$ be the sub-dg module of $\Omega^{\proj}$ generated by the $d$ rightmost generators in the top and bottom rows of (\ref{ProjResEqn}).
\item Let $\Omega^{\tail(d)}$ be the quotient dg module $\Omega^{\proj} / X_d$.
\item Let $f_d: X_d \to \Omega$ be the composite homomorphism
\[
f_d := X_d \xrightarrow{\incl_d} \Omega^{\proj} \xrightarrow{f} \Omega,
\] 
where $f$ is the quasi-isomorphism from Lemma \ref{ProjRes}.
\item Let $\Omega^d$ be the dg module 
\[
\xymatrix @C-2.4pc {X_d[1] \ar@/^1.1pc/[rrrrrr]^{f_d} & & & \oplus & & &\Omega}.
\]
\end{itemize}
\end{definition}

\begin{warning} When $d = 1$, the symbol $\Omega^1$ here should not be confused with the $\Omega^1$ above, which was a Type D structure over $\B(2,\{1\})$. 
\end{warning}

Our strategy will be to show that 
\[
\Hom_{\D(H(\Ccr(2,1,\{1\})))}(\Omega, \Omega^1\{d-1\}) \neq 0
\]
for all $d \geq 1$. First we show that $\Omega^1\{d-1\} \cong \Omega^d$ in $\D(H(\Ccr(2,1,\{1\})))$:

\begin{lemma}\label{ShortLemma}
We have isomorphisms
\[
\Omega^d \cong \Omega^{\tail(d)} \cong \Omega^{\tail(1)}\{d-1\} \cong \Omega^1\{d-1\}
\]
in $\D(H(\Ccr(2,1,\{1\})))$.
\end{lemma}

\begin{proof}
It suffices to prove the first isomorphism; the second follows immediately from the periodicity of (\ref{ProjResEqn}) and the definition of $\Omega^{\tail(d)}$, and the third is a special case of the first.

We can split the generators in (\ref{ProjResEqn}) into the $d$ rightmost generators in the top and bottom rows (generating $X_d$), and the rest of the generators (generating $\Omega^{\tail(d)}$). If we let $g_d$ denote the differential map in (\ref{ProjResEqn}) from the $(d+1)^{st}$ rightmost generators in the top and bottom rows to the $d^{th}$ rightmost generators in the top and bottom rows, then $g_d$ defines a dg homomorphism of degree $1$ from $\Omega^{\tail(d)}$ to $X_d$.

We see from (\ref{ProjResEqn}) that
\[
\Omega^{\proj} \cong \xymatrix @C-2.4pc {\Omega^{\tail(d)} \ar@/^1.1pc/[rrrrrr]^{g_d} & & & \oplus & & & X_d},
\]
so we have a distinguished triangle
\[
\Omega^{\tail(d)}[-1] \xrightarrow{g_d} X_d \xrightarrow{\incl_d} \Omega^{\proj} \rightarrow \Omega^{\tail(d)}.
\]
Rotating one step, we get
\[
X_d \xrightarrow{\incl_d} \Omega^{\proj} \xrightarrow{\quotient} \Omega^{\tail(d)} \rightarrow X_d[1].
\]
Hence we have isomorphisms
\begin{align*}
\Omega^{\tail(d)} &\cong \xymatrix @C-2.4pc {X_d[1] \ar@/^1.1pc/[rrrrrr]^{\incl_d} & & & \oplus & & & \Omega^{\proj}} \\
&\cong \xymatrix @C-2.4pc {X_d[1] \ar@/^1.1pc/[rrrrrr]^{f \circ \, \incl_d} & & & \oplus & & & \Omega} \\
&= \xymatrix @C-2.4pc {X_d[1] \ar@/^1.1pc/[rrrrrr]^{f_d} & & & \oplus & & & \Omega} \\
&= \Omega^d
\end{align*}
in $\D(H(\Ccr(2,1,\{1\})))$.
\end{proof}

\begin{lemma}\label{LongLemma}
Let $\incl_{\Omega,d}$ denote the inclusion map from $\Omega$ into 
\[
\Omega^d = \xymatrix @C-2.4pc {X_d[1] \ar@/^1.1pc/[rrrrrr]^{f_d} & & & \oplus & & & \Omega}.
\]
The composite map
\[
\Omega^{\proj} \xrightarrow{f} \Omega \xrightarrow{\incl_{\Omega,d}} \Omega^d
\]
defines a nonzero element of $\Hom_{\D(H(\Ccr(2,1,\{1\})))} (\Omega^{\proj}, \Omega^d)$ for all $d \geq 1$. 
\end{lemma}

\begin{proof}
First, $\Omega^{\proj}$ is homotopically projective. For a reference, see \cite[Section C.8]{DrinfeldDGQuotients}; since $\Omega^{\proj}$ is an operationally bounded Type D structure, its underlying dg module is semi-free and hence homotopically projective. Alternatively, this fact is proved in bordered terms in \cite[Corollary 2.3.25]{LOTBimod}. Thus, the morphism space in the derived category from $\Omega^{\proj}$ into any other object represented by a dg module $M$ is the homology of the space of dg module homomorphisms from $\Omega^{\proj}$ to $M$, or in other words the space of dg module homomorphisms modulo null-homotopic morphisms.

Let $W$ denote the vector space over $\F$ spanned by all bigrading-preserving module homomorphisms from $\Omega^{\proj}$ to $\Omega^d$, without regard for the differential; note that the map $\incl_{\Omega,d} \circ f$ is an element of $W$. Let $V$ denote the vector space over $\F$ spanned by all module homomorphisms from $\Omega^{\proj}$ to $\Omega^d$ that preserve intrinsic gradings and decrease homological gradings by one. 

There is a linear transformation $D: V \to W$ that sends a map $h \in V$ to 
\[
D(h) := \partial_{\Omega^d} \circ h + h \circ \partial_{\Omega^{\proj}} \in W,
\]
where $\partial$ denotes the differential on $\Omega^d$ and $\Omega^{\proj}$ (to avoid confusion with the index $d$). To prove the lemma, it suffices to show that $\incl_{\Omega,d} \circ f$ is not in the image of $D$.

\begin{claim*}
$V$ and $W$ are finite-dimensional over $\F$; in fact, $\dim V = \frac{3}{2}d(d+1)$ and $\dim W = \frac{3}{2}(d+1)(d+2) - 1$. 
\end{claim*}

\begin{proof}[Proof of claim]
We will find bases for $V$ and $W$, but first we must label the generators of $\Omega^{\proj}$ and $\Omega^d$. 

Let $x_1, x_2, x_3, \ldots$ be the generators of $\Omega^{\proj}$ in the top row of (\ref{ProjResEqn}), numbering from right to left, and let $y_1, y_2, y_3, \ldots$ be the generators of $\Omega^{\proj}$ in the bottom row of (\ref{ProjResEqn}). 

Analogously to (\ref{ProjResEqn}), we may write $\Omega^d$ as:

\begin{equation}\label{OmegaDDiagram}
\scalebox{0.9}{\parbox{\linewidth}{
\xymatrix{
Q'_{| \hskip 0.07in |\Bigcdot}\{(d-\frac{3}{2})e\}[1] \ar[r]^-{C_1} \ar[ddr]^-{1} \ar[dd]^-{L_2 R_2} & \cdots \ar[r]^-{C_1} \ar[ddr]^-{1} & Q'_{| \hskip 0.07in |\Bigcdot}\{\frac{e}{2}\}[1] \ar[r]^-{C_1 L_2} \ar[ddr]^-{1} \ar[dd]^-{L_2 R_2} & Q'_{|\Bigcdot|}[1] \ar[dd]^-{ R_2} \ar[dr]^-{\omega} & \\
 & & & & \Omega = Q'_{|\Bigcdot|} \save !R(.9) \ar@(ur,dr)^-{R_2 C_1 L_2} \restore \\
Q'_{| \hskip 0.07in |\Bigcdot}\{ (d-\frac{1}{2})e \} \ar[r]_-{C_1} & \cdots \ar[r]_-{C_1} & Q'_{| \hskip 0.07in |\Bigcdot}\{ \frac{3e}{2} \} \ar[r]_-{C_1} & Q'_{| \hskip 0.07in |\Bigcdot}\{ \frac{e}{2} \} \ar[ur]_-{C_1 L_2 \cdot \omega}
}
}}
\end{equation}

As opposed to $\Omega^{\proj}$, $\Omega^d$ is finitely generated (but not operationally bounded). Let $z_1, z_2, \ldots, z_d$ be the generators of $\Omega^d$ in the top row of (\ref{OmegaDDiagram}), numbering from right to left, and let $u_1, \ldots, u_d$ be the generators of $\Omega^d$ in the bottom row of (\ref{OmegaDDiagram}). Let $\omega$ be the far-right generator in (\ref{OmegaDDiagram}); $\omega$ generates $\Omega$ as a dg submodule of $\Omega^d$. 

For a top-row generator $x_i$ of $\Omega^{\proj}$, we describe the elements of $\Omega^d$ that have the same intrinsic degree as $x_i$ but whose homological degree is $1$ less than that of $x_i$. For $j \geq i$, there is a unique basis element $a$ of the quiver algebra $H(\Ccr(2,1,\{1\}))$ such that $a \cdot z_j$ has the correct bidegree, namely 
\[
a = \begin{cases}
(L_2 R_2)^{j-i} & i > 1 \\
R_2(L_2 R_2)^{j-i-1} & i = 1 \textrm{ and } j > 1, \\
1 & i = j = 1.
\end{cases}
\]  

For $j < i$, no element $a \cdot z_j$ has the correct bidegree $x_i$, and no $b \cdot u_j$ has the correct bidegree for any $j$. Indeed, for elements $a \cdot z_j$ to have the correct homological degree, $a$ must not be divisible by $C_1$, so the intrinsic degree of $a$ is a nonnegative multiple of $e$, and thus $j \geq i$ by the pattern of shifts in (\ref{ProjResEqn}) and (\ref{OmegaDDiagram}). No algebra element has negative homological degree, so no generator $b \cdot u_j$ has homological degree one less than $x_i$. Similarly, no generator $a \cdot \omega$ has the correct homological degree. In total, the subspace of $V$ spanned by morphisms that are zero except at $x_i$ has dimension $d-i+1$, for $1 \leq i \leq d$, and dimension zero otherwise.

For a bottom-row generator $y_i$ of $\Omega^{\proj}$, elements $b \cdot z_j$ of $\Omega^d$ where $j \geq i$ and 
\[
b = \begin{cases}
C_1 (L_2 R_2)^{j-i} & j > 1, \\
C_1 L_2 & i = j = 1
\end{cases}
\]
have the correct bidegree. Also, elements $c \cdot u_j$ with $j \geq i$ and
\[
c = (L_2 R_2)^{j-i}
\]
have the correct bidegree. Algebra elements $b \cdot z_j$ and $c \cdot u_j$ with $j < i$ cannot have the correct bidegree; for the generators $b \cdot z_j$, the intrinsic degree of $b$ would have to be $\leq -\frac{3}{2}$ (we set $e = 1$ for ease of notation here), which does not occur. For the generators $c \cdot u_j$, the intrinsic degree of $c$ would have to be $\leq -1$; this degree does not occur among algebra generators not divisible by $C_1$ (indivisibility by $C_1$ is required by the homological gradings). Similarly, for generators $a \cdot \omega$, the intrinsic degree of $a$ would have to be $\leq -\frac{1}{2}$, and this degree does not occur among algebra generators indivisible by $C_1$. In total, the subspace of $V$ spanned by morphisms that are zero except at $y_i$ has dimension $2(d-i+1)$, for $1 \leq i \leq d$, and dimension zero otherwise.

Any element of $V$ may be written uniquely as a linear combination of elements of $V$ that vanish except on one generator. Thus, we may add the dimensions of the subspaces computed above to get the dimension of $V$. The result is
\begin{align*}
\sum_{i=1}^d 3(d-i+1) &= 3 \sum_{i=1}^d i \\
&= \frac{3}{2}d\bigg(d+1\bigg).
\end{align*}

We compute the dimension of $W$ similarly. For a top-row generator $x_i$ of $\Omega^{\proj}$, elements $a \cdot z_j$ and $b \cdot u_j$ of $\Omega^d$ where $j \geq i-1$ and
\[
a = \begin{cases}
C_1 (L_2 R_2)^{j-i+1} & i > 1 \textrm{ and } j > 1 \\
C_1 L_2 & i = 2 \textrm{ and } j = 1 \\
R_2 C_1 (L_2 R_2)^{j-1} & i = 1 \textrm{ and } j > 1 \\
R_2 C_1 L_2 & i = 1 \textrm{ and } j = 1,
\end{cases}
\]
\[
b = \begin{cases}
(L_2 R_2)^{j-i+1} & i > 1 \\
R_2 (L_2 R_2)^{j-1} & i = 1
\end{cases}
\]
have the same bidegree as $x_i$. Also, $\omega$ has the same bidegree as $x_1$. No elements $a \cdot z_j$ with $j < i-1$ have the same bidegree as $x_i$, since no algebra element has intrinsic degree $\leq -\frac{3}{2}$. No elements $b \cdot u_j$ with $j < i-1$ have the same bidegree as $x_i$, since no algebra element that is not divisible by $C_1$ has intrinsic degree $\leq -1$. Similarly, no elements $a \cdot \omega$ have the same bidegree as $x_i$ for $i > 1$. Finally, no element $a \cdot \omega$ has the same bidegree as $x_1$ when $a$ is a nontrivial monomial in the quiver generators of the algebra, since all such monomials have bidegree different from $(0,0)$. In total, the subspace of $W$ spanned by morphisms that are zero except at $x_i$ has dimension $2(d-i+2)$, when $1 < i \leq d+1$, dimension $2d + 1$ when $i = 1$, and dimension zero otherwise.

For a bottom-row generator $y_i$ of $\Omega^{\proj}$, elements $c \cdot u_j$ of $\Omega^d$ where $j \geq i-1$ and
\[
c = C_1 (L_2 R_2)^{j-i+1}
\]
have the same bidegree as $y_i$. Also, $C_1 L_2 \cdot \omega$ has the same bidegree as $y_1$. No elements $b \cdot z_j$ have the same homological degree as $y_i$. Elements $c \cdot u_j$ with $j < i-1$ always have intrinsic degrees greater than the intrinsic degree of $y_i$. Similarly, elements $a \cdot \omega$ always have intrinsic degree greater than the intrinsic degree of $y_i$ when $i > 1$. Finally, no element $a \cdot \omega$, other than $C_1 L_2 \cdot \omega$, has the same bidegree as $y_1$, since no algebra basis element other than $C_1 L_2$ has the right idempotents, intrinsic degree $-1/2$, and homological degree $1$.  In total, the subspace of $W$ spanned by morphisms that are zero except at $y_i$ has dimension $d-i+2$ when $1 \leq i \leq d+1$ and dimension zero otherwise.

As with $V$, we may compute the dimension of $W$ by summing the dimensions we just computed. The result is
\begin{align*}
-1 + \sum_{i=1}^{d+1} 3(d-i+2) &= -1 + 3 \sum_{i=1}^{d+1} i\\
&= \frac{3}{2}\bigg(d+1\bigg)\bigg(d+2\bigg) - 1,
\end{align*}
proving the claim.
\end{proof}

We now give names to the basis elements constructed in the proof of the above claim; we start with $V$. For $j \geq i$, let $\alpha_{i,j}$ denote the basis element of $V$ that sends $x_i$ to $a \cdot z_j$, where $a$ is specified above. Let $\beta_{i,j}$ denote the basis element of $V$ that sends $y_i$ to $b \cdot z_j$, and let $\gamma_{i,j}$ denote the basis element that sends $y_i$ to $c \cdot u_j$, where $b$ and $c$ are specified above.

Similarly, if $j \geq i-1$ then let $\alpha'_{i,j}$, $\beta'_{i,j}$, and $\gamma'_{i,j}$ denote the basis elements of $W$ that send $x_i$ to $a \cdot z_j$, $x_i$ to $b \cdot u_j$, and $y_i$ to $c \cdot u_j$ respectively ($a$, $b$, and $c$ are determined uniquely as in the proof of the claim). Let $\beta'_{1,\omega}$ denote the basis element sending $x_1$ to $\omega$, and let $\gamma'_{1,\omega}$ denote the basis element sending $y_1$ to $C_1 L_2 \cdot \omega$. Note that the element $\incl_{\Omega,d} \circ f \in W$ may be expanded in the basis for $W$ as 
\[
\incl_{\Omega,d} \circ f = \beta'_{1,\omega} + \gamma'_{1,\omega}.
\] 

We evaluate $D$ on the basis vectors $\alpha_{i,j}$, $\beta_{i,j}$, and $\gamma_{i,j}$ of $V$:
\begin{claim*}\label{DValuesClaim}
For $j > 1$, we have
\begin{itemize}
\item $D(\alpha_{i,j}) = \alpha'_{i+1,j} + \beta'_{i,j} + \alpha'_{i,j-1} + \beta'_{i,j-1}$;
\item $D(\beta_{i,j}) = \alpha'_{i+1,j} + \alpha'_{i,j} + \gamma'_{i,j} + \gamma'_{i,j-1}$;
\item $D(\gamma_{i,j}) = \beta'_{i+1,j} + \gamma'_{i+1,j} + \beta'_{i,j} + \gamma'_{i,j-1}$.
\end{itemize}
For $i = j = 1$, we have
\begin{itemize}
\item $D(\alpha_{1,1}) = \alpha'_{2,1} + \beta'_{1,1} + \beta'_{1,\omega}$;
\item $D(\beta_{1,1}) = \alpha'_{2,1} + \alpha'_{1,1} + \gamma'_{1,1} + \gamma'_{1,\omega}$;
\item $D(\gamma_{1,1}) = \beta'_{2,1} + \gamma'_{2,1} + \beta'_{1,1} + \gamma'_{1,\omega}$.
\end{itemize}
\end{claim*}

\begin{proof}[Proof of claim]
For $\alpha_{i,j}$, there is one arrow incoming to $x_i$ in (\ref{ProjResEqn}), and there are three arrows outgoing from $z_j$ in (\ref{OmegaDDiagram}). The arrow from $x_{i+1}$ to $x_i$ gives the term $\alpha'_{i+1,j}$. The arrows from $z_j$ to $u_j$, $z_{j-1}$, and $u_{j-1}$ give the terms $\beta'_{i,j}$, $\alpha'_{i,j-1}$, and $\beta'_{i,j-1}$ respectively. For $\gamma_{i,j}$, the argument is similar. For $\beta_{i,j}$, there are three arrows incoming to $y_i$ and three arrows outgoing from $z_j$. However, one of the arrows incoming to $y_i$ and one of the arrows outgoing from $z_j$ have labels divisible by $C_1$. Since $\beta_{i,j}$ sends $y_i$ to $b \cdot z_j$ where $b$ is also divisible by $C_1$, and $C_1^2 = 0$, these two arrows do not contribute to $D(\beta_{i,j})$. Hence $D(\beta_{i,j})$ has four terms rather than six.
\end{proof}

Now we can finish the proof of Lemma \ref{LongLemma}. Counting instances of $\gamma'$ basis vectors in Claim \ref{DValuesClaim}, we see that $D(\alpha_{i,j})$, $D(\beta_{i,j})$, and $D(\gamma_{i,j})$ have two $\gamma'$ vectors each in their basis expansions. Thus, any $w \in W$ that lies in the image of $D$ must have an even number of $\gamma'$ vectors in its basis expansion; in a linear combination of various $D(\alpha_{i,j})$, $D(\beta_{i,j})$, and $D(\gamma_{i,j})$, some $\gamma'$ vectors may cancel in pairs, but these cancellations do not affect the parity of the number of $\gamma'$ vectors. 

Since $\incl_{\Omega,d} \circ f = \beta'_{1,\omega} + \gamma'_{1,\omega}$, which has an odd number of $\gamma'$ vectors in its basis expansion, we see that $\incl_{\Omega,d} \circ f$ is not in the image of $D$. Therefore, 
\[
\incl_{\Omega,d} \circ f \neq 0 \in \Hom_{\D(H(\Ccr(2,1,\{1\})))}(\Omega^{\proj}, \Omega^d).
\]
\end{proof}

\begin{proof}[Proof of Theorem \ref{NotCompact}]
By Lemma \ref{FormalityLemma} and Lemma \ref{RestOmegaLemma}, it suffices to show that $\Omega$ is not a compact object of $\D(H(\Ccr(2,1,\{1\})))$. By Lemma \ref{ProjRes}, Lemma \ref{ShortLemma} and Lemma \ref{LongLemma}, we see that
\[
\Hom_{\D(H(\Ccr(2,1,\{1\})))}(\Omega, \Omega^1\{d-1\}) \neq 0
\]
for all $d \geq 1$. We will show that this statement would be impossible if $\Omega$ were compact.

Let $\A$ be a dg algebra. As above, compact objects of the derived category of $\A$ are the same as perfect dg modules over $\A$. We will show that if $M$ is a perfect dg module over $\A$, and $N$ is any other dg module over $\A$ whose homology is contained in finitely many intrinsic degrees, then $\Hom_{\D(\A)}(M,N\{d\}) = 0$ for sufficiently large $|d|$. 

It would follow that $\Omega$ is not perfect, since the homology of $\Omega^1$ is contained in finitely many intrinsic degrees, namely $-\frac{3e}{2}$ and $-e$ (or $-\frac{3}{2}$ and $-1$ if we set $e = 1$). Indeed, canceling the arrow labeled $1$ in the three-term diagram (\ref{OmegaDDiagram}, $d = 1$), we see that $\Omega^1$ is homotopy equivalent to the dg bimodule $Q'_{| \hskip 0.07in |\Bigcdot}\{\frac{e}{2}\}$ with differential $\partial(1) = C_1 L_2 R_2$. By looking at bases, the homology of this dg bimodule is generated by $C_1$ and $R_2 C_1$ over $\F$. The intrinsic degrees of these elements, after shifting downwards by $\frac{e}{2}$, are $-\frac{3e}{2}$ and $-e$ respectively.

Let $\mathcal{C}$ denote the full subcategory of $\D(\A)$ consisting of objects $M$ such that for all $N$ with $H(N)$ contained in finitely many intrinsic degrees, 
\[
\Hom_{\D(\A)}(M,N\{d\}) = 0
\]
for sufficiently large $|d|$. We want to show $\A\textrm{-perf} \subset \mathcal{C}$; by the definition of $\A$-perf, it suffices to show that $\mathcal{C}$ contains $\A$ and is closed under shifts, mapping cones, and direct summands.

Below, $N$ is an arbitrary dg module over $\A$ whose homology is contained in finitely many intrinsic degrees. We have $\A \in \mathcal{C}$ because $\Hom_{\D(\A)}(\A,N\{d\}) = H_{0,0}(N\{d\})$, the summand of the homology of $N\{d\}$ in bidegree $(0,0)$, which by assumption is zero for $|d|$ sufficiently large. 

If $M \in \mathcal{C}$, then $M[1] \in \mathcal{C}$ because $\Hom_{\D(\A)}(M[1],N\{d\}) = \Hom_{\D(\A)}(M,N[1]\{d\})$ and the homology of $N[1]$ is still contained in finitely many intrinsic degrees. Similarly, $M\{1/2\} \in \mathcal{C}$ because 
\[
\Hom_{\D(\A)}(M\{1/2\},N\{d\}) = \Hom_{\D(\A)}(M,N\{d - 1/2\}) = \Hom_{\D(\A)}(M, N\{-1/2\}\{d\})
\]
and the homology of $N\{-1/2\}$ is still contained in finitely many intrinsic degrees. The inverse shift functors work in the same way.

For mapping cones, if 
\[
M_1 \to M_2 \to M_3 \rightsquigarrow
\]
is a distinguished triangle in $\D(\A)$ and $M_1$ and $M_2$ are contained in $\mathcal{C}$, then by (e.g.) \cite[Proposition 10.1.13]{KashiwaraSchapira} we have a long exact sequence:

\[
\cdots \to \Hom_{\D(\A)}(M_1[1], N\{d\}) \to \Hom_{\D(\A)}(M_3, N\{d\}) \to \Hom_{\D(\A)}(M_2, N\{d\}) \to \cdots.
\]

\noindent For $|d|$ large enough that the terms $\Hom_{\D(\A)}(M_1[1], N\{d\})$ and $\Hom_{\D(\A)}(M_2, N\{d\})$ are zero, we have
\[
\Hom_{\D(\A)}(M_3, N\{d\}) = 0
\]
as well.

Finally, if $M_1$ is a direct summand of $M$ in $\D(\A)$, write $M = M_1 \oplus M_2$ (this direct sum should be understood as a categorical direct sum in $\D(\A)$, not necessarily a literal direct sum of dg modules). If $M \in \mathcal{C}$, then
\[
\Hom_{\D(\A)}(M,N\{d\}) \cong \Hom_{\D(\A)}(M_1,N\{d\}) \oplus \Hom_{\D(\A)}(M_2, N\{d\})
\]
which is now a literal direct sum of abelian groups. We see that if 
\[
\Hom_{\D(\A)}(M,N\{d\}) = 0
\] 
for sufficiently large $|d|$, then the same is true for $\Hom_{\D(\A)}(M_1,N\{d\})$. Thus, 
\[
\A\textrm{-perf} \subset \mathcal{C},
\]
proving that $\Omega$ is not a perfect (i.e. compact) object of $\D(H(\Ccr(2,1,\{1\})))$. 
\end{proof}

\subsection{Bounded derived categories}\label{BddTypeDSect}

The noncompactness issues discussed above only arise for minima and maxima. For crossings, we may work with finitely generated bounded Type D structures, which automatically represent compact objects of the derived category. Similarly, we may work with finitely generated left- (and right-) bounded DA bimodules. Taking box tensor products with these bimodules sends bounded Type D structures to bounded Type D structures.

\begin{definition}\label{BddTypeDCatsDef}
If $\A$ is a dg algebra with homological grading by $\Z$ and intrinsic grading by $\frac{1}{2}\Z$, let $\D^b(\A)$ denote the homotopy category of Type D structures homotopy equivalent to finitely generated bounded Type D structures over $\A$ with intrinsic gradings by $\frac{1}{2}\Z$ (as in Warning \ref{PlusOrMinusOneWarn}, we take the differentials to increase the homological gradings by one, in contrast with \cite{OSzNew}, but this convention will not matter too much). The category $\D^b(\A)$ is triangulated; the functor $\Sigma_{\T}$ is the downward shift $[1]$ in the homological grading on the Type D structures.

Now suppose $\A$ has an intrinsic grading by a free abelian group $M_{\A}$. Let $M$ be another free abelian group and let $f: M_{\A} \to M$ be a group homomorphism. Define $\D^b_M(\A)$ to be the homotopy category of Type D structures homotopy equivalent to finitely generated bounded Type D structures with intrinsic gradings by the $M_{\A}$--module $M$. Again, $\D^b_M(\A)$ is a triangulated category.
\end{definition}

\begin{remark}
Up to equivalence, we could have taken the objects of $\D^b(\A)$ and $\D^b_M(\A)$ to be finitely generated bounded Type D structures. Definition \ref{BddTypeDCatsDef} chooses to work with the smallest strictly full subcategory of the homotopy category of Type D structures over $\A$ containing the finitely generated bounded Type D structures. A third option, equivalent to the other two, is the smallest strictly full subcategory of $\D(\A)$ containing the finitely generated bounded Type D structures. Yet another equivalent option is the smallest strictly full subcategory of the homotopy category of dg modules over $\A$ containing the finitely generated bounded Type D structures.
\end{remark}

\begin{remark}\label{TypeDvsDbRem}
We use the notation $\D^b(\A)$ in this context following \cite{KontsevichICM} and \cite{SeidelQuartic}.
\end{remark}

\begin{remark}
We expect that $\D^b(\Ccrl(n,\Sc))$ and $\D^b_{M_Y}(\Ccrl(n,\Sc))$ are equivalent to $\D_c(\Ccrl(n,\Sc))$ when the latter is defined with single gradings or multi-gradings by $M_Y = H^0(Y; \frac{1}{4}\Z)$ respectively. However, we have not been able to prove this equivalence. As in Section \ref{CompactSect}, $\D_c(\Ccrl(n,\Sc))$ is the same as $\Ccrl(n,\Sc)$-perf, which is equivalent to the homotopy category of Type D structures that admit a split-injection (in the homotopy category) into a finitely generated bounded Type D structure. 

Let $H(\Ccrl(n,\Sc)-\dgmod)$ denote the homotopy category of dg modules over $\Ccrl(n,\Sc)$. We see that $\D_c(\Ccrl(n,\Sc))$ is equivalent to the smallest thick subcategory of $H(\Ccrl(n,\Sc)-\dgmod)$ containing the finitely generated bounded Type D structures. The category $\D^b(\Ccrl(n,\Sc))$ is equivalent to the smallest strictly full triangulated subcategory of $H(\Ccrl(n,\Sc)-\dgmod)$ containing the finitely generated bounded Type D structures, and this triangulated subcategory is dense in $\D_c(\Ccrl(n,\Sc))$. By Thomason \cite[Theorem 2.1]{Thomason}, $\D^b(\Ccrl(n,\Sc))$ and $\D_c(\Ccrl(n,\Sc))$ disagree if and only if $K_0(\D^b(\Ccrl(n,\Sc)))$ is a proper subgroup of $K_0(\D_c(\Ccrl(n,\Sc)))$.

Thus, either $\D^b(\Ccrl(n,\Sc)) \sim \D_c(\Ccrl(n,\Sc))$, or $K_0(\D_c(\Ccrl(n,\Sc)))$ contains some elements that are not visible in the subgroup $K_0(\D^b(\Ccrl(n,\Sc)))$. In the latter case, it would be interesting to know what these elements are. At any rate, we will work with $K_0(\D^b(\Ccrl(n,\Sc)))$ rather than $K_0(\D_c(\Ccrl(n,\Sc)))$ in this paper. 
\end{remark}

If we give $\Ccrl(n,\Sc)$ single or multiple intrinsic gradings, we get corresponding module structures on $K_0(\D^b(\Ccrl(n,\Sc)))$:
\begin{definition}\label{ModuleStrsOnGG}
Let $\{-\}$ denote the downward shift in the intrinsic gradings. We define a $\Z[t^{\pm 1/2}]$--module structure on $K_0(\D^b(\Ccrl(n,\Sc)))$ by
\[
t^{\pm 1/2}[X] := [X\{\mp 1/2\}],
\]
and we define a $\Z[M]$--module structure on $K_0(\D^b_M(\Ccrl(n,\Sc)))$ by
\[
m[X] := [X\{-m\}]
\]
for $m \in M$. 
\end{definition}

\subsubsection{Triangulated functors from DA bimodules}\label{EFunFromDASect}

Let $X$ be a finitely generated left bounded DA bimodule over dg algebras $(\A',\A)$ with a single intrinsic grading by $\frac{1}{2}\Z$. We have a triangulated functor
\[
X \boxtimes -: \D^b(\A) \to \D^b(\A'). 
\]
This functor induces a map on $K_0$ by Definition \ref{ZLinInducedK0Map}. In fact, this map is linear over $\Z[t^{\pm 1/2}]$:
\begin{proposition}\label{SinglyGrK0MapLinear}
If $X$ has an intrinsic grading by $\frac{1}{2}\Z$, then the map
\[
[X \boxtimes -]: K_0(\D^b(\A)) \to K_0(\D^b(\A'))
\]
is linear over $\Z[t^{\pm 1/2}]$.
\end{proposition}

\begin{proof}
Let $[Y] \in K_0(\D^b(\A))$. We have
\begin{align*}
[X \boxtimes -](t^{1/2} \cdot [Y]) &= [X \boxtimes -]([Y\{-1/2\}]) \\
&= [X \boxtimes (Y\{-1/2\})] \\
&= [(X \boxtimes Y)\{-1/2\}] \\
&= t^{1/2} [X \boxtimes Y] \\
&= t^{1/2} [X \boxtimes -]([Y]).
\end{align*}
\end{proof}

Now suppose $\A'$ has intrinsic grading group $M_{\A'}$, $\A$ has intrinsic grading group $M_{\A}$, and $X$ has intrinsic grading group $M_X$, where we are given homomorphisms from $M_{\A'}$ and $M_{\A}$ into $M_X$. Let $M$ be another free abelian group equipped with a homomorphism from $M_{\A}$, and define $f$ to be the inclusion map from $M$ to 
\[
M_X \times_{M_{\A}} M := \frac{M_X \oplus M}{M_{\A}}.
\] 
By slight abuse of notation, we may identify $f$ with the induced homomorphism $f: \Z[M_{\A}] \to \Z[M]$ of group rings.
\begin{proposition}\label{MultiplyGrK0MapLinear}
The map
\[
[X \boxtimes -]: K_0(\D^b_{M_{\A}}(\A)) \to K_0(\D^b_M(\A'))
\]
intertwines the $\Z[M_{\A}]$--action on its domain with the $\Z[M]$--action on its codomain via the homomorphism $f \Z[M_{\A}] \to \Z[M]$.
\end{proposition}
\begin{proof}
Let $[Y] \in K_0(\D^b_{M_{\A}}(\A))$, and let $m \in M_{\A}$. We have
\begin{align*}
[X \boxtimes -](m \cdot [Y]) &= [X \boxtimes -]([Y\{m^{-1}\}]) \\
&= [X \boxtimes (Y\{m^{-1}\})] \\
&= [(X \boxtimes Y)\{f(m)^{-1}\}] \\
&= f(m) \cdot [X \boxtimes Y] \\
&= f(m) \cdot [X \boxtimes -]([Y]).
\end{align*}
\end{proof}

\subsubsection{Computing the Grothendieck group: single gradings}

The algebra $\Ccrl(n,\Sc)$ has a natural augmentation map
\[
\epsilon: \Ccrl(n,\Sc) \to \Ibrl(n),
\]
where $\Ibrl(n)$ is the idempotent ring of $\Ccrl(n,\Sc)$. Without going into too much detail, $\Ccrl(n,\Sc)$ is generated by certain elements $R_i$, $L_i$, $U_i$, and $C_i$ over $\Ibrl(n)$, and $\epsilon$ sends each of these generators to zero. The generators $C_1$ and $\widetilde{C}_1$ of Section \ref{CompactSect} are both considered $C_i$ generators.

\begin{definition}\label{XEpsDef}
The augmentation map $\epsilon$ defines a one-dimensional left (and right) bounded DA bimodule $X_{\epsilon}$ over $(\Ibrl(n), \Ccrl(n,\Sc))$. This bimodule has one generator, and the only nonzero DA actions have a single element $a$ of $\Ccrl(n,\Sc)$ acting on the right on the one generator of $X_{\epsilon}$ and outputting $\epsilon(a)$ on the left. 
\end{definition}

View both $\Ccrl(n,\Sc)$ and $\Ibrl(n)$ as dg algebras with homological gradings by $\Z$ and intrinsic gradings by $\frac{1}{2}\Z$. Then we may view $X_{\epsilon}$ as having the same grading structure; its one generator has bidegree $(0,0)$. Taking box tensor product with $X_{\epsilon}$, we get a triangulated functor
\[
X_{\epsilon} \boxtimes -: \D^b(\Ccrl(n,\Sc)) \to \D^b(\Ibrl(n))
\]
and thus a $\Z[t^{\pm 1/2}]$--linear map
\[
[X_{\epsilon} \boxtimes -]: K_0(\D^b(\Ccrl(n,\Sc))) \to K_0(\D^b(\Ibrl(n)))
\]
by Proposition \ref{SinglyGrK0MapLinear}. 

\begin{lemma}\label{SingleGrIdemDecatLemma}
The Grothendieck group $K_0(\D^b(\Ibrl(n)))$ is a free module of rank $2^n$ over $\Z[t^{\pm 1/2}]$. The classes $[S_{\Ib_{\x}}]$ of the one-dimensional modules $S_{\Ib_{\x}}$ over $\Ibrl(n)$, for $\x \subset [0,\ldots,n-1]$ or $\x \subset [1,\ldots,n]$, form a basis.
\end{lemma}

\begin{proof}
Any object $X$ of $\D^b(\Ibrl(n))$ is isomorphic to its homology, which is a Type D structure with zero Type D operation (i.e. $\delta^1 = 0$). Such a Type D structure is a direct sum of grading-shifted copies of $S_{\Ib_{\x}}$ for various $\x$, so $[X]$ can be written as a linear combination of the corresponding classes $[S_{\Ib_{\x}}]$.

Similarly, any distinguished triangle in $\D^b(\Ibrl(n))$ is isomorphic to the mapping cone on the zero morphism between two Type D structures, so the only relations in $K_0(\D^b(\Ibrl(n)))$ come from direct sums and shifts. If the classes $[S_{\Ib_{\x}}]$ were not linearly independent, then either a direct sum of shifts of some Type D structure $S_{\Ib_{\x}}$ would be isomorphic to zero in $\D^b(\Ibrl(n))$, or a direct sum of shifts of some $S_{\Ib_{\x}}$ would be isomorphic to a direct sum of shifts of various $S_{\Ib_{\x'}}$ with $\x' \neq \x$. Both of these options are impossible since homology is preserved under homotopy equivalence of Type D structures, so the classes are independent. 
\end{proof}

For every idempotent $\Ib$ of $\Ccrl(n,\Sc)$, recall that we have a (bounded) Type D structure $Q_{\Ib} := \Ccrl(n,\Sc) \cdot \Ib$ over $\Ccrl(n,\Sc)$ with one generator $\Ib$ and $\delta^1 = 0$. 

\begin{theorem}[cf. Theorem \ref{AlgDecatIntroThm}]\label{SingleAlgDecatBodyThm} Let $\x$ range over subsets of $[0,n-1]$ for $\Ccl(n,\Sc)$ or subsets of $[1,n]$ for $\Ccr(n,\Sc)$. The $2^n$ classes $[Q_{\Ib_{\x}}]$ form a basis, over $\Z[t^{\pm 1/2}]$, for $K_0(\D^b(\Ccrl(n,\Sc)))$.
\end{theorem}

\begin{proof} 
The proof follows Theorem 21 of Petkova \cite{PetkovaBordered}; our setting is easier since we have a homological grading by $\Z$. 

The classes $[Q_{\Ib_{\x}}]$ are linearly independent in $K_0(\D^b(\Ccrl(n,\Sc)))$ since their images $[S_{\Ib_{\x}}]$ under $[X_{\epsilon} \boxtimes -]$ form a basis for $K_0(\D^b(\Ibrl(n)))$ by Lemma \ref{SingleGrIdemDecatLemma}. 

We want to see that the classes $[Q_{\Ib_{\x}}]$ span $K_0(\D^b(\Ccrl(n,\Sc)))$. An arbitrary element of this Grothendieck group may be written as $[D]$, where $D$ is an operationally bounded Type D structure that has finitely many generators $\{y_1,\ldots,y_k\}$. We may assume each $y_i$ has a unique left idempotent $\x_i$.

By the boundedness of $D$, we may assume (after reordering generators if necessary) that $\delta^1(y_k) = 0$. Then 
\[
\Ccrl(n,\Sc) \cdot y_k \cong Q_{\Ib_{\x_k}} [-\Mas(y_k)] \{-\Alex^{\sing}(y_k) \}
\]
is a subcomplex of $\Ccrl(n,\Sc) \cdot D$. Here, by analogy with Type D structures coming from Heegaard Floer homology, we denote the homological grading of $y_k$ by $\Mas(y_k)$ and the half-integer intrinsic grading of $y_k$ by $\Alex^{\sing}(y_k)$.

Define 
\[
D' := \frac{D}{\Ccrl(n,\Sc) \cdot y_k}.
\] 
We have a distinguished triangle
\[
D'[-1] \to \Ccrl(n,\Sc) \cdot y_k \to D \rightarrow D'
\]
in $\D^b(\Ccrl(n,\Sc))$, so we can write 
\[
[D] = (-1)^{\Mas(y_k)} t^{\Alex^{\sing}(y_k)} [Q_{\Ib_{\x_k}}] + [D'].
\]
Thus, $D'$ can be generated by strictly fewer than $k$ elements, namely $y_1, \ldots, y_{k-1}$. By induction, the classes $[Q_{\Ib_{\x}}]$ span $K_0(\D^b(\Ccrl(n,\Sc)))$ as a module over $\Z[t^{\pm 1/2}]$.
\end{proof}

Theorem \ref{SingleAlgDecatBodyThm} allows us to make the identification
\[
K_0(\D^b(\Ccrl(n,\Sc))) \otimes_{\frac{\Z[t^{\pm 1/2}, q^{\pm 1}]}{q = t^{1/2}}} \C(q) \cong V^{\otimes \Sc}
\]
of $\C(q)$--modules, where $V^{\otimes \Sc}$ is the representation of $\Uq$ defined in the introduction. The explicit identification is given by
\[
[Q_{\Ib_{\x}}] \leftrightarrow l_{\x_{i_1}} \wedge \cdots \wedge l_{\x_{i_k}},
\]
where $l_{\x_{i_1}} \wedge \cdots \wedge l_{\x_{i_k}}$ is an element of the right (for $\Ccr(n,\Sc)$) or left (for $\Ccl(n,\Sc)$) modified basis for $V^{\otimes \Sc}$ as defined in Section \ref{BasesSect}.

In Section \ref{DADecatSect}, we will identify the maps $[X \boxtimes -]$ on these Grothendieck groups, where $X$ is the DA bimodule for a braid $\Gamma$, with the corresponding $\Uq$--linear maps associated to $\widetilde{\Gamma}$ in Section \ref{MapsForCrossingsSect} (using the geometric identification of tangles $\Gamma \leftrightarrow \widetilde{\Gamma}$ from Figure \ref{RotateReverseFig}). 

\subsubsection{Computing the Grothendieck group: multiple gradings}

Now we will consider the multiply graded setting. View both $\Ccrl(n,\Sc)$ and $\Ibrl(n)$ as dg algebras with homological gradings by $\Z$ and intrinsic gradings by $M_Y:= H^0(Y,\frac{1}{4}\Z)$, where $Y$ is (as usual) a zero-manifold consisting of $n$ points oriented according to $\Sc$. 

We use the DA bimodule $X_{\epsilon}$ from Definition \ref{XEpsDef} again. This time, we view the grading group of $X_{\epsilon}$ as $M_Y$, equipped with identity homomorphisms from the grading groups of the incoming algebra $\Ccrl(n,\Sc)$ and the outgoing algebra $\Ibrl(n)$. The single generator of $X_{\epsilon}$ has degree zero.

Let $M$ be a free abelian group equipped with a homomorphism from $M_Y$; we consider the category $\D^b_M(\Ccrl(n,\Sc))$. As in Section \ref{EFunFromDASect}, taking box tensor product with $X_{\epsilon}$ induces a triangulated functor
\[
X_{\epsilon} \boxtimes -: \D^b_M(\Ccrl(n,\Sc)) \to \D^b_M(\Ibrl(n)),
\]
which in turn induces a $\Z[M]$--linear map on Grothendieck groups by Proposition \ref{MultiplyGrK0MapLinear}. Note that we are using the identification $M = M_Y \times_{M_Y} M$.

\begin{lemma}\label{MultiIdemDecatLemma}
The Grothendieck group $K_0(\D^b_M(\Ibrl(n)))$ is a free module of rank $2^n$ over $\Z[M]$. The classes $[S_{\Ib_{\x}}]$ of the one-dimensional modules $S_{\Ib_{\x}}$ for $\x \subset [0,\ldots,n-1]$ or $\x \subset [1,\ldots,n]$ form a basis.
\end{lemma}

\begin{proof}
The proof is the same as for Lemma \ref{SingleGrIdemDecatLemma}.
\end{proof}

As in the singly graded setting, we have a bounded Type D structure $Q_{\Ib} := \Ccrl(n,\Sc) \cdot \Ib$ over $\Ccrl(n,\Sc)$ for each elementary idempotent $\Ib$ of $\Ccrl(n,\Sc)$. 

\begin{theorem}[cf. Theorem \ref{MultiAlgDecatIntroThm}]\label{MultiAlgDecatBodyThm} Let $\x$ range over subsets of $[0,n-1]$ for $\Ccl(n,\Sc)$ or subsets of $[1,n]$ for $\Ccr(n,\Sc)$. The $2^n$ classes $[Q_{\Ib_{\x}}]$ form a basis, over $\Z[M]$, for $K_0(\D^b_M(\Ccrl(n,\Sc)))$.
\end{theorem}

\begin{proof} 
The proof is analogous to the proof of Theorem \ref{SingleAlgDecatBodyThm}.
\end{proof}

\begin{corollary}\label{K0ofWIsTensorProd}
Let $M$ be a free abelian group equipped with a homomorphism from $M_Y = H^0(Y; \frac{1}{4}\Z)$. As $\Z[M]$--modules,
\[
K_0(\D^b_M(\Ccrl(n,\Sc))) \cong K_0(\D^b_{M_Y}(\Ccrl(n,\Sc))) \otimes_{\Z[M_Y]} \Z[M]. 
\]
\end{corollary}

\begin{proof}
By Theorem \ref{MultiAlgDecatBodyThm}, we see that both $K_0(\D^b_M(\Ccrl(n,\Sc)))$ and $K_0(\D^b_{M_Y}(\Ccrl(n,\Sc)))$ are free modules over $\Z[M_Y]$ and $\Z[M]$ respectively, with the same basis $\{Q_{\Ib_{\x}}\}$. 
\end{proof}

Specializing to $M = M_Y$, Theorem \ref{MultiAlgDecatBodyThm} allows us to make the identification
\[
K_0(\D^b_{M_Y}(\Ccrl(n,\Sc))) \cong \A^1_{P_Y}(Y),
\]
where $\A^1_{P_Y}(Y)$ is Viro's invariant for $Y$ colored by $P_Y$ as in Section \ref{ViroObjectInvsSect} (see also Definition \ref{UnivColorObjects}). As in the singly graded case, a basis element $[Q_{\Ib_{\x}}]$ of $K_0$ corresponds to the modified basis element $l_{\x} = l_{x_{i_1}} \wedge \cdots \wedge l_{x_{i_k}}$. 

Similarly, specializing to $M = M_W$ where $W$ is an oriented $1$--manifold and $M_W = H^1(W,\partial W; \frac{1}{4}\Z)$, we get
\[
K_0(\D^b_{M_W}(\Ccrl(n,\Sc))) \cong \A^1_{P_Y}(Y) \otimes_{\Z[M_Y]} \Z[M_W] \cong \A^1_{P_W}(Y),
\]
with the same identification of basis elements.

In Section \ref{DADecatSect}, we will identify the maps $[X \boxtimes -]$ on these Grothendieck groups, where $X$ is the DA bimodule for a braid $\Gamma$, with Viro's invariants $\A^1_{P_W}(\widetilde{\Gamma})$ discussed in Section \ref{ViroMapsForCrossingsSect} (using the geometric identification of tangles $\Gamma \leftrightarrow \widetilde{\Gamma}$ from Figure \ref{RotateReverseFig} again).

\subsubsection{Decategorifying DA bimodules for crossings: singly graded case}\label{SingleBraidDecatSect}
We compute the decategorification of DA bimodules for crossings abstractly here. Computations of actual matrices will be done in Section \ref{DADecatSect}. 

First we consider the singly graded case. Let $X$ be a finitely generated left bounded DA bimodule over $(\Ccrl(n,\Sc'),\Ccrl(n,\Sc))$ with a single intrinsic grading by $\frac{1}{2}\Z$. By Proposition \ref{SinglyGrK0MapLinear}, the functor $X \boxtimes -$ induces a $\frac{1}{2}\Z$--linear map on Grothendieck groups 
\[
[X \boxtimes -]: K_0(\D^b(\Ccrl(n,\Sc))) \to K_0(\D^b(\Ccrl(n,\Sc'))).
\]

\begin{lemma}\label{SingleGrMapFromGensLemma}
Let $X$ be as above. Choose a set of generators $x_{\alpha}$ for $X$, each of which is grading-homogeneous and has unique left and right idempotents $e_L(x_{\alpha})$ and $e_R(x_{\alpha})$. The matrix element for $[X \boxtimes -]$ with incoming basis element $[Q_{\Ib}]$ and outgoing basis element $[Q_{\Ib'}]$ is equal to 
\[
\sum_{x_{\alpha} \textrm{ with } e_R(x_{\alpha}) = \Ib, \,\, e_L(x_{\alpha}) = \Ib'} (-1)^{\Mas(x_{\alpha})} t^{\Alex^{\sing}(x_{\alpha})},
\]
where $\Mas(x_{\alpha}) \in \Z$ is the homological grading of $x_{\alpha}$ and $\Alex^{\sing}(x_{\alpha}) \in \frac{1}{2}\Z$ is the intrinsic grading of $x_{\alpha}$. 
\end{lemma}

\begin{proof}
We can break apart $[X \boxtimes Q_{\Ib}]$ as a sum, inductively, like in the proof of Theorem \ref{SingleAlgDecatBodyThm}. The coefficient of this sum on $[Q_{\Ib'}]$ is a weighted sum over the appropriate generators of $X$ as described.
\end{proof}

\begin{corollary}\label{BddTypeDKSCorr}
Let $X$ be the DA bimodule associated to a crossing. In the singly graded case, the matrix element for $[X \boxtimes -]$ with incoming basis element $[Q_{\Ib}]$ and outgoing basis element $[Q_{\Ib'}]$ is equal to 
\[
\sum_{\textrm{partial K.S. } \xi \textrm{ with } e_R(\xi) = \Ib, \,\, e_L(\xi) = \Ib'} (-1)^{\Mas(\xi)} t^{\Alex^{\sing}(\xi)},
\]
where $\Mas(\xi) \in \Z$ and $\Alex^{\sing}(\xi) \in \frac{1}{2}\Z$ are given by the formulas in Figure \ref{AlexMaslovDefFig}.
\end{corollary}

\begin{proof} This corollary follows from Lemma \ref{SingleGrMapFromGensLemma} and the fact that the generators of Ozsv{\'a}th--Szab{\'o}'s DA bimodules for crossings are partial Kauffman states, with gradings as specified.
\end{proof}

\subsubsection{Decategorifying DA bimodules for crossings: multiply graded case}\label{MultiBraidDecatSect}

Let $X$ be the DA bimodule associated to a crossing with underlying oriented $1$--manifold $W$, and suppose $\partial W = -Y_2 \sqcup Y_1$. The intrinsic grading group of $X$ is $M_W = H^1(W,\partial W; \frac{1}{4}\Z)$. This group comes equipped with homomorphisms from the grading groups $M_{Y_1} = H^0(Y_1; \frac{1}{4}\Z)$ and $M_{Y_2} = H^0(Y_2; \frac{1}{4}\Z)$ of $\Ccrl(n,\Sc)$ and $\Ccrl(n,\Sc')$ respectively. To $X$, we want to associate a $\Z[M_W]$--linear map from
\begin{gather*}
K_0(\D^b_{M_{Y_1}}(\Ccrl(n,\Sc))) \otimes_{\Z[M_{Y_1}]} \Z[M_W] \\
\tag*{\text{to}} K_0(\D^b_{M_{Y_2}}(\Ccrl(n,\Sc))) \otimes_{\Z[M_{Y_2}]} \Z[M_W].
\end{gather*}

By Proposition \ref{MultiplyGrK0MapLinear}, $X \boxtimes -$ induces an additive map
\[
[X \boxtimes -]: K_0(\D^b_{M_{Y_1}}(\Ccrl(n,\Sc))) \to K_0(\D^b_{M_W}(\Ccrl(n,\Sc')))
\]
that intertwines the $M_{Y_1}$--action on the input Grothendieck group and the $M_W$ action on the output Grothendieck group. Equivalently, we get a $\Z[M_W]$--linear map
\[
K_0(\D^b_{M_{Y_1}}(\Ccrl(n,\Sc))) \otimes_{\Z[M_{Y_1}]} \Z[M_W] \to K_0(\D^b_{M_W}(\Ccrl(n,\Sc'))).
\]

Finally, by Corollary \ref{K0ofWIsTensorProd}, we may rewrite the map on Grothendieck groups associated to $X \boxtimes -$ as a $\Z[M_W]$--linear map from
\begin{gather*}
K_0(\D^b_{M_{Y_1}}(\Ccrl(n,\Sc))) \otimes_{\Z[M_{Y_1}]} \Z[M_W] \\
\tag*{\text{to}} K_0(\D^b_{M_{Y_2}}(\Ccrl(n,\Sc'))) \otimes_{\Z[M_{Y_2}]} \Z[M_W].
\end{gather*}

As in the singly graded case, we get formulas for $[X \boxtimes -]$ as sums over generators of $X$:
\begin{lemma}\label{MultiTypeDSumOfGensLemma}
Let $X$ be a finitely generated left bounded DA bimodule over $(\Ccrl(n,\Sc'), \Ccrl(n,\Sc))$, with multiple intrinsic gradings by $M_W$. Choose a set of generators $x_{\alpha}$ for $X$, each of which is grading-homogeneous and has unique left and right idempotents $e_L(x_{\alpha})$ and $e_R(x_{\alpha})$. The matrix element for $[X \boxtimes -]$
with incoming basis element
\[
[Q_{\Ib}] \in K_0(\D^b_{M_1}(\Ccrl(n,\Sc))) \otimes_{\Z[M_{Y_1}]} \Z[M_W] 
\]
and outgoing basis element 
\[
[Q_{\Ib'}] \in K_0(\D^b_{M_2}(\Ccrl(n,\Sc'))) \otimes_{\Z[M_{Y_2}]} \Z[M_W],
\] 
is equal to 
\[
\sum_{x_{\alpha} \textrm{ with } e_R(x_{\alpha}) = \Ib, \,\, e_L(x_{\alpha}) = \Ib'} (-1)^{\Mas(x_{\alpha})} \Alex^{\multi}(x_{\alpha}),
\]
an element of $\Z[M_W]$, where $\Mas(x_{\alpha}) \in \Z$ is the homological grading of $x_{\alpha}$ and $\Alex^{\multi}(x_{\alpha}) \in M_W$ is the intrinsic grading of $x_{\alpha}$.
\end{lemma}

\begin{proof}
The proof is the same as for Lemma \ref{SingleGrMapFromGensLemma}.
\end{proof}

\begin{corollary}\label{MultiKSCorr}
Let $X$ be the DA bimodule for a crossing (or more generally an open braid) $\Gamma$ whose underlying oriented $1$--manifold is $W$, with multi-gradings by $M_W := H^1(W, \partial W; \frac{1}{2}\Z)$. The matrix element for $[X \boxtimes -]$ with incoming basis element 
\[
[Q_{\Ib}] \in K_0(\D^b_{M_{Y_1}}(\Ccrl(n,\Sc))) \otimes_{\Z[M_{Y_1}]} \Z[M_W] 
\]
and outgoing basis element 
\[
[Q_{\Ib'}] \in K_0(\D^b_{M_{Y_2}}(\Ccrl(n,\Sc'))) \otimes_{\Z[M_{Y_2}]} \Z[M_W],
\] 
is equal to 
\[
\sum_{\textrm{partial K.S. } \xi \textrm{ with } e_R(\xi) = \Ib, \,\, e_L(\xi) = \Ib'} (-1)^{\Mas(\xi)} \Alex^{\multi}(\xi),
\]
where $\Mas(\xi) \in \Z$ and $\Alex^{\multi}(\xi) \in M_W$ are defined in Figure \ref{OSzMultiFig} when $\Gamma$ is a single crossing. When $\Gamma$ is an open braid, $\Mas(\xi)$ and $\Alex^{\multi}(\xi)$ are sums of contributions from Figure \ref{OSzMultiFig} over the Kauffman corners of $\xi$.
\end{corollary}

\begin{proof}
This corollary follows from Lemma \ref{MultiTypeDSumOfGensLemma}.
\end{proof}

\subsection{Semisimple modules}\label{SemisimpleSect}

Finally, we describe another decategorification procedure that works for maxima and minima as well as for crossings.

\subsubsection{Triangulated categories of semisimple modules}

For each elementary idempotent $\Ib$ of $\Ccrl(n,\Sc)$, we have a one-dimensional module $S_{\Ib}$ over $\Ccrl(n,\Sc)$. Above, $S_{\Ib}$ was a Type D structure over $\Ibrl(n)$ with zero Type D operation; now we view $S_{\Ib}$ as a dg module over the full algebra $\Ccrl(n,\Sc)$ with zero differential and with all non-idempotent generators acting as zero. 

We will call $\{S_{\Ib}\}$ the ``elementary simple modules'' over $\Ccrl(n,\Sc)$. We call a dg module $S$ over $\Ccrl(n,\Sc)$ ``semisimple'' if it is isomorphic, in the homotopy category of dg modules over $\Ccrl(n,\Sc)$, to either a direct sum or an iterated mapping cone of the elementary simple modules. These two conditions are equivalent because any mapping cone of a morphism between direct sums of elementary simple modules may be viewed, up to isomorphism, as the mapping cone of a zero morphism between two (possibly smaller) sums of elementary simple modules.

\begin{remark}
We will not attempt to define or classify simple modules over $\Ccr(n,\Sc)$ or $\Ccl(n,\Sc)$ in general; we will just work with the elementary ones. In the related context of bordered strands algebras, all simple modules are one-dimensional; see \cite[Section 5]{LOTMappingClass}. Unfortunately, the argument in that paper does not apply to $\Ccr(n,\Sc)$ or $\Ccl(n,\Sc)$ since the augmentation ideals of these algebras are not nilpotent.
\end{remark}

\begin{definition}\label{SingleSSCatsDef}
Let $\SSp^{\sing}(\Ccrl(n,\Sc))$ denote the homotopy category of semisimple dg modules over $\Ccrl(n,\Sc)$, with intrinsic gradings by $\frac{1}{2}\Z$. It is a full triangulated subcategory of the homotopy category of dg modules over $\Ccrl(n,\Sc)$ with $\frac{1}{2}\Z$ intrinsic gradings. 
\end{definition}

Similarly, in the multi-graded case:

\begin{definition}\label{MultiSSCatsDef}
For a zero-manifold $Y$ consisting of $n$ points oriented according to $\Sc$, let $\SSp^{\multi}_Y(\Ccrl(n,\Sc))$ denote the homotopy category of semisimple dg modules over $\Ccrl(n,\Sc)$, with intrinsic gradings by $M_Y = H^0(Y; \frac{1}{4}\Z)$. It is a full triangulated subcategory of the homotopy category of dg modules over $\Ccrl(n,\Sc)$ with intrinsic gradings by $M_Y$.

Let $W$ be an oriented $1$--manifold with $\partial W = -Y_2 \sqcup Y_1$, where $Y_2$ consists of $n$ points oriented according to $\Sc$. Let $\SSp^{\multi}_W(\Ccrl(n,\Sc))$ denote the triangulated category defined as above with $M_Y$ replaced by $M_W := H^1(W,\partial W;\frac{1}{4}\Z)$. 
\end{definition}

\begin{theorem}[cf. Theorem \ref{AlgDecatIntroThm}, Theorem \ref{MultiAlgDecatIntroThm}]\label{SSAlgDecatProp}
The classes $[S_{\Ib}]$ form a basis for $K_0(\SSp^{\sing}(\Ccrl(n,\Sc)))$ over $\Z[t^{\pm 1/2}]$, a basis for $K_0(\SSp^{\multi}_Y(\Ccrl(n,\Sc)))$ over $\Z[H^0(Y;\frac{1}{4}\Z)]$, and a basis for $K_0(\SSp^{\multi}_W(\Ccrl(n,\Sc)))$ over $\Z[H^1(W,\partial W; \frac{1}{4}\Z)]$.
\end{theorem}

\begin{proof}
The triangulated categories 
\[
\SSp^{\sing}(\Ccrl(n,\Sc)), \,\, \SSp^{\multi}_Y(\Ccrl(n,\Sc)), \textrm{ and } \SSp^{\multi}_W(\Ccrl(n,\Sc))
\]
are equivalent to the categories $\D^b(\Ibrl(n))$ in Lemma \ref{SingleGrIdemDecatLemma} and 
\[
\D^b_{M_Y}(\Ibrl(n)), \,\, \D^b_{M_W}(\Ibrl(n))
\] 
in Lemma \ref{MultiIdemDecatLemma}. Thus, their Grothendieck groups are isomorphic (as $\Z$--modules, at first) to the ones computed in Lemma \ref{SingleGrIdemDecatLemma} and Lemma \ref{MultiIdemDecatLemma}. Since the equivalences of categories respect shifts in the intrinsic gradings, the isomorphisms are linear over $\Z[t^{\pm 1/2}]$, $\Z[M_Y]$, and $\Z[M_W]$ respectively.
\end{proof}

\begin{corollary}\label{SSYandWCorr}
We have
\[
K_0(\SSp^{\multi}_W(\Ccrl(n,\Sc))) \cong K_0(\SSp^{\multi}_Y(\Ccrl(n,\Sc))) \otimes_{\Z[H^0(Y;\frac{1}{4}\Z)]} \Z\bigg[H^1\bigg(W,\partial W;\frac{1}{4}\Z\bigg)\bigg]
\]
as modules over $\Z[H^1(W,\partial W; \frac{1}{4}\Z)]$.
\end{corollary}

\subsubsection{Unbounded Type D structures and maps on Grothendieck groups}

Now consider Type D structures that are not necessarily operationally bounded:
\begin{definition}
Let $W'$ be an oriented $1$--manifold with no closed components. Let $\TypeD^{\multi}_{W'}(\Ccrl(n,\Sc))$ denote the homotopy category of finitely generated (but possibly unbounded) Type D structures over $\Ccrl(n,\Sc)$ with intrinsic gradings by $H^1(W', \partial W'; \frac{1}{4}\Z)$.
\end{definition}

Let $W$ be another $1$--manifold such that the composition $W \cup W'$ makes sense and has no closed components. By \cite[Proposition 3.19]{OSzNew}, taking box tensor product with a finitely generated DA bimodule graded by $H^1(W, \partial W; \frac{1}{4}\Z)$ induces a functor
\[
\TypeD^{\multi}_{W'}(\Ccrl(n,\Sc)) \to \TypeD^{\multi}_{W \cup W'}(\Ccrl(m,\Sc')),
\]
even though the boundedness needed for this tensor product to be well-defined \emph{a priori} does not hold. 

In particular, we can take $W' = Y_1 \times I$ where $\partial W = -Y_2 \sqcup Y_1$. In this case, after some natural identifications, we get a functor
\[
\TypeD^{\multi}_{Y_1}(\Ccrl(n,\Sc)) \to \TypeD^{\multi}_{W}(\Ccrl(m,\Sc')).
\]

We do not attempt to compute $K_0$ of $\TypeD^{\multi}_W(\Ccrl(n,\Sc))$, or use the above functor to define a functor from $\SSp^{\multi}_{Y_1}(\Ccrl(n,\Sc))$ into $\SSp^{\multi}_W(\Ccrl(m,\Sc'))$. Rather, given a finitely generated DA bimodule graded by $H^1(W, \partial W, \frac{1}{4}\Z)$ where $W$ has no closed components, we directly construct a map from
\begin{gather*}
K_0(\SSp^{\multi}_{Y_1}(\Ccrl(n, \Sc))) \\
\tag*{\text{to}} K_0(\SSp^{\multi}_{W}(\Ccrl(m, \Sc'))):
\end{gather*}
\begin{definition}\label{InducedDAMapSemisimple}
Let $X$ be a finitely generated DA bimodule over $\Ccrl(m,\Sc')$ and $\Ccrl(n,\Sc)$, graded by $H^1(W, \partial W, \frac{1}{4}\Z)$ where $W$ has no closed components. Define a map 
\[
[X]: K_0(\SSp^{\multi}_{Y_1}(\Ccrl(n, \Sc))) \to K_0(\SSp^{\multi}_W(\Ccrl(m, \Sc')))
\]
as follows: the basis element $[S_{\Ib}]$ of $K_0(\SSp^{\multi}_{Y_1}(\Ccrl(n, \Sc)))$ is sent to
\[
\sum_{\Ib'} \sum_{i \in \Z, \,\, \gamma \in H^1(W, \partial W, \frac{1}{4}\Z)} [\Hom_{H(\Ccrl(m,\Sc')-\dgmod)}(X \boxtimes Q_{\Ib}, \,\, S_{\Ib'}[-i]\{-\gamma\}) \otimes_{\Z} S_{\Ib'}[-i]\{-\gamma\}].
\]
Note that the formula uses $Q_{\Ib}$ rather than $S_{\Ib}$. We view $Q_{\Ib}$ here as a Type D structure rather than an indecomposable projective module. The map $[X]$ is additive, and it intertwines the $\Z[H^0(Y_1; \frac{1}{4}\Z)]$ action on the incoming Grothendieck group with the $\Z[H^1(W, \partial W; \frac{1}{4}\Z)]$ action on the outgoing Grothendieck group via the usual homomorphism 
\[
d^*: \Z\bigg[H^0 \bigg(Y_1; \frac{1}{4}\Z \bigg) \bigg] \to \Z\bigg[H^1 \bigg(W, \partial W; \frac{1}{4}\Z \bigg) \bigg].
\]

Under certain conditions, we may also define a $\Z[t^{\pm 1/2}]$--linear map 
\[
[X]: K_0(\SSp^{\sing}(\Ccrl(n, \Sc))) \to K_0(\SSp^{\sing}(\Ccrl(m, \Sc')))
\]
by sending $[\widetilde{S}_{\Ib}]$ to
\[
\sum_{\Ib'} \sum_{i \in \Z, \,\, \gamma \in H^1(W, \partial W; \frac{1}{4}\Z)} [\Hom_{H(\Ccrl(m,\Sc')-\dgmod)}(X \boxtimes Q_{\Ib}, S_{\Ib'}[-i]\{-\gamma\}) \otimes_{\Z} \widetilde{S}_{\Ib'}[-i]\{-\gamma(W)\}],
\]
where $\widetilde{S}_{\Ib}$ and $\widetilde{S}_{\Ib'}$ are $S_{\Ib}$ and $S_{\Ib'}$ viewed as objects of the singly-graded categories $\SSp^{\sing}(\Ccrl(n,\Sc))$ and $\SSp^{\sing}(\Ccrl(m,\Sc'))$, and $\gamma(W) \in \frac{1}{4}\Z$ is the pairing of $\gamma$ with the fundamental class of $W$. For this map to make sense, we must make the assumption that $\gamma(W) \in \frac{1}{2}\Z \subset \frac{1}{4}\Z$ whenever 
\[
\Hom_{H(\Ccrl(m,\Sc')-\dgmod)}(X \boxtimes Q_{\Ib}, S_{\Ib'}[-i]\{-\gamma\}) \neq 0.
\]
\end{definition}

\begin{remark}
As motivation for Definition \ref{InducedDAMapSemisimple}, we can think of objects of the category $\SSp^{\multi}_W(\Ccrl(n,\Sc))$ as ``test functions.'' Objects of $\TypeD^{\multi}_W(\Ccrl(n,\Sc))$ define ``distributions'' via the pairing
\[
\Hom_{H(\Ccrl(n,\Sc)-\dgmod)}(-, -).
\]
Furthermore, since $\SSp^{\multi}_W(\Ccrl(n,\Sc))$ is self-dual in an appropriate sense, we may identify these ``distributions'' with objects of $\SSp^{\multi}_W(\Ccrl(n,\Sc))$. (We see that we discard information when we view unbounded Type D structures distributionally in this way.) We will not attempt to make this reasoning precise, and we will always go directly from objects of $\TypeD^{\multi}_W(\Ccrl(n,\Sc))$ to classes in $K_0(\SSp^{\multi}_W(\Ccrl(n,\Sc)))$.
\end{remark}

Following Section \ref{MultiBraidDecatSect}, we can view $[X]$ as a $\Z[H^1(W, \partial W; \frac{1}{4}\Z)]$--linear map from 
\begin{gather*}
K_0(\SSp^{\multi}_{Y_1}(\Ccrl(n, \Sc))) \otimes_{\Z[H^0(Y_1; \frac{1}{4}\Z)]} \Z\bigg[H^1\bigg(W,\partial W \frac{1}{4}\Z\bigg)\bigg] \\
\tag*{\text{to}} K_0(\SSp^{\multi}_W(\Ccrl(m, \Sc'))).
\end{gather*}
Using Corollary \ref{SSYandWCorr}, we arrive at a $\Z[H^1(W, \partial W; \frac{1}{4}\Z)]$--linear map $[X]$ from 
\begin{gather*}
K_0(\SSp^{\multi}_{Y_1}(\Ccrl(n, \Sc))) \otimes_{\Z[H^0(Y_1; \frac{1}{4}\Z)]} \Z\bigg[H^1\bigg(W,\partial W; \frac{1}{4}\Z\bigg)\bigg] \\
\tag*{\text{to}} K_0(\SSp^{\multi}_{Y_2}(\Ccrl(m, \Sc'))) \otimes_{\Z[H^0(Y_2; \frac{1}{4}\Z)]} \Z\bigg[H^1\bigg(W,\partial W; \frac{1}{4}\Z\bigg)\bigg].
\end{gather*}

\subsubsection{Expanding maps as sums over generators}
As in Section \ref{SingleBraidDecatSect} and \ref{MultiBraidDecatSect}, we may write the maps $[X]$ associated to DA bimodules as sums over generators of $X$. When $X$ is Ozsv{\'a}th--Szab{\'o}'s DA bimodule for a tangle diagram, these generators are partial Kauffman states.

\begin{lemma}\label{SSSumOfGensLemma}
Let $X$ be as in Definition \ref{InducedDAMapSemisimple} and choose a set of generators $x_{\alpha}$ for $X$, each of which is grading-homogeneous and has unique left and right idempotents $e_L(x_{\alpha})$ and $e_R(x_{\alpha})$. In the multi-graded case with grading group $H^1(W, \partial W; \frac{1}{4}\Z)$ for a $1$--manifold $W$, the matrix element for $[X]$ with incoming basis element $[S_{\Ib}]$ and outgoing basis element $[S_{\Ib'}]$ is equal to 
\[
\sum_{x_{\alpha} \textrm{ with } e_R(x_{\alpha}) = \Ib, \,\, e_L(x_{\alpha}) = \Ib'} (-1)^{\Mas(x_{\alpha})} \Alex^{\multi}(x_{\alpha}),
\]
where $\Mas(x_{\alpha}) \in \Z$ is the homological grading of $x_{\alpha}$ and 
\[
\Alex^{\multi}(x_{\alpha}) \in H^1\bigg(W, \partial W; \frac{1}{4}\Z \bigg)
\]
is the intrinsic multi-grading of $x_{\alpha}$. 

In the singly graded case, the matrix element for $[X]$ with incoming basis element $[\widetilde{S}_{\Ib}]$ and outgoing basis element $[\widetilde{S}_{\Ib'}]$ is equal to 
\[
\sum_{x_{\alpha} \textrm{ with } e_R(x_{\alpha}) = \Ib, \,\, e_L(x_{\alpha}) = \Ib'} (-1)^{\Mas(x_{\alpha})} t^{\Alex^{\sing}(x_{\alpha})},
\]
where $\Alex^{\sing}(x_{\alpha}) \in \frac{1}{2}\Z$ is the single intrinsic grading of $x_{\alpha}$. 
\end{lemma}

\begin{proof}
In the multiply graded case, we may write $[X]([S_{\Ib}])$ as 
\[
\sum_{\Ib'} \sum_{i \in \Z, \,\, \gamma \in H^1(W, \partial W; \frac{1}{4}\Z)} (-1)^i ( \rk \Hom_{H(\Ccrl(m,\Sc'))-\dgmod}(X \boxtimes Q_{\Ib}, \,\, S_{\Ib'}[-i]\{-\gamma\}))\gamma \cdot [S_{\Ib'}].
\]
Thus, the coefficient on $[S_{\Ib'}]$ is 
\[
\sum_{i,\gamma} (-1)^i (\rk \Hom_{H(\Ccrl(m,\Sc'))-\dgmod}(X \boxtimes Q_{\Ib}, \,\, S_{\Ib'}[-i]\{-\gamma\}))\gamma.
\]
The box tensor product $X \boxtimes Q_{\Ib}$ is a projective left dg module over $\Ccrl(m,\Sc')$ with generators 
\[
\{x_{\alpha} \,\,\vert\,\, e_R(x_{\alpha}) = \Ib\}.
\]
There is a basis for $\Hom_{H(\Ccrl(m,\Sc')-\dgmod)}(X \boxtimes Q_{\Ib}, S_{\Ib'}[-i]\{-\gamma\})$ consisting of homomorphisms that send one generator $x_{\alpha}$ of $X$ with
\begin{itemize}
\item $e_R(x_{\alpha}) = \Ib$,
\item $e_L(x_{\alpha}) = \Ib'$,
\item $\Mas(x_{\alpha}) = i$,
\item $\Alex^{\multi}(x_{\alpha}) = \gamma$
\end{itemize}
to the generator of $S_{\Ib'}[-i]\{-\gamma\}$ and send all other generators of $X$ to zero. We see that the sum defining the matrix entry above is equal to a sum over generators $x_{\alpha}$ of $X$ with $e_R(x_{\alpha}) = \Ib$ and $e_L(x_{\alpha}) = \Ib'$, weighted by the homological and intrinsic gradings of the generators, as in the statement of the proposition. 

The singly-graded case follows from the multiply-graded case by pairing classes in $H^1(W, \partial W; \frac{1}{4}\Z)$ with the fundamental class of $W$, assuming as usual that this pairing produces elements of $\frac{1}{2}\Z \subset \frac{1}{4}\Z$ for all $\gamma$ such that $\gamma = \Alex^{\multi}(x_{\alpha})$ for some generator $x_{\alpha}$.
\end{proof}

\begin{corollary}\label{PartialKSDecatLemma}
In Definition \ref{InducedDAMapSemisimple}, let $X$ be Ozsv{\'a}th--Szab{\'o}'s DA bimodule for a tangle with no closed components. In the multiply graded setting, the matrix element for $[X]$ with incoming basis element $[S_{\Ib}]$ and outgoing basis element $[S_{\Ib'}]$ is equal to 
\[
\sum_{\textrm{partial K.S. } \xi \textrm{ with } e_R(\xi) = \Ib, \,\, e_L(\xi) = \Ib'} (-1)^{\Mas(\xi)} \Alex^{\multi}(\xi),
\]
an element of $\Z[H^1(W, \partial W; \frac{1}{4}\Z)]$.

In the singly graded setting, the matrix element for $[X]$ with incoming basis element $[\widetilde{S}_{\Ib}]$ and outgoing basis element $[\widetilde{S}_{\Ib'}]$ is equal to 
\[
\sum_{\textrm{partial K.S. } \xi \textrm{ with } e_R(\xi) = \Ib, \,\, e_L(\xi) = \Ib'} (-1)^{\Mas(\xi)} t^{\Alex^{\sing}(\xi)},
\]
an element of $\Z[t^{\pm 1/2}]$.
\end{corollary}

\subsubsection{Invariance of the map on Grothendieck groups}

Suppose we have two diagrams representing the same tangle. Ozsv{\'a}th and Szab{\'o} show that the corresponding DA bimodules are quasi-isomorphic. In the above setup, quasi-isomorphism invariance is not enough to guarantee that the maps induced by the DA bimodules using Definition \ref{InducedDAMapSemisimple} are the same.

In fact, Corollary \ref{PartialKSDecatLemma} gives us a simple expression for the maps in terms of Kauffman states, which transform predictably under Reidemeister moves, so we could show the maps are tangle invariants directly. However, this approach requires an involved computation. Instead, we will make use of the comparison results in Section \ref{DADecatSect} and the Reidemeister invariance of representation-theoretic maps (in the singly-graded case) and Viro's maps (in the multiply graded case) to show that the maps $[X]$ of Definition \ref{InducedDAMapSemisimple} are independent of the diagram representing the tangle; see Corollary \ref{InvarianceCorr} below.

\section{Computing the decategorifications of tangle bimodules}\label{DADecatSect}

We now use Corollary \ref{BddTypeDKSCorr}, Corollary \ref{MultiKSCorr}, and Corollary \ref{PartialKSDecatLemma} to compute the linear maps on Grothendieck groups associated to Ozsv{\'a}th--Szab{\'o}'s DA bimodules for crossings, maxima, and minima. For crossings, note that the maps from Section \ref{BddTypeDSect} and Section \ref{SemisimpleSect} have the same formulas in terms of partial Kauffman states, so we will not need to distinguish between them.

\subsection{DA bimodules for crossings}\label{ExplicitDecatCrossingSect}

\subsubsection{Generic case} 

First we describe the generic case $1 < i < n$; the special cases $i = 1$ and $i = n$ will be discussed in Section \ref{CompareCrossingsSpecialSect}.

Our goal is to compute the matrix element of the map $[\Pc^i_{r/l} \boxtimes -]$ corresponding to any incoming idempotent $\Ib_{\y}$ of $\Ccrl(n,\Sc)$ and any outgoing idempotent $\Ib_{\x}$ of $\Cc_{r/l}(n,\Sc')$. In the singly-graded case, for each basis element
\[
[Q_{\Ib_{\y}}] \in K_0(\D^b(\Ccrl(n,\Sc))),
\]
we want to express $[\Pc^i_{r/l} \boxtimes Q_{\Ib_{\y}}]$ in terms of the basis elements 
\[
\{[Q_{\Ib_{\x}}]\} \in K_0(\D^b(\Ccrl(n,\Sc'))).
\]
The matrix we get is also the matrix for $[\Pc^i_{r/l}]$ in terms of the bases $\{S_{\Ib_{\y}}\}$ and $\{S_{\Ib_{\x}}\}$ of $K_0(\SSp^{\sing}(\Ccrl(n,\Sc)))$ and $K_0(\SSp^{\sing}(\Ccrl(n,\Sc')))$. 

Consider an idempotent $\Ib_{\y}$ of $\Cc_{r/l}(n,\Sc)$ with dots in regions $i_1, \ldots, i_k$ (i.e. $\x = \{i_1,\ldots,i_k\}$). Relative to the three regions $\{i-1,i,i+1\}$, $\Ib_{\y}$ takes the form of one of eight ``local idempotents,'' namely 
\begin{align*}
&\{ \varnothing = \hskip 0.221in | \hskip 0.221in | \hskip 0.221in, \\
&A = \hskip 0.09in \Bigcdot | \hskip 0.221in | \hskip 0.221in, \quad B = \hskip 0.2in | \Bigcdot | \hskip 0.221in, \quad C = \hskip 0.2in | \hskip 0.221in | \Bigcdot, \\
&AB = \hskip 0.09in \Bigcdot | \Bigcdot | \hskip 0.221in, \quad AC = \Bigcdot | \hskip 0.221in | \Bigcdot, \quad BC = \hskip 0.221in | \Bigcdot | \Bigcdot, \\
&ABC = \hskip 0.09in \Bigcdot | \Bigcdot | \Bigcdot \}.
\end{align*}

\begin{proposition}\label{DecatCorr}
Suppose strands $i$ and $i+1$ are oriented \OrUU. The matrix for $[\Pc^i_{r/l} \boxtimes -]$ and $[\Pc^i_{r/l}]$ in terms of local idempotents is
\[
\kbordermatrix{
& \varnothing & A & B & C & AB & AC & BC & ABC \\
\varnothing & t^{1/2} &  &  &  &  &  &  &  \\
A &  & t^{1/2} & 0 & 0 &  &  &  &  \\
B &  & 1 & -t^{-1/2} & 1 &  &  &  &  \\
C &  & 0 & 0 & t^{1/2} &  &  &  &  \\
AB &  &  &  &  & -t^{-1/2} & 1 & 0 &  \\
AC &  &  &  &  & 0 & t^{1/2} & 0 &  \\
BC &  &  &  &  & 0 & 1 & -t^{-1/2} &  \\
ABC &  &  &  &  &  &  &  & -t^{-1/2}
}.
\]

If strands $i$ and $i+1$ are oriented \OrDU, then the matrix for $[\Pc^i_{r/l} \boxtimes -]$ and $[\Pc^i_{r/l}]$ is
\[
\kbordermatrix{
& \varnothing & A & B & C & AB & AC & BC & ABC \\
\varnothing & 1 &  &  &  &  &  &  &  \\
A &  & 1 & 0 & 0 &  &  &  &  \\
B &  & -t^{1/2} & 1 & t^{-1/2} &  &  &  &  \\
C &  & 0 & 0 & 1&  &  &  &  \\
AB &  &  &  &  & 1 & t^{-1/2} & 0 &  \\
AC &  &  &  &  & 0 & 1 & 0 &  \\
BC &  &  &  &  & 0 & -t^{1/2} & 1 &  \\
ABC &  &  &  &  &  &  &  & 1
}.
\]

If strands $i$ and $i+1$ are oriented \OrUD, then the matrix for $[\Pc^i_{r/l} \boxtimes -]$ and $[\Pc^i_{r/l}]$ is
\[
\kbordermatrix{
& \varnothing & A & B & C & AB & AC & BC & ABC \\
\varnothing & 1 &  &  &  &  &  &  &  \\
A &  & 1 & 0 & 0 &  &  &  &  \\
B &  & t^{-1/2} & 1 & -t^{1/2} &  &  &  &  \\
C &  & 0 & 0 & 1 &  &  &  &  \\
AB &  &  &  &  & 1 & -t^{1/2} & 0 &  \\
AC &  &  &  &  & 0 & 1 & 0 &  \\
BC &  &  &  &  & 0 & t^{-1/2} & 1  &  \\
ABC &  &  &  &  &  &  &  & 1
}.
\]

Finally, if strands $i$ and $i+1$ are oriented \OrDD, then the matrix for $[\Pc^i_{r/l} \boxtimes -]$ and $[\Pc^i_{r/l}]$ is
\[
\kbordermatrix{
& \varnothing & A & B & C & AB & AC & BC & ABC \\
\varnothing & -t^{-1/2} &  &  &  &  &  &  &  \\
A &  & -t^{-1/2} & 0 & 0 &  &  &  &  \\
B &  & 1 & t^{1/2} & 1 &  &  &  &  \\
C &  & 0 & 0 & -t^{-1/2} &  &  &  &  \\
AB &  &  &  &  & t^{1/2} & 1 & 0 &  \\
AC &  &  &  &  & 0 & -t^{-1/2} & 0 &  \\
BC &  &  &  &  & 0 & 1 & t^{1/2} &  \\
ABC &  &  &  &  &  &  &  & t^{1/2}
}.
\]
\end{proposition}

\begin{proof}
We use Corollary \ref{BddTypeDKSCorr}. Given incoming and outgoing idempotents, the relevant partial Kauffman states are shown in Figure \ref{PartialKSFig}. In each of the four orientation cases, we can read off the Alexander and Maslov gradings from Figure \ref{AlexMaslovDefFig}; these gradings determine the matrix entries. The proof is thus an entry-by-entry check which we omit for space. Compared with the proofs in Section \ref{MapsInModifiedBasisSect}, this check is straightforward.
\end{proof}

We will not write out the matrices for $[\Nc^i_{r/l} \boxtimes -]$ and $[\Nc^i_{r/l}]$; they may be obtained as inverses of the matrices in Proposition \ref{DecatCorr}.

In the multigraded case, for a basis element 
\[
[Q_{\Ib_{\y}}] \in K_0(\D^b_{Y_1}(\Ccrl(n,\Sc))) \otimes_{\Z[H^0(Y_1;\frac{1}{4}\Z)]} \Z\bigg[H^1 \bigg(W, \partial W; \frac{1}{4}\Z \bigg) \bigg],
\]
we want to express $[\Pc^i_{r/l} \boxtimes Q_{\Ib_{\y}}]$ in terms of the basis elements 
\[
\{[Q_{\Ib_{\x}}]\} \in K_0(\D^b_{Y_2}(\Ccrl(n,\Sc'))) \otimes_{\Z[H^0(Y_2;\frac{1}{4}\Z)]} \Z\bigg[H^1 \bigg(W, \partial W; \frac{1}{4}\Z \bigg) \bigg].
\]
The resulting matrix is also the matrix for $[\Pc^i_{r/l}]$ in terms of the bases $\{S_{\Ib_{\x}}\}$ and $\{S_{\Ib_{\x}}\}$ of 
\begin{gather*}
K_0(\SSp^{\multi}_{Y_1}(\Ccrl(n,\Sc))) \otimes_{\Z[H^0(Y_1;\frac{1}{4}\Z)]} \Z \bigg[H^1 \bigg(W, \partial W; \frac{1}{4}\Z \bigg) \bigg] \\
\tag*{\text{and}} K_0(\SSp^{\multi}_{Y_2}(\Ccrl(n,\Sc'))) \otimes_{\Z[H^0(Y_2;\frac{1}{4}\Z)]} \Z\bigg[H^1\bigg(W, \partial W; \frac{1}{4}\Z\bigg) \bigg]
\end{gather*}
respectively. 

\begin{proposition}\label{MultiDecatComputation}
Suppose strands $i$ and $i+1$ are oriented \OrUU. The matrix for $[\Pc^i_{r/l} \boxtimes -]$ in terms of local idempotents, in the multiply graded setting, is
\[
\scalebox{0.85}{\parbox{\linewidth}{
\kbordermatrix{
& \varnothing & A & B & C & AB & AC & BC & ABC \\
\varnothing & \sigma_i^{1/4} \sigma_{i+1}^{1/4} &  &  &  &  &  &  &  \\
A &  & \sigma_i^{1/4} \sigma_{i+1}^{1/4} & 0 & 0 &  &  &  &  \\
B &  & \sigma_i^{-1/4} \sigma_{i+1}^{1/4} & -\sigma_i^{-1/4} \sigma_{i+1}^{-1/4} & \sigma_i^{1/4} \sigma_{i+1}^{-1/4} &  &  &  &  \\
C &  & 0 & 0 & \sigma_i^{1/4} \sigma_{i+1}^{1/4} &  &  &  &  \\
AB &  &  &  &  & -\sigma_i^{-1/4} \sigma_{i+1}^{-1/4} & \sigma_i^{1/4} \sigma_{i+1}^{-1/4} & 0 &  \\
AC &  &  &  &  & 0 & \sigma_i^{1/4} \sigma_{i+1}^{1/4} & 0 &  \\
BC &  &  &  &  & 0 & \sigma_i^{-1/4} \sigma_{i+1}^{1/4} & -\sigma_i^{-1/4} \sigma_{i+1}^{-1/4} &  \\
ABC &  &  &  &  &  &  &  & -\sigma_i^{-1/4} \sigma_{i+1}^{-1/4}
}.
}}
\]

If strands $i$ and $i+1$ are oriented \OrDU, then the matrix for $[\Pc^i_{r/l} \boxtimes -]$ is
\[
\scalebox{0.85}{\parbox{\linewidth}{
\kbordermatrix{
& \varnothing & A & B & C & AB & AC & BC & ABC \\
\varnothing & \sigma_i^{-1/4} \sigma_{i+1}^{1/4} &  &  &  &  &  &  &  \\
A &  & \sigma_i^{-1/4} \sigma_{i+1}^{1/4} & 0 & 0 &  &  &  &  \\
B &  & -\sigma_i^{1/4} \sigma_{i+1}^{1/4} & \sigma_i^{1/4} \sigma_{i+1}^{-1/4} & \sigma_i^{-1/4} \sigma_{i+1}^{-1/4} &  &  &  &  \\
C &  & 0 & 0 & \sigma_i^{-1/4} \sigma_{i+1}^{1/4} &  &  &  &  \\
AB &  &  &  &  & \sigma_i^{1/4} \sigma_{i+1}^{-1/4} & \sigma_i^{-1/4} \sigma_{i+1}^{-1/4} & 0 &  \\
AC &  &  &  &  & 0 & \sigma_i^{-1/4} \sigma_{i+1}^{1/4} & 0 &  \\
BC &  &  &  &  & 0 & -\sigma_i^{1/4} \sigma_{i+1}^{1/4} & \sigma_i^{1/4} \sigma_{i+1}^{-1/4} &  \\
ABC &  &  &  &  &  &  &  & \sigma_i^{1/4} \sigma_{i+1}^{-1/4}
}.
}}
\]

If strands $i$ and $i+1$ are oriented \OrUD, then the matrix for $[\Pc^i_{r/l} \boxtimes -]$ is
\[
\scalebox{0.85}{\parbox{\linewidth}{
\kbordermatrix{
& \varnothing & A & B & C & AB & AC & BC & ABC \\
\varnothing & \sigma_i^{1/4} \sigma_{i+1}^{-1/4} &  &  &  &  &  &  &  \\
A &  & \sigma_i^{1/4} \sigma_{i+1}^{-1/4} & 0 & 0 &  &  &  &  \\
B &  & \sigma_i^{-1/4} \sigma_{i+1}^{-1/4} & \sigma_i^{-1/4} \sigma_{i+1}^{1/4} & -\sigma_i^{1/4} \sigma_{i+1}^{1/4} &  &  &  &  \\
C &  & 0 & 0 & \sigma_i^{1/4} \sigma_{i+1}^{-1/4} &  &  &  &  \\
AB &  &  &  &  & \sigma_i^{-1/4} \sigma_{i+1}^{1/4} & -\sigma_i^{1/4} \sigma_{i+1}^{1/4} & 0 &  \\
AC &  &  &  &  & 0 & \sigma_i^{1/4} \sigma_{i+1}^{-1/4} & 0 &  \\
BC &  &  &  &  & 0 & \sigma_i^{-1/4} \sigma_{i+1}^{-1/4} & \sigma_i^{-1/4} \sigma_{i+1}^{1/4}  &  \\
ABC &  &  &  &  &  &  &  & \sigma_i^{-1/4} \sigma_{i+1}^{1/4}
}.
}}
\]

Finally, if strands $i$ and $i+1$ are oriented \OrDD, then the
matrix for $[\Pc^i_{r/l} \boxtimes -]$ is
\[
\scalebox{0.85}{\parbox{\linewidth}{
\kbordermatrix{
& \varnothing & A & B & C & AB & AC & BC & ABC \\
\varnothing & -\sigma_i^{-1/4} \sigma_{i+1}^{-1/4} &  &  &  &  &  &  &  \\
A &  & -\sigma_i^{-1/4} \sigma_{i+1}^{-1/4} & 0 & 0 &  &  &  &  \\
B &  & \sigma_i^{1/4} \sigma_{i+1}^{-1/4} & \sigma_i^{1/4} \sigma_{i+1}^{1/4} & \sigma_i^{-1/4} \sigma_{i+1}^{1/4} &  &  &  &  \\
C &  & 0 & 0 & -\sigma_i^{-1/4} \sigma_{i+1}^{-1/4} &  &  &  &  \\
AB &  &  &  &  & \sigma_i^{1/4} \sigma_{i+1}^{1/4} & \sigma_i^{-1/4} \sigma_{i+1}^{1/4} & 0 &  \\
AC &  &  &  &  & 0 & -\sigma_i^{-1/4} \sigma_{i+1}^{-1/4} & 0 &  \\
BC &  &  &  &  & 0 & \sigma_i^{1/4} \sigma_{i+1}^{-1/4} & \sigma_i^{1/4} \sigma_{i+1}^{1/4} &  \\
ABC &  &  &  &  &  &  &  & \sigma_i^{1/4} \sigma_{i+1}^{1/4}
}.
}}
\]
\end{proposition}

\begin{proof}
We use Corollary \ref{MultiKSCorr}. Given incoming and outgoing idempotents, the relevant partial Kauffman states are shown in Figure \ref{PartialKSFig} as before. In each of the four orientation cases, we can read off the Alexander and Maslov gradings from Figure \ref{OSzMultiFig}; these gradings determine the matrix entries.
\end{proof}

\subsubsection{Special cases}\label{CompareCrossingsSpecialSect}
Now we consider the cases $i = 1$ and $i = n$. When $i = 1$, $\Pc^1_l$ works no differently than for generic $i$. However, we need to discuss $\Pc^1_r$ separately.

As above, we want to compute the matrix element of the map $[\Pc^1_r \boxtimes -]$ corresponding to any incoming idempotent $\Ib_{\y}$ and outgoing idempotent $\Ib_{\x}$ of $\Cc_r(n,\Sc)$. Let $\Ib_{\y}$ be an idempotent of $\Cc_r(n,\Sc)$, where $\y \subset [1,\ldots,n]$. Relative to the three regions $\{0,1,2\}$, $\Ib_{\y}$ has local form $\varnothing$, $B$, $C$, or $BC$ as labeled above.

\begin{proposition}
In terms of local idempotents, the matrix for $[\Pc^1_r \boxtimes -]$ is the $4 \times 4$ submatrix of the relevant $8 \times 8$ matrix from Proposition \ref{DecatCorr} or Proposition \ref{MultiDecatComputation}, in which only the entries whose row-label and column-label are $\varnothing$, $B$, $C$, or $BC$ are considered.
\end{proposition}

\begin{proof} 
This proposition follows from Corollary \ref{BddTypeDKSCorr} in the singly graded case and Corollary \ref{MultiKSCorr} in the multi-graded case. 
\end{proof}

Similarly, when $i = n$, $\Pc^n_r$ works no differently than for generic $i$, but we need to discuss $\Pc^n_l$ separately. Let $\Ib_{\y}$ be an idempotent of $\Cc_l(n,\Sc)$, where $\y \subset [0,\ldots,n-1]$. Relative to the three regions $\{n-2,n-1,n\}$, $\Ib_{\y}$ has local form $\varnothing$, $A$, $B$, or $AB$.

\begin{proposition}
In terms of local idempotents, the matrix for $[\Pc^n_l \boxtimes -]$ is the $4 \times 4$ submatrix of the relevant $8 \times 8$ matrix from Proposition \ref{DecatCorr} or Proposition \ref{MultiDecatComputation}, in which only the entries whose row-label and column-label are $\varnothing$, $A$, $B$, or $AB$ are considered.
\end{proposition}

\begin{proof}
Again, this proposition follows from Corollary \ref{BddTypeDKSCorr} in the singly graded case and Corollary \ref{MultiKSCorr} in the multi-graded case.
\end{proof}

\subsection{DA bimodules for maxima and minima}\label{ExplicitDecatMaxMinSect}
\begin{figure} \centering
\includegraphics[scale=0.625]{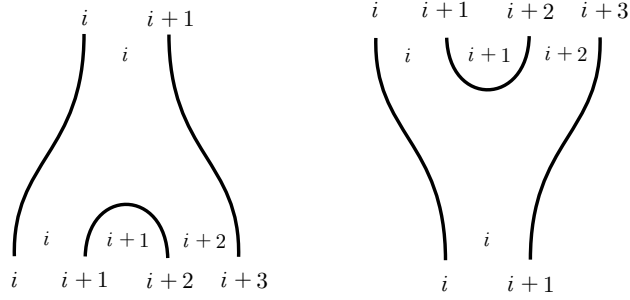}
\caption{Labels on endpoints and regions for maxima and minima.}
\label{MinMaxLabelsFig}
\end{figure}

Now we decategorify $\Omega^{i+1}_{r/l}$ and $\mho^{i+1}_{r/l}$. We must use the decategorification scheme of Section \ref{SemisimpleSect}; we want to compute the maps $[\Omega^{i+1}_{r/l}]$ and $[\mho^{i+1}_{r/l}]$. See Figure \ref{MinMaxLabelsFig} for a summary of the labels on endpoints and regions. 

\subsubsection{Maxima}
We first consider the generic case $1 \leq i \leq n-1$, and then cover the special cases $i = 0$ and $i = n$. Let $\Ib_{\y}$ be an idempotent of $\Cc_{r/l}(n,\Sc)$. We either have $i \notin \y$, in which case we say $\Ib_{\y}$ has local form $\varnothing$, or $i \in \y$, in which case we say $\Ib_{\y}$ has local form $A$. For ``outgoing idempotents'' $\Ib_{\x}$ of $\Cc_{r/l}(n+2,\Sc')$, $\Ib_{\x}$ may have local form $\varnothing$, $A$, $B$, $C$, $AB$, $AC$, $BC$, or $ABC$ depending on the restriction of $\x$ to regions $\{i,i+1,i+2\}$. 

\begin{proposition}\label{OmegaDecatProp}
Let $\Ib_{\y}$ be an idempotent of $\Ccrl(n,\Sc)$, and let $1 \leq i \leq n-1$.
\begin{itemize}
\item
If the local form of $\Ib_{\y}$ is $\varnothing$, then $\Omega^{i+1}_{r/l}$ has one generator with incoming idempotent $\Ib_{\y}$. This generator has bidegree $(0,0)$ and outgoing idempotent $\Ib_{\x}$, where
\[
\x := (\y \cap [0,\ldots,i-1]) \cup \{i+1\} \cup \{j+2: j \in \y \cap [i+1,\ldots,n]\}.
\]

\item
If the local form of $\Ib_{\y}$ is $A$, then $\Omega^{i+1}_{r/l}$ has two generators with incoming idempotent $\Ib_{\y}$. Both generators have bidegree $(0,0)$; one has outgoing idempotent $\Ib_{\x'}$, where
\[
\x' := (\y \cap [0,\ldots,i-1]) \cup \{i,i+1\} \cup \{j+2: j \in \y \cap [i+1,\ldots,n]\},
\]
and the other has outgoing idempotent $\Ib_{\x''}$, where
\[
\x'' := (\y \cap [0,\ldots,i-1]) \cup \{i+1,i+2\} \cup \{j+2: j \in \y \cap [i+1,\ldots,n]\}.
\]
\end{itemize}
Thus, in terms of local idempotents, the matrix for $[\Omega^{i+1}_{r/l}]$ is
\[
\kbordermatrix{
& \varnothing & A  \\
\varnothing & 0 & 0   \\
A & 0 & 0\\
B & 1  & 0\\
C & 0 & 0\\
AB & 0 & 1 \\
AC & 0 & 0\\
BC & 0 & 1 \\
ABC & 0 & 0
},
\]
in either the singly graded or the multiply graded case.
\end{proposition}

\begin{proof}
This proposition follows from Corollary \ref{PartialKSDecatLemma}, as well as the definition of $\Omega^{i+1}$ in \cite[Section 8]{OSzNew}. The grading is discussed \cite[Section 7.1]{OSzNew} in the context of DD bimodules; the gradings remain zero in the DA versions.
\end{proof}

When $i = 0$ and we are considering $\Cc_r$ rather than $\Cc_l$, incoming idempotents $\Ib_{\y}$ of $\Cc_r(n,\Sc)$ must have local form $\varnothing$. Outgoing idempotents $\Ib_{\x}$ of $\Cc_l(n+2,\Sc')$ must have local form $\varnothing$, $B$, $C$, or $BC$.

When $i = n$ and we are considering $\Cc_l$, incoming idempotents $\Ib_{\y}$ of $\Cc_l(n,\Sc)$ must also have local form $\varnothing$. Outgoing idempotents $\Ib_{\x}$ of $\Cc_l(n+2,\Sc')$ must have local form $\varnothing$, $A$, $B$, or $AB$.

\begin{corollary}
The matrices for $[\Omega^1_r]$ and $[\Omega^{n+1}_l]$ are submatrices of the relevant matrices from Proposition \ref{OmegaDecatProp}. Only entries whose row-label and column label are $\varnothing$, $B$, $C$, or $BC$ (respectively $\varnothing$, $A$, $B$, or $AB$) are included in the submatrices.
\end{corollary}

\subsubsection{Minima}
We deal with $\mho^1_{r/l}$ first, and then consider the general case. 

Let $\Ib_{\y}$ be an idempotent of $\Cc_{r/l}(n+2,\Sc)$. Relative to the regions $\{0,1,2\}$, $\Ib_{\y}$ has local form $\varnothing$, $A$, $B$, $C$, $AB$, $AC$, $BC$, or $ABC$. If we are using $\Ccr(n+2,\Sc)$, then $\Ib_{\y}$ has local form $\varnothing$, $B$, $C$, or $BC$. Outgoing idempotents $\Ib_{\x}$ of $\Ccl(n,\Sc')$ may have local form $\varnothing$ or $A$ depending on the restriction of $\x$ to region $i$; outgoing idempotents of $\Ccr(n,\Sc')$ may only have local form $\varnothing$.

\begin{proposition}\label{Mho1Decat}
Let $\Ib_{\y}$ be an idempotent of $\Ccrl(n+2,\Sc)$.
\begin{itemize}
\item
If the local form of $\Ib_{\y}$ is $A$, then $\mho^1_l$ has one generator with incoming idempotent $\Ib_{\y}$. This generator has bidegree $(0,0)$ and outgoing idempotent $\Ib_{\x}$, where
\[
\x :=  \{j-2: j \in \y \cap [3,\ldots,n]\}.
\]

\item
If the local form of $\Ib_{\y}$ is $C$, then $\mho^1_{r/l}$ has one generator with incoming idempotent $\Ib_{\y}$. This generator has bidegree $(0,0)$ and outgoing idempotent $\Ib_{\x}$, where
\[
\x := \{j-2: j \in \y \cap [3,\ldots,n]\}
\]
as in the previous case. 

\item
If the local form of $\Ib_{\y}$ is $AC$, then $\mho^1_l$ has one generator with incoming idempotent $\Ib_{\y}$. This generator has bidegree $(0,0)$ and outgoing idempotent $\Ib_{\x}$, where
\[
\x := \{0\} \cup \{j-2: j \in \y \cap [3,\ldots,n]\}.
\]

\item
For all other local forms of $\Ib_{\y}$, $\mho^1_{r/l}$ has no corresponding generators.

\end{itemize}
Thus, the matrix for $[\mho^1_l]$ in terms of local idempotents is
\[
\kbordermatrix{
& \varnothing & A & B & C & AB & AC & BC & ABC \\
\varnothing & 0 & 1 & 0 & 1 & 0 & 0 & 0 & 0 \\
A & 0 & 0 & 0 & 0 & 0 & 1 & 0 & 0 \\
},
\]
in both the singly graded and the multiply graded cases, except when $n = 0$ in which case we get the $1 \times 4$ submatrix obtained by discarding all rows and columns except those labeled $\varnothing$, $A$, $B$, or $AB$. For all $n$, the matrix for $[\mho^1_r]$ is the $1 \times 4$ submatrix obtained by discarding all rows and columns except those labeled $\varnothing$, $B$, $C$, or $BC$.
\end{proposition}

\begin{proof} 
This proposition follows from Corollary \ref{PartialKSDecatLemma} and the definition of $\mho^1$ in \cite[Section 9.1]{OSzNew}.
\end{proof}

At this point we can prove the following subcase of Theorems \ref{MainThmFinalForm} and \ref{MultiMainThmFinalForm}, which will be useful in dealing with general minima:
\begin{proposition}\label{MainThmSubcase}
For the triangulated category $\T' = \SSp^{\sing}(\Ccrl(n,\Sc))$ with single intrinsic grading, make the identification 
\[
K_0(\T') \otimes_{\frac{\Z[t^{\pm 1/2},q^{\pm 1}]}{q=t^{1/2}}} \C(q) \cong V^{\otimes \Sc}
\]
using the correspondence
\begin{equation}\label{BasisIdentSubThmSS}
[S_{\Ib_{\x}}] \leftrightarrow l_{\x}.
\end{equation}
Similarly, for the triangulated category $\T' = \SSp^{\multi}_Y(\Ccrl(n,\Sc))$ with multiple intrinsic gradings by $M_Y$, make the identification 
\[
K_0(\T') \cong \A^1_{P_Y}(Y)
\]
using the correspondence
\begin{equation}\label{MultiBasisIdentSubThmSS}
[S_{\Ib_{\x}}] \leftrightarrow l_{\x}.
\end{equation}
Let $\Gamma$ be an $(n,n)$--tangle consisting of a single crossing between strands $i$ and $i+1$, positive as a braid generator. Let $\widetilde{\Gamma}$ be the tangle obtained by rotating $\Gamma$ and reversing the orientations on strands as in Figure \ref{RotateReverseFig}.
\begin{itemize}
\item
Under the identification of (\ref{BasisIdentSubThmSS}), the linear map $[\Pc_{r/l}^i]$ from Section \ref{SemisimpleSect} in the singly graded case agrees with the map associated to $\widetilde{\Gamma}$ in Section \ref{MapsForCrossingsSect} (with modified-basis matrix given in Section \ref{MapsInModifiedBasisSect}). 
\item
Under the identification of (\ref{MultiBasisIdentSubThmSS}), the linear map $[\Pc_{r/l}^i]$ from Section \ref{SemisimpleSect} in the multiply graded case agrees with the map associated to $\widetilde{\Gamma}$ in Section \ref{ViroMapsForCrossingsSect} (with modified-basis matrix given in Section~\ref{ViroMapsModifiedBasisSect}). 
\end{itemize}
Now let $\Gamma$ be an $(n,n+2)$--tangle consisting of a single minimum point between strands $1$ and $2$. Let $\widetilde{\Gamma}$ be defined as usual.
\begin{itemize}
\item
Under the identification of (\ref{BasisIdentSubThmSS}), the linear map $[\mho^1_{r/l}]$ from Section \ref{SemisimpleSect} in the singly graded case agrees with the map associated to $\widetilde{\Gamma}$ in Section \ref{MinimaMaximaSect} (with modified-basis matrix given in Section \ref{MapsInModifiedBasisSect}). 
\item
Under the identification of (\ref{MultiBasisIdentSubThmSS}), the linear map $[\mho^1_{r/l}]$ from Section \ref{SemisimpleSect} in the multiply graded case is equal to $t^{-2}$ times the map associated to $\widetilde{\Gamma}$ in Section \ref{ViroMinMaxSect} (with modified-basis matrix given in Section~\ref{ViroMapsModifiedBasisSect}). 
\end{itemize}
\end{proposition}

\begin{proof}
First we consider single gradings. In the generic case $1 < i < n-1$, the result for crossings follows by comparing the four matrices in Proposition \ref{DualCanonicalRForm} with the four matrices in Proposition \ref{DecatCorr}, making the substitution $q = t^{1/2}$. The geometric correspondence $\Gamma \leftrightarrow \widetilde{\Gamma}$ means that the second matrix in one of the propositions should be compared with the third matrix in the other proposition, and vice-versa. The special cases $i = 1$ and $i = n-1$ follow by comparing submatrices. The result for $\mho^1_l$ follows by comparing the second matrix in Proposition \ref{MinMaxModifiedBasis} with the matrix in Proposition \ref{Mho1Decat}; the result for $\mho^1_r$ follows by comparing $1 \times 4$ submatrices of these $2 \times 8$ matrices.

Now we consider multiple gradings. In the generic case $1 < i < n-1$, the result for crossings follows by comparing the four matrices in Proposition \ref{ViroCrossingMatrices} with the four matrices in Proposition \ref{MultiDecatComputation}, making the substitutions $\sigma_i^{1/4} = t_i$. The geometric correspondence $\Gamma \leftrightarrow \widetilde{\Gamma}$ again means that the second matrix in one of the propositions should be compared with the third matrix in the other proposition, and vice-versa. It also implies that the labeling of strands as $i$ and $i+1$ is swapped in $\widetilde{\Gamma}$ when compared to $\Gamma$, so in one of the two propositions, one should interchange $i$ and $i+1$ everywhere. The special cases $i = 1$ and $i = n-1$ follow by comparing submatrices. The result for $\mho^1_l$ follows by comparing the second matrix in Proposition \ref{ViroMinMaxProp} (note the factor of $t^2$) with the matrix in Proposition \ref{Mho1Decat}. The result for $\mho^1_r$ again follows by comparing submatrices.
\end{proof}

\begin{proposition}\label{MhoDecatProp}
Let $1 \leq i \leq n-1$. The matrix for $[\mho^{i+1}_{r/l}]$ in terms of local idempotents is
\[
\kbordermatrix{
& \varnothing & A & B & C & AB & AC & BC & ABC \\
\varnothing & 0 & 1 & 0 & 1 & 0 & 0 & 0 & 0 \\
A & 0 & 0 & 0 & 0 & 0 & 1 & 0 & 0 \\
},
\]
in both the singly graded and the multiply graded cases.
\end{proposition}

\begin{proof}
In \cite[Section 9.3]{OSzNew}, Ozsv{\'a}th and Szab{\'o} define $\mho^{i+1}$ for general $i$ as a tensor product of $\mho^1$ with some crossing bimodules. Thus, $[\mho^{i+1}_{r/l}]$ is a weighted sum over partial Kauffman states of a tangle diagram $\Gamma$ such as the one in \cite[Figure 29]{OSzNew}, by Corollary \ref{PartialKSDecatLemma}; see also Figure \ref{GeneralMinFig} above.

By Proposition \ref{MainThmSubcase}, $[\mho^{i+1}_{r/l}]$ in the singly graded case is equal to the $\Uq$--linear map associated to $\widetilde{\Gamma}$. The tangle $\widetilde{\Gamma}$ is isotopic to a crossingless tangle $\widetilde{\Gamma'}$ with a single maximum point between strands $i+1$ and $i+2$. By invariance of the $\Uq$--linear map under Reidemeister moves (only R2 moves are needed here), $\widetilde{\Gamma}$ and $\widetilde{\Gamma'}$ induce the same $\Uq$--linear map. Thus, $[\mho^{i+1}_{r/l}]$ agrees with the $\Uq$--linear map associated to $\widetilde{\Gamma'}$, which is given by the second matrix in Proposition \ref{MinMaxModifiedBasis}.

In the multiply graded case, the same argument works except that we need to keep track of the scalar factors: $[\mho^{i+1}_{r/l}]$ is $t^{-2}$ times Viro's map for $\widetilde{\Gamma}$, which equals $t^{-2}$ times Viro's map for $\widetilde{\Gamma'}$ by the Reidemeister invariance of Viro's maps for tangles (again, only R2 moves are needed). Thus, $[\mho^{i+1}_{r/l}]$ is represented by $t^{-2}$ times the second matrix in Proposition \ref{ViroMinMaxProp}.
\end{proof}

When $i = n$, incoming idempotents $\Ib_{\y}$ of $\Cc_l(n+2,\Sc)$ must have local form $\varnothing$, $A$, $B$, or $AB$. Outgoing idempotents $\Ib_{\x}$ of $\Cc_l(n,\Sc')$ must have local form $\varnothing$.

\begin{corollary} 
The matrix for $[\mho^{n+1}_l]$ is a $1 \times 4$ submatrix of the matrix from Proposition \ref{MhoDecatProp}. Only entries whose row-label and column label are $\varnothing$, $A$, $B$, or $AB$ are included in the submatrix.
\end{corollary}

\subsection{The terminal Type A structure}\label{TTASect}
We will not attempt to fit $\widehat{t\mho}$ into the decategorification schemes above. We simply note that the Kauffman-state formulas we have been using, applied to $\widehat{t\mho}$, give the following matrices:
\begin{itemize}

\item $\widehat{t\mho}_{\Ccr(2,\{2\})}$ (for orientations \MinLR) is assigned the matrix
\[
\kbordermatrix{
& \varnothing & B & C & BC  \\
\varnothing & 0 & 1 & 0 & 0
}
\]

\item $\widehat{t\mho}_{\Ccl(2,\{2\})}$ (for orientations \MinLR) is assigned the matrix
\[
\kbordermatrix{
& \varnothing & A & B & AB  \\
\varnothing & 0 & 1 & 1 & 0
}
\]

\item $\widehat{t\mho}_{\Ccr(2,\{1\})}$ (for orientations \MinRL) is assigned the matrix
\[
\kbordermatrix{
& \varnothing & B & C & BC  \\
\varnothing & 0 & 1 & 1 & 0
}
\]

\item $\widehat{t\mho}_{\Ccl(2,\{1\})}$ (for orientations \MinRL) is assigned the matrix
\[
\kbordermatrix{
& \varnothing & A & B & AB  \\
\varnothing & 0 & 0 & 1 & 0
}
\]
\end{itemize}
We do not know whether these matrices have analogues in the representation theory of $\Uq$. 

\subsection{Comparing with representation theory}\label{RepThyCompareSect}

\begin{figure} \centering
\includegraphics[scale=0.625]{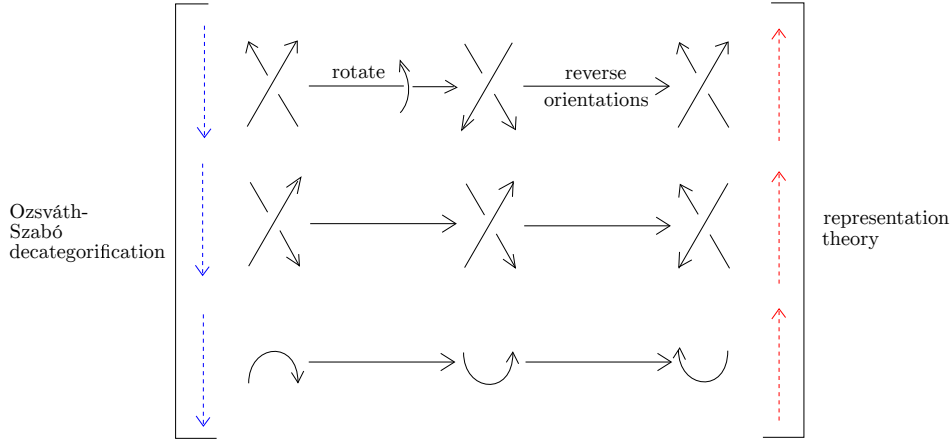}
\caption{Changing tangles $\Gamma$ to $\widetilde{\Gamma}$.}
\label{RotateReverseFig}
\end{figure}

Finally, we state the remaining theorems from the introduction precisely and give proofs. The geometric identification of tangles $\Gamma \leftrightarrow \widetilde{\Gamma}$ we use is shown in Figure \ref{RotateReverseFig}.

We start with the singly-graded case.

\begin{theorem}[cf. Theorem~\ref{IntroMainThm}]\label{MainThmFinalForm} 
For the triangulated category 
\[
\T = \D^b(\Ccrl(n,\Sc))
\]
with single intrinsic grading, make the identification
\[
K_0(\T) \otimes_{\frac{\Z[t^{\pm 1/2},q^{\pm 1}]}{q=t^{1/2}}} \C(q) \cong V^{\otimes \Sc}
\]
using the correspondence
\begin{equation}\label{BasisIdentFinalThm}
[Q_{\Ib_{\x}}] \leftrightarrow l_{\x},
\end{equation}
where $\{l_{\x}\}$ is the right (for $\Ccr(n,\Sc)$) or left (for $\Ccl(n,\Sc)$) modified basis for $V^{\otimes \Sc}$ from Definition~\ref{ModifiedBasisDef}.

Similarly, for the triangulated category $\T' = \SSp^{\sing}(\Ccrl(n,\Sc))$ with single intrinsic grading, make the identification 
\[
K_0(\T') \otimes_{\frac{\Z[t^{\pm 1/2},q^{\pm 1}]}{q=t^{1/2}}} \C(q) \cong V^{\otimes \Sc}
\]
using the correspondence
\begin{equation}\label{BasisIdentFinalThmSS}
[S_{\Ib_{\x}}] \leftrightarrow l_{\x}.
\end{equation}

In each of the following cases, for a tangle $\Gamma$, let $\widetilde{\Gamma}$ be the tangle obtained by rotating $\Gamma$ and reversing the orientations on strands as in Figure \ref{RotateReverseFig}.

\begin{itemize}
\item
Let $\Gamma$ be an $(n,n)$--tangle consisting of a single crossing between strands $i$ and $i+1$. Under the identification of (\ref{BasisIdentFinalThm}), the linear map $[\Pc^i_{r/l} \boxtimes -]$ or $[\Nc^i_{r/l} \boxtimes -]$ (as appropriate) from Section \ref{BddTypeDSect} agrees with the map associated to $\widetilde{\Gamma}$ in Section \ref{MapsForCrossingsSect} (with modified-basis matrix given in Section~\ref{MapsInModifiedBasisSect}). The same is true for the maps $[\Pc^i_{r/l}]$ and $[\Nc^i_{r/l}]$ from Section \ref{SemisimpleSect}, under the identification of (\ref{BasisIdentFinalThmSS}).

\item
Let $\Gamma$ be an $(n+2,n)$--tangle with no crossings and strands $i+1$ and $i+2$ matched on the bottom. Under the identification of (\ref{BasisIdentFinalThmSS}), the linear map $[\Omega^{i+1}_{r/l}]$ agrees with the map associated to $\widetilde{\Gamma}$ in Section \ref{MinimaMaximaSect} (with modified-basis matrix given in Section~\ref{MapsInModifiedBasisSect}).

\item
Let $\Gamma$ be an $(n,n+2)$--tangle with no crossings and strands $i+1$ and $i+2$ matched on top. Under the identification of (\ref{BasisIdentFinalThmSS}), the linear map $[\mho^{i+1}_{r/l}]$ agrees with the map associated to $\widetilde{\Gamma}$ in Section \ref{MinimaMaximaSect} (with modified-basis matrix given in Section~\ref{MapsInModifiedBasisSect}).

\end{itemize}
\end{theorem}

\begin{proof}
The first item, for positive braid generators, follows from Proposition \ref{MainThmSubcase}; recall that $[\Pc^i_{r/l} \boxtimes -]$ has the same matrix as $[\Pc^i_{r/l}]$. The result for negative braid generators follows by taking inverses of matrices. 

The second and third items, in the generic case $0 < i < n$, follow from comparing the first and second matrices in Proposition \ref{MinMaxModifiedBasis} with the matrices in Proposition \ref{OmegaDecatProp} and Proposition \ref{MhoDecatProp} respectively. These comparisons are compatible with $\Gamma \leftrightarrow \widetilde{\Gamma}$. The special cases $i = 0$ and $i = n$ follow by comparing submatrices.
\end{proof}

We have an analogous theorem in the multi-graded setting:

\begin{theorem}[cf. Theorem~\ref{MultiIntroMainThm}]\label{MultiMainThmFinalForm} 
Let $Y$ be a zero-manifold consisting of $n$ points oriented by $\Sc$. For the triangulated category 
\[
\T = \D^b_{M_Y}(\Ccrl(n,\Sc))
\]
with multiple intrinsic gradings by $M_Y := H^0(Y; \frac{1}{4}\Z)$, make the identification
\[
K_0(\T) \cong \A^1_{P_Y}(Y)
\]
using the correspondence
\begin{equation}\label{MultiBasisIdentFinalThm}
[Q_{\Ib_{\x}}] \leftrightarrow l_{\x},
\end{equation}
where $\{l_{\x}\}$ is the right (for $\Ccr(n,\Sc)$) or left (for $\Ccl(n,\Sc)$) modified basis for $\A^1_{P_Y}(Y)$ from Definition~\ref{ViroModifiedBasisDef}.

Similarly, for the triangulated category $\T' = \SSp^{\multi}_Y(\Ccrl(n,\Sc))$ with multiple intrinsic gradings by $M_Y$, make the identification 
\[
K_0(\T') \cong \A^1_{P_Y}(Y)
\]
using the correspondence
\begin{equation}\label{MultiBasisIdentFinalThmSS}
[S_{\Ib_{\x}}] \leftrightarrow l_{\x}.
\end{equation}

In each of the following cases, for a tangle $\Gamma$, let $\widetilde{\Gamma}$ be the tangle obtained by rotating $\Gamma$ and reversing the orientations on strands as in Figure \ref{RotateReverseFig}.

\begin{itemize}
\item
Let $\Gamma$ be an $(n,n)$--tangle consisting of a single crossing between strands $i$ and $i+1$. Under the identification of (\ref{MultiBasisIdentFinalThm}), the linear map $[\Pc^i_{r/l} \boxtimes -]$ or $[\Nc^i_{r/l} \boxtimes -]$ (as appropriate) from Section \ref{BddTypeDSect} agrees with the map associated to $\widetilde{\Gamma}$ in Section \ref{ViroMapsForCrossingsSect} (with modified-basis matrix given in Section~\ref{ViroMapsModifiedBasisSect}). The same is true for the maps $[\Pc^i_{r/l}]$ and $[\Nc^i_{r/l}]$ from Section \ref{SemisimpleSect}, under the identification of (\ref{MultiBasisIdentFinalThmSS}).

\item
Let $\Gamma$ be an $(n+2,n)$--tangle with no crossings and strands $i+1$ and $i+2$ matched on the bottom. Under the identification of (\ref{MultiBasisIdentFinalThmSS}), the linear map $[\Omega^{i+1}_{r/l}]$ agrees with the map associated to $\widetilde{\Gamma}$ in Section \ref{ViroMinMaxSect} (with modified-basis matrix given in Section~\ref{ViroMapsModifiedBasisSect}), up to a scalar factor.

\item
Let $\Gamma$ be an $(n,n+2)$--tangle with no crossings and strands $i+1$ and $i+2$ matched on top. Under the identification of (\ref{MultiBasisIdentFinalThmSS}), the linear map $[\mho^{i+1}_{r/l}]$ agrees with the map associated to $\widetilde{\Gamma}$ in Section \ref{ViroMinMaxSect} (with modified-basis matrix given in Section~\ref{ViroMapsModifiedBasisSect}), up to a scalar factor.

\end{itemize}
\end{theorem}

\begin{proof}
As before, the first item for positive braid generators follows from Proposition \ref{MainThmSubcase}. The result for negative braid generators follows by taking inverses of matrices. 

The second and third items follow from comparing the first and second matrices in Proposition \ref{ViroMinMaxProp} with the matrices in Proposition \ref{OmegaDecatProp} and Proposition \ref{MhoDecatProp} respectively; note the scalar factors of $t^{-2}$ and $t^{2}$ in Proposition \ref{ViroMinMaxProp}. These comparisons are compatible with $\Gamma \leftrightarrow \widetilde{\Gamma}$. The special cases $i = 0$ and $i = n$ follow by comparing submatrices.
\end{proof}

\begin{corollary}\label{InvarianceCorr}
The map $[X]$ of Definition \ref{InducedDAMapSemisimple}, when $X$ is Ozsv{\'a}th--Szab{\'o}'s DA bimodule for a tangle diagram $\Gamma$ with no closed components, is an invariant of the underlying tangle of $\Gamma$. 
\end{corollary}

\begin{proof}
This corollary follows from Theorem \ref{MainThmFinalForm} in the singly graded case and Theorem \ref{MultiMainThmFinalForm} in the multiply graded case, together with the Reidemeister invariance of $\Uq$--linear maps induced by tangles colored with $V$ and $V^*$, on which the ribbon element of $\Uq$ acts trivially (see \cite[Theorem 4.7]{Ohtsuki}) and the Reidemeister invariance of Viro's maps induced by tangles whose coloring has $T$--component zero (see \cite[Section 2.9]{Viro}).
\end{proof}

\bibliographystyle{alpha}
\bibliography{biblio}

\end{document}